\documentclass[a4paper, 11pt]{article}
\pdfoutput=1 

\usepackage[T1]{fontenc} 
\usepackage[utf8]{inputenc} 

\usepackage{amsfonts,amsmath,amssymb,mathtools,dsfont} 
\usepackage{microtype} 
\usepackage[usenames, dvipsnames]{xcolor} 

\usepackage{fancyhdr}

\usepackage[labelsep = period, labelfont = bf, justification = centering]{caption}
\usepackage{float,graphicx,subcaption} 

\usepackage{booktabs,makecell} 
\usepackage{enumitem,mdwlist} 
\setlist[itemize]{topsep=0ex,itemsep=0ex,parsep=0.4ex}
\setlist[enumerate]{topsep=0ex,itemsep=0ex,parsep=0.4ex}

\usepackage{etoolbox}
\usepackage[super]{nth} 
\allowdisplaybreaks

\usepackage[left = 2cm, right = 2cm, top = 2.5cm, bottom = 2.5cm, headsep = 12pt, headheight = 15pt]{geometry}
\usepackage{parskip}

\usepackage[noBBpl]{mathpazo} 
\DeclareFontFamily{U}{matha}{\hyphenchar\font45}
\DeclareFontShape{U}{matha}{m}{n}{
	<5> <6> <7> <8> <9> <10> gen * matha
	<10.95> matha10 <12> <14.4> <17.28> <20.74> <24.88> matha12
}{}
\DeclareSymbolFont{matha}{U}{matha}{m}{n}

\DeclareMathSymbol{\specialuparrow}{\mathrel}{matha}{"D2}
\DeclareMathSymbol{\specialrightarrow}{\mathrel}{matha}{"D1}

\DeclareFontFamily{U} {cmmi}{}

\DeclareFontShape{U}{cmmi}{m}{n}{
	<-6> cmmi5
	<6-7> cmmi6
	<7-8> cmmi7
	<8-9> cmmi8
	<9-10> cmmi9
	<10-12> cmmi10
	<12-> cmmi12}{}

\DeclareSymbolFont{Xcmmi} {U} {cmmi}{m}{n}

\DeclareMathSymbol{\mu}{\mathord}{Xcmmi}{'026}
\DeclareMathSymbol{\rho}{\mathord}{Xcmmi}{'032}
\DeclareMathSymbol{\varphi}{\mathord}{Xcmmi}{'047}

\DeclareFontFamily{U} {cmr}{}

\DeclareFontShape{U}{cmr}{m}{n}{
	<-6> cmr5
	<6-7> cmr6
	<7-8> cmr7
	<8-9> cmr8
	<9-10> cmr9
	<10-12> cmr10
	<12-> cmr12}{}

\DeclareSymbolFont{Xcmr} {U} {cmr}{m}{n}
\DeclareMathSymbol{\Delta}{\mathord}{Xcmr}{'001}
\DeclareMathSymbol{\Upsilon}{\mathord}{Xcmr}{'007}
\DeclareMathSymbol{\Omega}{\mathord}{Xcmr}{'012}

\usepackage[hyphens]{url} 
\usepackage[linktoc = all, hidelinks, colorlinks, unicode=true, pagebackref=true]{hyperref} 
\usepackage[capitalise, compress, nameinlink, noabbrev]{cleveref} 

\hypersetup{linkcolor={blue!70!black}, citecolor={green!70!black}, urlcolor={blue!70!black}}

\crefname{section}{\S}{\S\S} 
\Crefname{section}{Section}{Sections} 
\crefname{subsection}{\S}{\S\S} 
\Crefname{subsection}{Subsection}{Subsections} 
\crefname{page}{page}{pages}
\Crefname{page}{Page}{Pages}

\usepackage{pgf,tikz}
\usepackage{mathrsfs}
\usetikzlibrary{decorations.pathreplacing}
\usetikzlibrary[patterns]
\usetikzlibrary{arrows.meta}
\usepackage{tkz-euclide}

\tikzset{font={\fontsize{10pt}{12}\selectfont}}
\tkzSetUpPoint[fill=black, size = 3pt]
\tkzSetUpLine[color=black, line width=0.6pt]

\usepackage{amsthm,thmtools,thm-restate}

\declaretheoremstyle[
spaceabove = .5\baselineskip\@plus.2\baselineskip\@minus.2\baselineskip, 
spacebelow = .2\baselineskip\@plus.2\baselineskip\@minus.2\baselineskip,
headfont = \normalfont\bfseries,
notefont = \mdseries, 
notebraces = {(}{)},
bodyfont = \normalfont\itshape,
postheadspace = .5em,
headpunct = .
]{bolditalic}

\declaretheoremstyle[
spaceabove = .5\baselineskip\@plus.2\baselineskip\@minus.2\baselineskip, 
spacebelow = .2\baselineskip\@plus.2\baselineskip\@minus.2\baselineskip,
headfont = \normalfont\bfseries,
notefont = \mdseries, 
notebraces = {(}{)},
bodyfont = \normalfont,
postheadspace = .5em,
headpunct = .
]{boldnormal}

\declaretheoremstyle[
spaceabove = .2\baselineskip\@plus.2\baselineskip\@minus.2\baselineskip, 
spacebelow = .5\baselineskip\@plus.2\baselineskip\@minus.2\baselineskip,
headfont = \normalfont\itshape,
notefont = \mdseries, 
notebraces = {}{},
bodyfont = \normalfont,
postheadspace = .5em,
headpunct = .,
qed = \qedsymbol
]{proofstyle}

\declaretheoremstyle[
spaceabove = .5\baselineskip\@plus.2\baselineskip\@minus.2\baselineskip, 
spacebelow = .2\baselineskip\@plus.2\baselineskip\@minus.2\baselineskip,
headfont = \normalfont\bfseries,
notefont = \mdseries, 
notebraces = {(}{)},
bodyfont = \normalfont,
postheadspace = .5em,
headpunct = .,
qed = \qedsymbol
]{solutionstyle}

\renewcommand*{\backref}[1]{}
\renewcommand*{\backrefalt}[4]{
	\ifcase #1 Not cited.%
	\or $\specialuparrow$#2%
	\else $\specialuparrow$#2%
	\fi%
}

\declaretheorem[name = Conjecture, numberwithin = section, style = bolditalic, refname = {Conjecture,Conjectures}, Refname = {Conjecture,Conjectures}]{conjecture}
\declaretheorem[name = Corollary, numberlike = conjecture, style = bolditalic, refname = {Corollary,Corollaries}, Refname = {Corollary,Corollaries}]{corollary}
\declaretheorem[name = Definition, numberlike = conjecture, style = bolditalic, refname = {Definition,Definitions}, Refname = {Definition,Definitions}]{definition}

\declaretheorem[name = Lemma, numberlike = conjecture, style = bolditalic, refname = {Lemma,Lemmas}, Refname = {Lemma,Lemmas}]{lemma}

\declaretheorem[name = Theorem, numberlike = conjecture, style = bolditalic, refname = {Theorem,Theorems}, Refname = {Theorem,Theorems}]{theorem}
\declaretheorem[name = Theorem, numbered = no, style = bolditalic, refname = {Theorem,Theorems}, Refname = {Theorem,Theorems}]{theorem*}

\declaretheorem[name = Claim, numberwithin = conjecture, style = bolditalic, refname = {Claim,Claims}, Refname = {Claim,Claims}]{claim}

\declaretheorem[name = Remark, numberlike = conjecture, style = boldnormal, refname = {Remark,Remarks}, Refname = {Remark,Remarks}]{remark}

\declaretheorem[name = Proof, numbered = no, style = proofstyle, refname = {Proof,Proofs}, Refname = {Proof,Proofs}]{Proof}

\renewcommand{\epsilon}{\varepsilon}

\DeclarePairedDelimiter{\abs}{\lvert}{\rvert}

\DeclarePairedDelimiter{\set}{\{}{\}}

\newcommand*{\bR}{\mathbb{R}}

\newcommand*{\ZZN}{\mathbb{Z}_{\geqslant 0}}

\newcommand*{\cF}{\mathcal{F}}

\newcommand*{\cO}{\mathcal{O}}

\newcommand{\defn}[1]{\textcolor{Maroon}{\emph{#1}}}

\begin{document}

\author{Freddie Illingworth\footnotemark[1]}

\title{\bf The chromatic profile of locally bipartite graphs}

\date{24 May 2022}

\maketitle

\begin{abstract}
	In 1973, Erd\H{o}s and Simonovits asked whether every $n$-vertex triangle-free graph with minimum degree greater than $1/3 \cdot n$ is 3-colourable. This question initiated the study of the chromatic profile of triangle-free graphs: for each $k$, what minimum degree guarantees that a triangle-free graph is $k$-colourable. This problem has a rich history which culminated in its complete solution by Brandt and Thomass\'{e}. Much less is known about the chromatic profile of $H$-free graphs for general $H$.
	
	Triangle-free graphs are exactly those in which each neighbourhood is one-colourable. Locally bipartite graphs, first mentioned by \L uczak and Thomass\'{e}, are the natural variant of triangle-free graphs in which each neighbourhood is bipartite. Here we study the chromatic profile of locally bipartite graphs. We show that every $n$-vertex locally bipartite graph with minimum degree greater than $4/7 \cdot n$ is 3-colourable ($4/7$ is tight) and with minimum degree greater than $6/11 \cdot n$ is 4-colourable. Although the chromatic profiles of locally bipartite and triangle-free graphs bear some similarities, we will see there are striking differences.
\end{abstract}
	
\renewcommand{\thefootnote}{\fnsymbol{footnote}} 
	
	
\footnotetext[1]{Mathematical Institute, University of Oxford (\textsf{illingworth@maths.ox.ac.uk}). Research carried out while at DPMMS, University of Cambridge. Research supported by EPSRC grant 2114463.}
	
\renewcommand{\thefootnote}{\arabic{footnote}} 
	
\section{Introduction}

The Andr\'{a}sfai-Erd\H{o}s-S\'{o}s theorem is a classical result which considers the chromatic number of dense graphs (those with large minimum degree) that do not contain the clique $K_{r + 1}$ (for a short proof see Brandt~\cite{Brandt2003AES}). It can be viewed as a minimum degree analogue of Erd\H{o}s and Simonovits's stability theorem~\cite{Erdos1967, Erdos1968, Simonovits1968} for the structure of $K_{r + 1}$-free graphs with close to the maximum number of edges.

\begin{theorem}[Andr\'{a}sfai-Erd\H{o}s-S\'{o}s,~\cite{AES1974}]\label{AES}
	Let $r \geqslant 2$ and $G$ be a $K_{r + 1}$-free graph with $n$ vertices and minimum degree greater than
	\begin{equation*}
		\biggl(1 - \frac{1}{r - 1/3}\biggr) \cdot n.
	\end{equation*}
	Then $G$ is $r$-colourable. Furthermore $1 - 1/(r - 1/3)$ is tight.
\end{theorem}

At around the same time, Erd\H{o}s and Simonovits~\cite{ErdosSimonovits1973} implicitly posed a very general problem: for a fixed graph $H$ and positive integer $k$, compute
\begin{equation*}
	\delta_{\chi}(H, k) = \inf\set{d \colon \text{if }\delta(G) \geqslant d \abs{G} \text{ and } G \text{ is } H\text{-free}, \text{then } \chi(G) \leqslant k}.
\end{equation*}
\Cref{AES} says exactly that $\delta_{\chi}(K_{r + 1}, r) = 1 - 1/(r - 1/3)$. The values of $\delta_{\chi}(H, k)$ as $k$ varies form the \defn{chromatic profile} of the family of $H$-free graphs. Erd\H{o}s and Simonovits stated that determining the chromatic profile for general $H$ was `too complicated' and this sentiment was reiterated recently by Allen, B\"{o}ttcher, Griffiths, Kohayakawa and Morris~\cite{ABGKM2013} -- indeed, for many $H$, even the chromatic number of the $H$-free graph with the most edges is unknown. The chromatic profile of triangle-free graphs has been extensively studied~\cite{AES1974,Brandt1999,BrandtThomasse2005,CJK1997,Haggkvist1982,Jin1995,Luczak2006,Thomassen2002} and is now known (and has been extended to $K_{r + 1}$-free graphs~\cite{GoddardLyle2010,Nikiforov2010}). \Cref{spec4triangle} summarises this concisely.
\begin{table}[H]
	\centering
	\begin{tabular}{r|cccc}
		\toprule
		$\delta(G)\abs{G}^{-1} >$ & $2/5$ & $10/29$ & $1/3$ & $1/3 - \varepsilon$ \\
		\midrule
		$\chi(G) \leqslant $ & 2 & 3 & 4 & $\infty$ \\
		\bottomrule
	\end{tabular}
	\caption{Chromatic profile of a triangle-free graph $G$.}\label{spec4triangle}
\end{table}
Clearly the sequence $\delta_{\chi}(H, k)$ is decreasing and so must tend to a limit, called the \defn{chromatic threshold} of $H$-free graphs. \Cref{spec4triangle} exhibits an interesting phenomenon: $\delta_{\chi}(K_{3}, k) = 1/3$ for all $k \geqslant 4$ and so the chromatic threshold of triangle-free graphs is $1/3$. More generally, the chromatic threshold of $H$-free graphs is defined as
\begin{align*}
	\delta_{\chi}(H) & = \inf_{k} \delta_{\chi}(H, k) \\
	& = \inf \set{d \colon \exists C = C(H, d) \text{ such that if } \delta(G) \geqslant d \abs{G} \text{ and } G \text{ is } H \text{-free}, \text{then } \chi(G) \leqslant C}.
\end{align*}
Much more is known about the chromatic threshold than the chromatic profile. Indeed, in~\cite{ABGKM2013}, Allen, B\"{o}ttcher, Griffiths, Kohayakawa, and Morris determined the chromatic threshold of $H$-free graphs for every graph $H$. They noted that there are other natural families of graphs which are beyond the reach of their techniques. For a family of graphs, $\cF$, one defines the chromatic profile and threshold analogously.
\begin{align*}
	\delta_{\chi}(\cF, k) & = \inf\set{d \colon \text{if }\delta(G) \geqslant d \abs{G} \text{ and } G \in \cF, \text{then } \chi(G) \leqslant k}, \\
	\delta_{\chi}(\cF) & = \inf_{k} \delta_{\chi}(\cF, k) \\
	& = \inf \set{d \colon \exists C = C(\cF, d) \text{ such that if } \delta(G) \geqslant d \abs{G} \text{ and } G \in \cF, \text{then } \chi(G) \leqslant C}.
\end{align*}
In~\cite{Illingworth2023mindegstab} we extend the Andr\'{a}sfai-Erd\H{o}s-S\'{o}s theorem to non-complete graphs. A certain family of graphs plays a central role. This family is a natural extension of triangle-free graphs and has previously appeared in the literature. Triangle-free graphs are exactly those in which each neighbourhood is independent (i.e.\ 1-colourable). Graphs in which each neighbourhood is 2-colourable (or 3-colourable \ldots) are termed \defn{locally bipartite} (\defn{locally tripartite} \ldots). Locally bipartite graphs were first mentioned a decade ago by \L uczak and Thomass\'{e}~\cite{LuczakThomasse2010} who asked for their chromatic threshold, conjecturing it was $1/2$. This was confirmed by Allen, B\"{o}ttcher, Griffiths, Kohayakawa, and Morris~\cite{ABGKM2013}. However, in contrast to the well-studied case of triangle-free graphs, the chromatic profile of locally bipartite graphs, and more generally that of locally $b$-partite graphs ($b \geqslant 3$), has not previously been examined. In this paper we focus on the chromatic profile of locally bipartite graphs, deferring the chromatic threshold and profile of locally $b$-partite graphs to~\cite{Illingworth2022localbpart}.

Locally bipartite graphs, just like triangle-free ones, exhibit a spectrum of thresholds. There are, however, some interesting differences which we explore in \cref{sec:trianglefree}. Our understanding of their profile is summarised in our main result, \cref{spec4localbip}. The graphs $\overline{C}_{7}$ and $H_{2}^{+}$ can be seen in \cref{fig:graphs} where they are discussed more thoroughly -- for now it suffices to note that they are both small 4-chromatic locally bipartite graphs which are edge-maximal with respect to local bipartiteness.
\begin{theorem}[locally bipartite graphs]\label{spec4localbip}
	Let $G$ be a locally bipartite graph.
	\begin{itemize}[noitemsep]
		\item If $\delta(G) > 4/7 \cdot \abs{G}$, then $G$ is $3$-colourable.
		\item If $\delta(G) > 5/9 \cdot \abs{G}$, then there is a homomorphism $G \to \overline{C}_{7}$.
		\item There is an absolute constant $\varepsilon > 0$ such that if $\delta(G) > (5/9 - \varepsilon) \cdot \abs{G}$, then there is either a homomorphism $G \to \overline{C}_{7}$ or $G \to H_{2}^{+}$.
		\item If $\delta(G) > 6/11 \cdot \abs{G}$, then $G$ is $4$-colourable.
	\end{itemize}
	Furthermore $4/7$ and $5/9$ are tight, as demonstrated by balanced blow-ups of $\overline{C}_{7}$ and suitable blow-ups of $H_{2}^{+}$ \textnormal{(}see \cref{fig:weightings} on \cpageref{fig:weightings}\textnormal{)}, respectively.
\end{theorem}
The first bullet point corresponds to the $r = 2$ case of \cref{AES}: every $n$-vertex triangle-free graph with minimum degree greater than $2/5 \cdot n$ is bipartite. That result has a very short proof, while the first bullet point requires a more substantial argument whose sketch precedes \cref{sec:H02C7}. \Cref{spec4localbip} is a summary of our understanding, see the start of \cref{sec:localbip} for further details.

\Cref{spec4localbip} gives the following information about the profile of the family of locally bipartite graphs, which we denote by $\cF_{1, 2}$:
\begin{equation*}
	\delta_{\chi}(\cF_{1, 2}, 3) = 4/7, \quad \delta_{\chi}(\cF_{1, 2}, 4) \leqslant 6/11.
\end{equation*}
As mentioned above, Allen, B\"{o}ttcher, Griffiths, Kohayakawa, and Morris~\cite{ABGKM2013} and \L uczak and Thomass\'{e}~\cite{LuczakThomasse2010} showed that $\delta_{\chi}(\cF_{1, 2}) = 1/2$. In~\cite{Illingworth2022localbpart}, we will extend \cref{spec4localbip} below $6/11$ to give more structural (but not colourability) results for locally bipartite graphs (essential for our analysis of locally $b$-partite graphs and extension of \cref{AES}).

\subsection{The graphs}\label{sec:graphs}

Throughout the paper the following graphs will appear frequently and here we note a few of their properties to acquaint the reader.

\begin{figure}[H]
	\centering
	\begin{subfigure}{.2\textwidth}
		\centering
		\begin{tikzpicture}
			\foreach \pt in {0,1,...,6} 
			{
				\tkzDefPoint(\pt*360/7 + 90:1){v_\pt}
			} 
			\tkzDrawPolySeg(v_0,v_1,v_2,v_3,v_4,v_5,v_6,v_0) 
			\tkzDrawPolySeg(v_5,v_0,v_2)
			\tkzDrawPolySeg(v_1,v_3)
			\tkzDrawPolySeg(v_4,v_6)
			\tkzDrawPoints(v_0,v_...,v_6)
		\end{tikzpicture}
		\subcaption*{$H_{0}$}
	\end{subfigure}
	\begin{subfigure}{.2\textwidth}
		\centering
		\begin{tikzpicture}
			\foreach \pt in {0,1,...,6} 
			{
				\tkzDefPoint(\pt*360/7 + 90:1){v_\pt}
			} 
			\tkzDrawPolySeg(v_0,v_1,v_2,v_3,v_4,v_5,v_6,v_0) 
			\tkzDrawPolySeg(v_3,v_5,v_0,v_2,v_4)
			\tkzDrawPolySeg(v_6,v_1)
			\tkzDrawPoints(v_0,v_...,v_6)
		\end{tikzpicture}
		\subcaption*{$H_{1}$}
	\end{subfigure}
	\begin{subfigure}{.2\textwidth}
		\centering
		\begin{tikzpicture}
			\foreach \pt in {0,1,...,6} 
			{
				\tkzDefPoint(\pt*360/7 + 90:1){v_\pt}
			} 
			\tkzDrawPolySeg(v_0,v_1,v_2,v_3,v_4,v_5,v_6,v_0)
			\tkzDrawPolySeg(v_1,v_3,v_5,v_0, v_2,v_4,v_6)
			\tkzDrawPoints(v_0,v_...,v_6)
		\end{tikzpicture}
		\subcaption*{$H_{2}$}
	\end{subfigure}
	\begin{subfigure}{.2\textwidth}
		\centering
		\begin{tikzpicture}
			\foreach \pt in {0,1,...,6} 
			{
				\tkzDefPoint(\pt*360/7 + 90:1){v_\pt}
			} 
			\tkzDrawPolySeg(v_0,v_1,v_2,v_3,v_4,v_5,v_6,v_0) 
			\tkzDrawPolySeg(v_0,v_2,v_4,v_6,v_1,v_3,v_5,v_0)
			\tkzDrawPoints(v_0,v_...,v_6)
		\end{tikzpicture}
		\subcaption*{$C_{7}^{2} = \overline{C}_{7}$}
	\end{subfigure}
	
	\medskip
	
	\begin{subfigure}{.2\textwidth}
		\centering
		\begin{tikzpicture}
			\foreach \pt in {0,1,...,6} 
			{
				\tkzDefPoint(\pt*360/7 + 90:1){v_\pt}
			} 
			\tkzDefPoint(0,0){u}
			\tkzDrawPolySeg(v_0,v_1,v_2,v_3,v_4,v_5,v_6,v_0)
			\tkzDrawPolySeg(v_1,v_3,v_5,v_0, v_2,v_4,v_6)
			\tkzDrawSegments(u,v_0 u,v_2 u,v_5)
			\tkzDrawPoints(v_0,v_...,v_6)
			\tkzDrawPoint(u)
		\end{tikzpicture}
		\subcaption*{$H_{2}^{+}$}
	\end{subfigure}
	\begin{subfigure}{.2\textwidth}
		\centering
		\begin{tikzpicture}
			\foreach \pt in {0,1,...,6} 
			{
				\tkzDefPoint(\pt*360/7 + 90:1){v_\pt}
			} 
			\tkzDefPoint(0,0){u}
			\tkzDrawPolySeg(v_0,v_1,v_2,v_3,v_4,v_5,v_6,v_0) 
			\foreach \pt in {0,1,...,6} 
			{
				\tkzDrawSegment(u,v_\pt)
			} 
			\tkzDrawPoints(v_0,v_...,v_6)
			\tkzDrawPoint(u)
		\end{tikzpicture}
		\subcaption*{$W_{7}$}
	\end{subfigure}
	\caption{The graphs appearing in \cref{spec4localbip} and its proof.}\label{fig:graphs}
\end{figure}

\begin{itemize}[noitemsep]
	\item All graphs shown are 4-chromatic and all bar $W_{7}$ are locally bipartite.
	\item The graph $H_{0}$ is isomorphic to the \defn{Moser Spindle} -- the smallest 4-chromatic unit distance graph. $H_{0}$ is also the smallest 4-chromatic locally bipartite graph and so it is natural that it should play such an integral part in many of our results. The graph $\overline{C}_{7}$ is the complement (and also the square) of the 7-cycle.
	\item Adding a single edge to $H_{0}$ while maintaining local bipartiteness can give rise to two non-isomorphic graphs, one of which is $H_{1}$. The other will appear fleetingly in \cref{sec:H02C7}. Adding a single edge to $H_{1}$ while maintaining local bipartiteness gives rise to a unique (up to isomorphism) graph -- $H_{2}$. There is only one way to add a single edge to $H_{2}$ and maintain local bipartiteness -- this gives $\overline{C}_{7}$. $H_{2}^{+}$ is $H_{2}$ with a degree 3 vertex added.
	\item $\overline{C}_{7}$ and $H_{2}^{+}$ are both edge-maximal locally bipartite graphs.
	\item $W_{7}$ is called the \defn{7-wheel}. More generally, a single vertex joined to all the vertices of a $k$-cycle is called a \defn{$k$-wheel} and is denoted by $W_{k}$. We term any edge from the central vertex to the cycle a \defn{spoke} of the wheel and any edge of the cycle a \defn{rim} of the wheel. Note that a graph is locally bipartite exactly if it does not contain any odd wheel (there is no such nice characterisation for a graph being locally tripartite, locally 4-partite, \ldots).
\end{itemize}

The following observation gives a useful link between local bipartiteness and some of these graphs. We will use it frequently when copies of $H_0$, $H_1$, $H_2$ or $\overline{C}_7$ appear.

\begin{remark}\label{rmk:5nbs}
	Any five vertices of $H_0$ contain a triangle or a 5-cycle. In particular, if $G$ is a locally bipartite graph, then any vertex has at most four neighbours in any copy of $H_0$ appearing in $G$.
\end{remark}

It will sometimes be useful to check whether there is a homomorphism from one locally bipartite graph to another. This is done with the following lemma. Two vertices of a graph are called \defn{twins} if they have the same neighbourhood (in particular, twins cannot be adjacent) and we say a graph is \defn{twin-free} if no two vertices have the same neighbourhood.

\begin{lemma}\label{homscores}
	Let $F$ be a twin-free, edge-maximal locally bipartite graph and let $G$ be a locally bipartite graph. If there is a homomorphism $F \to G$, then $F$ is an induced subgraph of $G$.
\end{lemma}

\begin{Proof}
	Let $\varphi \colon F \to G$ be a homomorphism. If $\varphi$ is injective, then $F$ is a subgraph of $G$. But $G$ is locally bipartite and $F$ is edge-maximal locally bipartite, so any copy of $F$ appearing in $G$ must be induced. If $\varphi$ is not injective, then there are distinct vertices $u, v$ of $F$ with $\varphi(u) = \varphi(v)$. The neighbourhoods of $u$ and $v$ are not the same, so we may assume there is a vertex $w$ in $F$ with $w$ adjacent to $v$ but not $u$. But then $\varphi$ is a homomorphism from $F + uw$ to $G$: $\varphi(u) = \varphi(v)$ and $\varphi(v)$ is adjacent to $\varphi(w)$.
	
	However, $F$ is edge-maximal locally bipartite and $uw$ is not an edge of $F$, so $F + uw$ is not locally bipartite. Hence $\varphi$ is a homomorphism from a graph which is not locally bipartite, $F + uw$, to a locally bipartite graph, $G$, which is absurd.
\end{Proof}

\subsection{Comparison with triangle-free graphs and open questions}\label{sec:trianglefree}

Here we explain the structural results behind the complete determination of the chromatic profile of triangle-free graphs and compare them to those in \cref{spec4localbip}.

To properly discuss the chromatic profile of triangle-free and locally bipartite graphs it is useful to define blow-ups and weighted graphs. Given a graph $G$, a \defn{blow-up} of $G$ is a graph obtained by replacing each vertex $v$ of $G$ by a non-empty independent set $I_{v}$ and each edge $uv$ by a complete bipartite graph between classes $I_{u}$ and $I_{v}$. We note in passing that a graph has the same chromatic and clique numbers as any of its blow-ups and furthermore if $H$ is blow-up of $G$, then $G$ is locally bipartite if and only if $H$ is. We say a vertex $v$ has been blown-up by $n$ if $\abs{I_{v}} = n$. It is often helpful to think of this as \defn{weighting} vertex $v$ by $n$. To be precise, a \defn{weighted graph} $(G, \omega)$ is a graph $G$ together with a weighting $\omega \colon V(G) \to \bR^{+}$ and so $(G, \omega)$ can be viewed as the blow-up of $G$ in which each vertex $v$ has been blown-up by $\omega(v)$. An unweighted graph can be viewed as a weighted graph where each vertex has weight 1. If a weighted graph has a pair of twins (see \cref{sec:graphs}), then merging those vertices and giving the new vertex the sum of their weights produces an equivalent graph with the same total weight.

Suppose we start with a triangle-free graph $G$. We can repeatedly add edges to $G$ and merge twins to obtain a twin-free, edge-maximal triangle-free weighted graph $(H, \omega)$ whose total weight, $\omega(H)$, equals $\abs{G}$ and whose minimum degree, $\delta(H, \omega)$, is at least $\delta(G)$ (the degree of a vertex in $(H, \omega)$ is the total weight of its neighbours). Note that there is a homomorphism $G \to H$. In particular, to understand the chromatic profile of triangle-free graphs, one only needs to understand the twin-free, edge-maximal triangle-free graphs $H$ which have a weighting $\omega$ with $\delta(H, \omega) > 1/3 \cdot \omega(H)$ (we will refer to this last property as \defn{$H$ beating $1/3$}). The above reasoning holds if we replace ``triangle-free'' by ``locally bipartite'' and replace $1/3$ by $1/2$ (the corresponding chromatic threshold). Hence, we are particularly interested in twin-free, edge-maximal locally bipartite graphs which beat $1/2$. \Cref{homscores} applies and so for two such graph $G$ and $H$, there is a homomorphism $G \to H$ if and only if $G$ is an induced subgraph of $H$.

In the triangle-free case, the endeavour of finding all such graphs was implicitly started by Haggkvist~\cite{Haggkvist1982}, continued by Chen, Jin and Koh~\cite{CJK1997}, and finished by Brandt and Thomass\'{e}~\cite{BrandtThomasse2005}: there are two important sequences of triangle-free graphs, the 3-colourable Andr\'{a}sfai graphs~\cite{Andrasfai1962} ($\Gamma_{1} = K_{2}, \Gamma_{2} = C_{5}, \Gamma_{3}, \dotsc$) and the 4-chromatic Vega graphs~\cite{BrandtPisanski1998} (which we denote by $\Upsilon_{j}$). Brandt and Thomass\'{e} showed that the twin-free, edge-maximal triangle-free graphs beating $1/3$ are exactly the Andr\'{a}sfai and Vega graphs and so every triangle-free graph with $\delta(G) > 1/3 \cdot \abs{G}$ has a homomorphism to one of these (and hence is 4-colourable). Furthermore, for any $c > 1/3$, only finitely many graphs of each sequence beat $c$. In particular, if a triangle-free $G$ has $\delta(G)/\abs{G} \geqslant c$ for $c > 1/3$, then there is a homomorphism from $G$ to some early $\Gamma_{i}$ or to some early $\Upsilon_{j}$. 

\Cref{spec4localbip} effectively shows that $K_{3}, \overline{C}_{7}, H_{2}^{+}$ play the same role for locally bipartite graphs as the first three Andr\'{a}sfai graphs do for triangle-free graphs. Furthermore, \cref{spec4localbip}, together with \cref{homscores}, shows that they are the only twin-free, edge-maximal locally bipartite graphs which beat $5/9 - \varepsilon$ (in fact, we believe this is true down to $6/11$). These results display similarities with the triangle-free case but also give a couple of striking differences. Firstly, the Andr\'{a}sfai graphs are nested, while $H_{2}^{+}$ does not contain $\overline{C}_{7}$ and so, by \cref{homscores}, there is not even a homomorphism from one to the other. Secondly, the Andr\'{a}sfai and Vega graphs have weightings in which all vertices have the same degree as is expected for extremal examples. However, $H_{2}^{+}$ has no such weighting and, in fact, its $n$-vertex weighting with greatest minimum degree ($5/9 \cdot n$) has $1/9 \cdot n$ vertices with degree $2/3 \cdot n$ (this is shown in \cref{fig:weightings} on \cpageref{fig:weightings}).

It is natural to ask what graphs come after $H_{2}^{+}$. There is an infinite nested sequence of twin-free, edge-maximal locally bipartite graphs all beating $1/2$: define $\Delta_{\ell}$ as the complement of $C_{4\ell - 1}^{\ell - 1}$ (this is very natural as $\Gamma_{i}$ is the complement of $C_{3i - 1}^{i - 1}$). Then $\Delta_{\ell}$ has $4\ell - 1$ vertices, is $(2\ell)$-regular, is 4-chromatic (its independence number is $\ell$), and is edge-maximal locally bipartite (the addition of any edge gives a 4-clique). Note that $\Delta_{2} = \overline{C}_{7}$. In fact, $\Delta_{3}$ satisfies $\delta(\Delta_{3})/\abs{\Delta_{3}} = 6/11$ suggesting it is the next key graph when extending \cref{spec4localbip} below $6/11$. Unlike the triangle-free case, the $\Delta_{\ell}$ are not the only 4-chromatic twin-free, edge-maximal locally bipartite graphs beating $1/2$. Indeed, $H_{2}^{+}$ is not a $\Delta_{\ell}$ and nor is the graph shown in \cref{fig:AnotherCounterexample} on \cpageref{fig:AnotherCounterexample}. Intriguingly, neither of these graphs is contained in (nor, by \cref{homscores}, has a homomorphism to) any $\Delta_{\ell}$, since no $\Delta_{\ell}$ contains an induced $H_{2}$ (no neighbourhood in $\Delta_{\ell}$ contains two edges with no edges between). It would be interesting to have an infinite sequence of such non-$\Delta_{\ell}$ graphs. Also, for each $c > 1/2$, are there only finitely many twin-free, edge-maximal locally bipartite graphs beating $c$ (for triangle-free graphs this was first shown by \L uczak~\cite{Luczak2006})?

A final question is whether there are any locally bipartite graphs beating $1/2$ that are not 4-colourable -- such graphs would be the analogue of Vega graphs in the triangle-free case. If there were none, then $\delta_{\chi}(\cF_{1, 2}, 5) = 1/2 = \delta_{\chi}(\cF_{1, 2})$ and so the chromatic profile of locally bipartite graphs would have only two thresholds ($4/7$ and $1/2$) compared to three ($2/5$, $10/29$ and $1/3$) for triangle-free graphs.

\subsection{Notation}\label{sec:notation}

Let $G$ be a graph and $X \subset V(G)$. We write \defn{$\Gamma(X)$} for $\cap_{v \in X} \Gamma(v)$ (the common neighbourhood of the vertices of $X$) and \defn{$d(X)$} for $\abs{\Gamma(X)}$. We often omit set parentheses so $\Gamma(u, v) = \Gamma(u) \cap \Gamma(v)$ and $d(u, v) = \abs{\Gamma(u, v)}$. We write \defn{$G_{X}$} for $G[\Gamma(X)]$ so, for example, $G_{u, v}$ is the induced graph on the common neighbourhood of vertices $u$ and $v$. We make frequent use of the fact that for two vertices $u$ and $v$ of $G$
\begin{equation*}
	d(u, v) = d(u) + d(v) - \abs{\Gamma(u) \cup \Gamma(v)} \geqslant d(u) + d(v) - \abs{G} \geqslant 2 \delta(G) - \abs{G}.
\end{equation*}
Given a set of vertices $X \subset V(G)$, we write \defn{$e(X, G)$} for the number of ordered pairs of vertices $(x, v)$ with $x \in X$, $v \in G$ and $xv$ an edge in $G$. In particular, $e(X, G)$ counts each edge in $G[X]$ twice and each edge from $X$ to $G - X$ once and satisfies
\begin{equation*}
	e(X, G) = \sum_{x \in X} d(x) = \sum_{v \in G} \abs{\Gamma(v) \cap X}.
\end{equation*}
We generalise this notation to vertex weightings which will appear in many of our arguments. We will take a set of vertices $X \subset V(G)$ and assign weights $\omega \colon X \to \ZZN$ to the vertices of $X$. Then we define
\begin{equation*}
	\omega(X, G) = \sum_{x \in X} \omega(x) d(x) = \sum_{v \in G} \text{Total weight of the neighbours of } v \text{ in } X.
\end{equation*}
We will often use the word \defn{circuit} (as opposed to cycle) in our arguments. A circuit is a sequence of (not necessarily distinct) vertices $v_{1}, v_{2}, \dotsc, v_{\ell}$ with $\ell > 1$, $v_{i}$ adjacent to $v_{i + 1}$ (for $i = 1, 2, \dotsc, \ell - 1$) and $v_{\ell}$ adjacent to $v_{1}$. Note that in a locally bipartite graph the neighbourhood of any vertex does not contain an odd circuit (and, of course, does not contain an odd cycle). We use circuit to avoid considering whether some pairs of vertices are distinct when it is unnecessary to do so.

For two graphs $G$ and $H$, we say there is a \defn{homomorphism $G \to H$} if there is a map $\varphi \colon V(G) \to V(H)$ such that for every edge $uv$ of $G$, $\varphi(u)\varphi(v)$ is an edge of $H$. Note that there is a homomorphism $G \to H$ if and only if $G$ is a subgraph of some blow-up of $H$. In particular, if there is a homomorphism $G \to H$, then $\chi(G) \leqslant \chi(H)$ and moreover if $H$ is locally bipartite, then $G$ is also.

\section{\texorpdfstring{\Cref{spec4localbip}}{Theorem 2} without homomorphisms}\label{sec:localbip}

\Cref{spec4localbip} follows from \cref{main4localbip,lemma4H,hom2C,hom2H,hom2Hepsilon} which we state here. \Cref{main4localbip} establishes containing $H_{0}$ as an obstruction to a locally bipartite graph being 3-colourable. \Cref{lemma4H} leverages $H_{0}$ up to $H_{2}^{+}$ and $\overline{C}_{7}$. We prove these two results in this section: they give the required starting structure for the proofs of the homomorphism results, \cref{hom2C,hom2H,hom2Hepsilon}, which we carry out in \cref{sec:hom}.

\begin{restatable}{theorem}{mainforlocalbip}\label{main4localbip}
	Let $G$ be a locally bipartite graph. If $\delta(G) > 6/11 \cdot \abs{G}$, then $G$ is either $3$-colourable or contains $H_{0}$.
\end{restatable}

\begin{theorem}\label{lemma4H}
	Let $G$ be a locally bipartite graph which contains $H_{0}$.
	\begin{itemize}[noitemsep]
		\item Firstly, it must be the case that $\delta(G) \leqslant 4/7 \cdot \abs{G}$.
		\item Secondly, if $\delta(G) > 5/9 \cdot \abs{G}$, then $G$ contains $\overline{C}_{7}$.
		\item Thirdly, if $\delta(G) > 6/11 \cdot \abs{G}$, then $G$ contains $H_{2}^{+}$ or $\overline{C}_{7}$.
	\end{itemize}
\end{theorem}

The graph $\overline{C}_{7}$ (and any of its balanced blow-ups) show that 4/7 is tight. Below are weightings (blow-ups) of $H_{2}^{+}$ and $H_{2}$ -- in the former, $0^{+}$ represents some tiny positive weight (we have not deleted the vertices entirely just given them a small weight relative to the rest). These blow-ups show that $5/9$ and $6/11$ are tight respectively (it follows from \cref{homscores} that there is no homomorphism from $H_{2}^{+}$ to $\overline{C}_{7}$, and no homomorphism from either of $H_{2}^{+}$ or $\overline{C}_{7}$ to $H_{2}$).

\begin{figure}[H]
	\centering
	\begin{subfigure}{.33\textwidth}
		\centering
		\begin{tikzpicture}
			\foreach \pt in {0,1,...,6} 
			{
				\tkzDefPoint(\pt*360/7 + 90:1){v_\pt}
			} 
			\tkzDefPoint(0,0){u}
			\tkzDrawPolySeg(v_0,v_1,v_2,v_3,v_4,v_5,v_6,v_0)
			\tkzDrawPolySeg(v_1,v_3,v_5,v_0, v_2,v_4,v_6)
			\tkzDrawSegments(u,v_0 u,v_2 u,v_5)
			\tkzDrawPoints(v_0,v_...,v_6)
			\tkzDrawPoint(u)
			\tkzLabelPoint[above](v_0){2}
			\tkzLabelPoint[left](v_1){$0^{+}$}
			\tkzLabelPoint[left](v_2){2}
			\tkzLabelPoint[below](v_3){1}
			\tkzLabelPoint[below](v_4){1}
			\tkzLabelPoint[right](v_5){2}
			\tkzLabelPoint[right](v_6){$0^{+}$}
			\tkzLabelPoint[below](u){1}
		\end{tikzpicture}
		\subcaption*{$H_{2}^{+}$ weighted}
	\end{subfigure}
	\begin{subfigure}{.33\textwidth}
		\centering
		\begin{tikzpicture}
			\foreach \pt in {0,1,...,6} 
			{
				\tkzDefPoint(\pt*360/7 + 90:1){v_\pt}
			} 
			\tkzDrawPolySeg(v_0,v_1,v_2,v_3,v_4,v_5,v_6,v_0)
			\tkzDrawPolySeg(v_1,v_3,v_5,v_0, v_2,v_4,v_6)
			\tkzDrawPoints(v_0,v_...,v_6)
			\tkzLabelPoint[above](v_0){3}
			\tkzLabelPoint[left](v_1){1}
			\tkzLabelPoint[left](v_2){2}
			\tkzLabelPoint[below](v_3){1}
			\tkzLabelPoint[below](v_4){1}
			\tkzLabelPoint[right](v_5){2}
			\tkzLabelPoint[right](v_6){1}
		\end{tikzpicture}
		\subcaption*{$H_{2}$ weighted}
	\end{subfigure}
	\caption{Weightings of $H^{+}_{2}$ and $H_{2}$.}\label{fig:weightings}
\end{figure}

\begin{restatable}{theorem}{homtoC}\label{hom2C}
	Let $G$ be a locally bipartite graph. If $\delta(G) > 6/11 \cdot \abs{G}$ and $G$ contains $\overline{C}_{7}$, then there is a homomorphism $G \to \overline{C}_{7}$.
\end{restatable}

\vspace{-10pt}

\begin{restatable}{theorem}{homtoH}\label{hom2H}
	Let $G$ be a locally bipartite graph. If $\delta(G) > 6/11 \cdot \abs{G}$, then $G$ is 4-colourable.
\end{restatable}

\vspace{-10pt}

\begin{restatable}{theorem}{homtoHepsilon}\label{hom2Hepsilon}
	There is an $\varepsilon > 0$ such that if $G$ is a locally bipartite graph with $\delta(G) > (5/9 - \varepsilon) \cdot \abs{G}$ and $G$ does not contain $\overline{C}_{7}$, then there is a homomorphism $G \to H_{2}^{+}$.
\end{restatable}

\begin{remark}
	We make no attempt to optimise the proof to obtain the `best value' of $\varepsilon$ as we believe that it is in fact possible to replace $5/9 - \varepsilon$ by $6/11$ (but our arguments do not achieve this).
\end{remark}

\begin{Proof}[of \cref{spec4localbip}]
	Let $G$ be a locally bipartite graph. \Cref{main4localbip,lemma4H} together show that 
	\begin{itemize}[noitemsep]
		\item If $\delta(G) > 4/7 \cdot \abs{G}$, then $G$ is 3-colourable.
		\item If $\delta(G) > 5/9 \cdot \abs{G}$, then $G$ is either 3-colourable or contains $\overline{C}_{7}$.
		\item If $\delta(G) > 6/11 \cdot \abs{G}$, then $G$ is either 3-colourable or contains $H_{2}^{+}$ or contains $\overline{C}_{7}$.
	\end{itemize}
	The first bullet point of \cref{spec4localbip} is immediate and the second follows from \cref{hom2C}.
	
	Let $\varepsilon > 0$ as in \cref{hom2Hepsilon}. Suppose $\delta(G) > (5/9 - \varepsilon) \cdot \abs{G}$. If $G$ contains $\overline{C}_{7}$, then there is a homomorphism $G \to \overline{C}_{7}$, by \cref{hom2C}. Otherwise, by \cref{hom2Hepsilon}, there is a homomorphism $G \to H_{2}^{+}$. This gives the third bullet point of \cref{spec4localbip}. The final bullet point of \cref{spec4localbip} follows immediately from \cref{hom2H}. 
	
	The balanced blow-up of $\overline{C}_{7}$ on $n$ vertices has minimum degree at least $4 \lfloor n/7 \rfloor$ and is 4-chromatic showing that $4/7$ is tight. Finally the blow-up of $H_{2}^{+}$ displayed in \cref{fig:weightings} shows that $5/9$ is tight (by \cref{homscores}, there is no homomorphism from $\overline{C}_{7}$ to $H_{2}^{+}$).
\end{Proof}

In \cref{sec:H02C7} we carry out a careful edge-counting/vertex weighting argument which proves \cref{lemma4H}. The proofs of \cref{hom2C,hom2H,hom2Hepsilon} are deferred to \cref{sec:hom}. We now introduce a key definition that will be crucial for our proofs and particularly for the proof of \cref{main4localbip}. To motivate this, consider a locally bipartite, $H_{0}$-free graph $G$ with $\delta(G) > 6/11 \cdot \abs{G}$. To prove \cref{main4localbip} we need to show that $G$ is $3$-colourable. We may as well assume that $G$ is edge-maximal: that is, the addition of any edge to $G$ introduces either a copy of $H_{0}$ or creates a vertex with a non-bipartite neighbourhood. Thus, any non-edge of $G$ is either a missing edge of a $K_{4}$, a missing rim of an odd wheel, a missing spoke of an odd wheel, or a missing edge of an $H_{0}$. This motivates a key definition.

\begin{definition}[dense and sparse]\label{def:dense}
	A pair of \emph{non-adjacent, distinct} vertices $u, v$ in a graph $G$ is \defn{dense} if $G_{u, v}$ contains an edge and \defn{sparse} if $G_{u, v}$ does not contain an edge.
\end{definition}

First note that every pair of distinct vertices in any graph is exactly one of `adjacent', `dense' or `sparse'. Another way to view being dense is as being the missing edge of a $K_{4}$. Locally bipartite graphs are $K_4$-free so any pair of distinct vertices with an edge in their common neighbourhood must be non-adjacent and so dense. Our initial observations above show that a sparse pair in the edge-maximal $G$ is either a missing edge of an $H_{0}$, a missing rim of an odd wheel or a missing spoke of an odd wheel. In \cref{sec:oddspokes} we will rule out the possibility that a sparse pair is the missing spoke of an odd wheel and in \cref{sec:proof4localbip} we will complete the proof of \cref{main4localbip}.

We finish the introduction to this section by collecting four simple but very effective lemmas about dense pairs of vertices. The second of these exhibits an edge-counting method that we will use frequently. These give us some control over dense pairs and much of the proof of \cref{main4localbip} involves understanding in what configurations sparse pairs appear.

\begin{lemma}\label{lemma4I}
	Let $G$ be a graph with $\delta(G) > 1/2 \cdot \abs{G}$ and let $I$ be any largest independent set in $G$. Then, for every distinct $u, v \in I$, the pair $u, v$ is dense.
\end{lemma}

\begin{Proof}
	Fix distinct $u, v \in I$. Note that $\Gamma(u), \Gamma(v) \subset V(G) \setminus I$ so $\abs{\Gamma(u) \cup \Gamma(v)} \leqslant \abs{G} - \abs{I}$. Hence
	\begin{align*}
		\abs{\Gamma(u) \cap \Gamma(v)} & = d(u) + d(v) - \abs{\Gamma(u) \cup \Gamma(v)} \geqslant 2 \delta(G) - (\abs{G} - \abs{I}) \\
		& = \abs{I} + 2 \delta(G) - \abs{G} > \abs{I}.
	\end{align*}
	But $I$ is a largest independent set in $G$ and so $\Gamma(u) \cap \Gamma(v)$ is not independent: $G_{u, v}$ contains an edge so $u, v$ is dense.
\end{Proof}

\begin{lemma}\label{lemma4sparse}
	Let $G$ be a graph with $\delta(G) > 1/2 \cdot \abs{G}$ and suppose $C$ is an induced 4-cycle in $G$. Then at least one of the non-edges of $C$ is a dense pair.
\end{lemma}

\begin{Proof}
	Suppose the result does not hold. We have an induced 4-cycle $C = v_{1}v_{2}v_{3}v_{4}$ in $G$ with edges $v_{1}v_{2}$, $v_{2}v_{3}$, $v_{3}v_{4}$, $v_{4}v_{1}$ where the pairs $v_{1}, v_{3}$ and $v_{2}, v_{4}$ are both sparse. Note that any vertex has at most two neighbours in $C$. Indeed if $u$ is adjacent to both $v_{1}$ and $v_{3}$, then $u$ cannot be adjacent to either $v_{2}$ or $v_{4}$, as the pair $v_{1}, v_{3}$ is sparse; similarly, if $u$ is adjacent to both $v_{2}$ and $v_{4}$, then $u$ cannot be adjacent to either $v_{1}$ or $v_{3}$. Counting the edges between $C$ and $G$ from both sides gives
	\begin{equation*}
		4\delta(G) \leqslant d(v_{1}) + d(v_{2}) + d(v_{3}) + d(v_{4}) = e(C, G) \leqslant 2 \cdot \abs{G},
	\end{equation*}
	which contradicts $\delta(G) > 1/2 \cdot \abs{G}$.
\end{Proof}

\begin{lemma}\label{lemma4dense}
	Let $G$ be a locally bipartite graph which does not contain $H_{0}$. For any vertex $v$ of $G$,
	\begin{equation*}
		D_{v} \coloneqq \set{u \colon \text{the pair } u, v \text{ is dense}}
	\end{equation*}
	is an independent set of vertices.
\end{lemma}

\begin{Proof}
	Suppose that in fact there are distinct vertices $v$, $u_{1}$ and $u_{2}$ with the pairs $v, u_{1}$ and $v, u_{2}$ both dense and with $u_{1}$ adjacent to $u_{2}$. Let $x_{1}x_{2}$ be an edge in the common neighbourhood of $v$ and $u_{1}$ and $x_{3}x_{4}$ be an edge in the common neighbourhood of $v$ and $u_{2}$.
	
	If $\set{x_{1}, x_{2}} = \set{x_{3}, x_{4}}$, then $u_{1}u_{2}x_{1}x_{2}$ is a $K_{4}$ in $G$. If $\set{x_{1}, x_{2}}$ and $\set{x_{3}, x_{4}}$ have one element in common, say $x_{1} = x_{3}$, then $G_{x_{1}}$ contains the 5-cycle $vx_{2}u_{1}u_{2}x_{4}$. Finally, if $\set{x_{1}, x_{2}}$ and $\set{x_{3}, x_{4}}$ are disjoint, then $G[\set{v, x_{1}, x_{2}, u_{1}, u_{2}, x_{3}, x_{4}}]$ contains a copy of $H_{0}$.
\end{Proof}

We can combine \cref{lemma4I,lemma4dense} to give the following which will play a crucial role in finishing the proof of \cref{main4localbip}.

\begin{lemma}\label{lemma4DuI}
	Let $G$ be an $H_{0}$-free, locally bipartite graph with $\delta(G) > 1/2 \cdot \abs{G}$. Let $I$ be any largest independent set in $G$. For any distinct vertices $u, v$ with $u \in I$\textnormal{:} $v \in I$ if and only if the pair $u, v$ is dense.
\end{lemma}

\begin{Proof}
	Consider the set $\{u\} \cup D_{u}$ consisting of $u$ and all the vertices which form a dense pair with $u$. It suffices to show that $I = \{u\} \cup D_{u}$. By \cref{lemma4I}, $I \subset \{u\} \cup D_{u}$. However, by the definition of dense and \cref{lemma4dense}, $\{u\} \cup D_{u}$ is an independent set. Hence, by the maximality of $I$, we have $I = \{u\} \cup D_{u}$.
\end{Proof}

\subsection{From \texorpdfstring{$H_{0}$}{H0} to \texorpdfstring{$\overline{C}_{7}$}{C7}}\label{sec:H02C7}

In this subsection we prove \cref{lemma4H}. The strategy is to start with a copy of $H_0$ and consider the edges between it and the rest of $G$. Using the high minimum degree we are able to find a vertex with the correct neighbours in the copy of $H_0$ so that a copy of $H_1$ is present. We then play the same game to get a copy of $H_2$ and a copy of $H_2^+$ or $\overline{C}_7$.

For ease of reading we split the proof of \cref{lemma4H} into a sequence of claims from $6/11$ up to $4/7$. Each claim corresponds to a bullet point of \cref{lemma4H} and we are addressing the bullet points in reverse order. The proof of the first claim is by far the longest. 

\begin{claim}\label{claim:H02H21C7}
	Let $G$ be a locally bipartite graph containing $H_{0}$. If $\delta(G) > 6/11 \cdot \abs{G}$, then $G$ contains $H_{2}^{+}$ or $\overline{C}_{7}$.
\end{claim}

\begin{Proof}
	This proof has the following structure. We will first show that $G$ contains $H_{1}$, then that it contains $H_{2}$ and finally that it contains one of $H_{2}^{+}$ or $\overline{C}_{7}$. We label a copy of $H_{0}$ in $G$ as below and let $X = \set{a_{0}, a_{1}, \dotsc, a_{6}}$. Our first aim is to show that $G$ contains $H_{1}$.
	
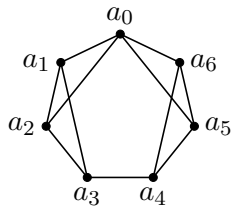
\begin{figure}[H]
		\centering
		\begin{tikzpicture}
			\foreach \pt in {0,1,...,6} 
			{
				\tkzDefPoint(\pt*360/7 + 90:1){v_\pt}
			} 
			\tkzDrawPolySeg(v_0,v_1,v_2,v_3,v_4,v_5,v_6,v_0) 
			\tkzDrawPolySeg(v_5,v_0,v_2)
			\tkzDrawPolySeg(v_1,v_3)
			\tkzDrawPolySeg(v_4,v_6)
			\tkzDrawPoints(v_0,v_...,v_6)
			\tkzLabelPoint[above](v_0){$a_{0}$}
			\tkzLabelPoint[left](v_1){$a_{1}$}
			\tkzLabelPoint[left](v_2){$a_{2}$}
			\tkzLabelPoint[below](v_3){$a_{3}$}
			\tkzLabelPoint[below](v_4){$a_{4}$}
			\tkzLabelPoint[right](v_5){$a_{5}$}
			\tkzLabelPoint[right](v_6){$a_{6}$}
		\end{tikzpicture}
	\end{figure}
	
	Let $U_{4}$ be the set of vertices with exactly four neighbours in $X$. \Cref{rmk:5nbs} says that no vertex has five neighbours in a copy of $H_{0}$ so all other vertices have at most three neighbours in $X$, and hence
	\begin{equation*}
		42/11 \cdot \abs{G} < 7 \delta (G) \leqslant e(X, G) \leqslant 4 \abs{U_{4}} + 3(\abs{G} - \abs{U_{4}}) = 3 \abs{G} + \abs{U_{4}},
	\end{equation*}
	and so 
	\begin{equation*}
		\abs{U_{4}} > 9/11 \cdot \abs{G}.
	\end{equation*}
	Now $\abs{U_{4}} + d(a_{0}) > \abs{G}$ and so some vertex $v$ is adjacent to $a_{0}$ and has four neighbours in $X$. Note that $v$ cannot be adjacent to both $a_{1}$, $a_{2}$ as otherwise $va_{0}a_{1}a_{2}$ is a $K_{4}$, so by symmetry we may assume that $v$ is not adjacent to $a_{1}$. Similarly we may assume that $v$ is not adjacent to $a_{5}$. But $v$ has four neighbours in $X$ so must be adjacent to at least one of $a_{2}$, $a_{6}$ -- by symmetry we may assume $v$ is adjacent to $a_{2}$.
	
	There are two possibilities: $v$ is adjacent to $a_{0}, a_{2}, a_{3}, a_{4}$, or $v$ is adjacent to $a_{0}, a_{2}, a_{6}$ and one of $a_{3}$, $a_{4}$. In the latter case we may assume by symmetry that $v$ is adjacent to $a_3$. Hence there are two possibilities for $\Gamma(v) \cap X$: $\set{a_0, a_2, a_3, a_4}$ and $\set{a_0, a_2, a_3, a_6}$. In both cases $v$ cannot be any $a_i$ except for possibly $a_1$.
	
	If $\Gamma(v) \cap X = \set{a_0, a_2, a_3, a_4}$, then $G$ contains $H_{1}$ (replace $a_{1}$ by $v$). If $\Gamma(v) \cap X = \set{a_0, a_2, a_3, a_6}$, then $G$ contains the following graph where, in particular, $v$ is not adjacent to $a_4$.
	\begin{figure}[H]
		\centering
		\begin{tikzpicture}
			\foreach \pt in {0,1,...,6} 
			{
				\tkzDefPoint(\pt*360/7 + 90:1){v_\pt}
			} 
			\tkzDrawPolySeg(v_0,v_1,v_2,v_3,v_4,v_5,v_6,v_0) 
			\tkzDrawPolySeg(v_4,v_6,v_1,v_3)
			\tkzDrawPolySeg(v_5,v_0,v_2)
			\tkzDrawPoints(v_0,v_...,v_6)
			\tkzLabelPoint[above](v_0){$a_{0}$}
			\tkzLabelPoint[left](v_1){$v$}
			\tkzLabelPoint[left](v_2){$a_{2}$}
			\tkzLabelPoint[below](v_3){$a_{3}$}
			\tkzLabelPoint[below](v_4){$a_{4}$}
			\tkzLabelPoint[right](v_5){$a_{5}$}
			\tkzLabelPoint[right](v_6){$a_{6}$}
		\end{tikzpicture}
	\end{figure}
	Vertex $a_6$ is not adjacent to $a_3$ else $G_{a_6}$ contains the 5-cycle $a_{0} v a_{3} a_{4} a_{5}$. Hence $va_{6}a_{4}a_{3}$ is an induced 4-cycle in $G$. By \cref{lemma4sparse}, at least one of the pairs $v, a_{4}$ and $a_{3}, a_{6}$ is dense. By symmetry we may assume that $v, a_{4}$ is dense: let $a'_{2}a'_{3}$ be an edge in $G_{v, a_{4}}$. Vertex $v$ is not adjacent to $a_5$, so $a_5$ is neither $a'_2$ nor $a'_3$. Note that $a_0$ is not adjacent to $a_4$ (else $a_{0} a_{4} a_{5} a_{6}$ is a $K_{4}$) and so $a_{0}$ is neither $a'_{2}$ nor $a'_{3}$. If $a_{6} = a'_{2}$, then $G_{a_{6}}$ contains the 5-cycle $a'_{3}va_{0}a_{5}a_{4}$, which is impossible. Similarly $a_{6} \neq a'_{3}$. Hence, $a'_{2}, a'_{3}$ are distinct from $a_{4}, a_{5}, a_{6}, a_{0}, v$ and so $G[\{a_{6}, a_{0}, v, a'_{2}, a'_{3}, a_{4}, a_{5}\}]$ contains a copy of $H_{1}$ (with apex $a_{6}$). Hence in all cases $G$ contains a copy of $H_{1}$.
	
	We now show that $G$ contains a copy of $H_{2}$. Consider a copy of $H_{1}$ with vertices $X = \{a_{0}, a_{1}, \dotsc, a_{6}\}$ as shown below and again let $U_{4}$ be the set of vertices with 4 neighbours in $X$. As before, $\abs{U_{4}} > 9/11 \cdot \abs{G}$.
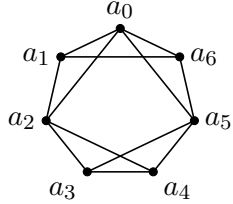
	\begin{figure}[H]
		\centering
		\begin{tikzpicture}
			\foreach \pt in {0,1,...,6} 
			{
				\tkzDefPoint(\pt*360/7 + 90:1){v_\pt}
			} 
			\tkzDrawPolySeg(v_0,v_1,v_2,v_3,v_4,v_5,v_6,v_0) 
			\tkzDrawPolySeg(v_3,v_5,v_0,v_2,v_4)
			\tkzDrawPolySeg(v_6,v_1)
			\tkzDrawPoints(v_0,v_...,v_6)
			\tkzLabelPoint[above](v_0){$a_{0}$}
			\tkzLabelPoint[left](v_1){$a_{1}$}
			\tkzLabelPoint[left](v_2){$a_{2}$}
			\tkzLabelPoint[below left](v_3){$a_{3}$}
			\tkzLabelPoint[below right](v_4){$a_{4}$}
			\tkzLabelPoint[right](v_5){$a_{5}$}
			\tkzLabelPoint[right](v_6){$a_{6}$}
		\end{tikzpicture}
	\end{figure}
	No vertex is adjacent to all of $a_{0}, a_{2}, a_{5}$, else $G_{a_{0}}$ contains an odd circuit. In particular, $\Gamma(a_{0}, a_{2})$, $\Gamma(a_{0}, a_{5})$ are disjoint. Each of these sets has size at least $2 \delta(G) - \abs{G} > 1/11 \cdot \abs{G}$. But then $\abs{U_{4}} + d(a_{0}, a_{2}) + d(a_{0}, a_{5}) > \abs{G}$ and so there is some vertex $v \in \Gamma(a_{0}, a_{2}) \cup \Gamma(a_{0}, a_{5})$ with four neighbours in $X$. By symmetry, we may assume $v$ is adjacent to both $a_{0}$ and $a_{2}$ and not to $a_{5}$. Also $v$ cannot be adjacent to $a_{1}$ otherwise $va_{0}a_{1}a_{2}$ is a $K_{4}$. Similarly $v$ cannot be adjacent to both $a_{3}$ and $a_{4}$. Hence, $v$ is adjacent to $a_{0}, a_{2}, a_{6}$ and one of $a_{3}, a_{4}$. By symmetry, we may assume $v$ is adjacent to $a_{3}$. If $v = a_{5}$, then $a_{2}a_{5}$ is an edge, so $G$ contains a $K_{4}$. If $v = a_{4}$, then $a_{0}$ has five neighbours in $X$, a contradiction. If $v$ is neither $a_{4}$ nor $a_{5}$, then $G$ contains $H_{2}$ (replace $a_{1}$ by $v$).
	
	We finally show that $G$ contains a copy of $H_{2}^{+}$ or $\overline{C}_{7}$. Consider a copy of $H_{2}$ in $G$ with vertices $X = \set{a_{0}, a_{1}, \dotsc, a_{6}}$. We assign weights $\omega \colon X \to \ZZN$ as shown in the diagram below, so, for example, $\omega(a_{0}) = 3$ and $\omega(a_{1}) = 1$ (recall this notation from \cref{sec:notation}). For each vertex $v \in G$, let $f(v)$ be the total weight of the neighbours of $v$ in $X$.
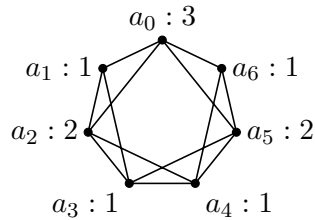
	\begin{figure}[H]
		\centering
		\begin{tikzpicture}
			\foreach \pt in {0,1,...,6} 
			{
				\tkzDefPoint(\pt*360/7 + 90:1){v_\pt}
			} 
			\tkzDrawPolySeg(v_0,v_1,v_2,v_3,v_4,v_5,v_6,v_0)
			\tkzDrawPolySeg(v_1,v_3,v_5,v_0, v_2,v_4,v_6)
			\tkzDrawPoints(v_0,v_...,v_6)
			\tkzLabelPoint[above](v_0){$a_{0} \colon 3$}
			\tkzLabelPoint[left](v_1){$a_{1} \colon 1$}
			\tkzLabelPoint[left](v_2){$a_{2} \colon 2$}
			\tkzLabelPoint[below left](v_3){$a_{3} \colon 1$}
			\tkzLabelPoint[below right](v_4){$a_{4} \colon 1$}
			\tkzLabelPoint[right](v_5){$a_{5} \colon 2$}
			\tkzLabelPoint[right](v_6){$a_{6} \colon 1$}
		\end{tikzpicture}
	\end{figure}
	
	Now
	\begin{equation*}
		6 \abs{G} < 11 \delta(G) \leqslant \omega(X, G) =  \sum_{v \in G} f(v),
	\end{equation*}
	so some vertex $v$ has $f(v) \geqslant 7$. But all vertices have at most four neighbours in a copy of $H_{2}$ so either $v$ is adjacent to all of $a_{0}$, $a_{2}$, $a_{5}$ or $v$ is adjacent to $a_{0}$, to exactly one of $a_{2}$ and $a_{5}$, and to exactly two of $a_{1}$, $a_{3}$, $a_{4}$, $a_{6}$.
	
	First suppose that $v$ is adjacent to all of $a_{0}$, $a_{2}$ and $a_{5}$. Note that $v$ cannot be in $X$. Indeed, if $v = a_{1}$, then $G_{a_{5}}$ contains the 5-cycle $a_{0}va_{3}a_{4}a_{6}$ and similarly if $v = a_{6}$. On the other hand, if $v = a_{3}$ then $a_{0}a_{1}a_{2}v$ is a $K_{4}$ and similarly if $v = a_{4}$. Hence $v$ together with $H_{2}$ gives a copy of $H_{2}^{+}$ in $G$.
	
	Second suppose that $v$ is adjacent to $a_{0}$, to one of $a_{2}$ and $a_{5}$, and to two of $a_{1}$, $a_{3}$, $a_{4}$, $a_{6}$. By symmetry we may assume that $v$ is adjacent to $a_{2}$ and not to $a_{5}$. Then $v$ is not adjacent to $a_{1}$ else $va_{0}a_{1}a_{2}$ is a $K_{4}$ and $v$ is not adjacent to $a_{4}$ else $va_{0}a_{1}a_{3}a_{4}$ is an odd circuit in $G_{a_{2}}$. Thus $v$ is adjacent to $a_{0}, a_{2}, a_{3}$ and $a_{6}$. Note $v$ is neither $a_{4}$ nor $a_{5}$ as $v$ does not have five neighbours in $X$. Thus $G[X \setminus \set{a_{1}} \cup \set{v}]$ contains a copy of $\overline{C}_{7}$.
\end{Proof}

\begin{claim}
	Let $G$ be a locally bipartite graph containing $H_{0}$. If $\delta(G) > 5/9 \cdot \abs{G}$, then $G$ contains $\overline{C}_{7}$.
\end{claim}

\begin{Proof}
	From \cref{claim:H02H21C7}, $G$ contains $H_{2}$. Let $X = \set{a_{0}, a_{1}, \dotsc, a_{6}}$ be a copy of $H_{2}$ in $G$.
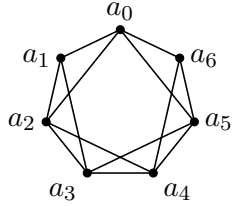
	\begin{figure}[H]
		\centering
		\begin{tikzpicture}
			\foreach \pt in {0,1,...,6} 
			{
				\tkzDefPoint(\pt*360/7 + 90:1){v_\pt}
			} 
			\tkzDrawPolySeg(v_0,v_1,v_2,v_3,v_4,v_5,v_6,v_0)
			\tkzDrawPolySeg(v_1,v_3,v_5,v_0, v_2,v_4,v_6)
			\tkzDrawPoints(v_0,v_...,v_6)
			\tkzLabelPoint[above](v_0){$a_{0}$}
			\tkzLabelPoint[left](v_1){$a_{1}$}
			\tkzLabelPoint[left](v_2){$a_{2}$}
			\tkzLabelPoint[below left](v_3){$a_{3}$}
			\tkzLabelPoint[below right](v_4){$a_{4}$}
			\tkzLabelPoint[right](v_5){$a_{5}$}
			\tkzLabelPoint[right](v_6){$a_{6}$}
		\end{tikzpicture}
	\end{figure}
	
	Let $U_{4}$ be the set of vertices with exactly four neighbours in $X$. All other vertices have at most three neighbours in $X$ so $\abs{U_{4}} \geqslant 7 \delta(G) - 3 \abs{G}$. Thus
	\begin{equation*}
		\abs{U_{4}} + \abs{\Gamma(a_{0}, a_{2})} \geqslant 7 \delta(G) - 3 \abs{G} + 2\delta(G) - \abs{G} = 9 \delta(G) - 4 \abs{G} > \abs{G},
	\end{equation*}
	so $U_{4}$ and $\Gamma(a_{0}, a_{2})$ are not disjoint: there is a vertex $v$ which is adjacent to both $a_{0}$ and $a_{2}$ and has four neighbours in $X$. Vertex $v$ is not adjacent to $a_{1}$ otherwise $va_{0}a_{1}a_{2}$ is a $K_{4}$. Also $v$ is not adjacent to $a_{4}$ otherwise $G_{a_{2}}$ contains the odd circuit $va_{0}a_{1}a_{3}a_{4}$. Note $v$ is not adjacent to both $a_{5}$ and $a_{6}$ otherwise $va_{0}a_{5}a_{6}$ is a $K_{4}$. Hence $v$ is adjacent to $a_{3}$. But then $v$ is not adjacent to $a_{5}$ else $G_{a_{5}}$ contains the odd circuit $va_{0}a_{6}a_{4}a_{3}$.
	
	Therefore $\Gamma(v) \cap X = \set{a_{0}, a_{2}, a_{3}, a_{6}}$. In particular, $v$ is neither $a_{4}$ nor $a_{5}$. Thus $G[X \setminus \set{a_{1}} \cup \set{v}]$ contains $\overline{C}_{7}$.
\end{Proof}

\begin{claim}
	Let $G$ be a locally bipartite graph containing $H_{0}$. Then $\delta(G) \leqslant 4/7 \cdot \abs{G}$.
\end{claim}

\begin{Proof}
	Let $X$ be a set of seven vertices in $G$ with $G[X]$ containing $H_{0}$. By \cref{rmk:5nbs}, every vertex has at most four neighbours in $X$ so
	\begin{equation*}
		7 \delta(G) \leqslant e(X, G) \leqslant 4 \abs{G}. \qedhere
	\end{equation*}
\end{Proof}

\subsection{Ruling out sparse pairs being spokes of odd wheels}\label{sec:oddspokes}

In this subsection we make a start on the proof of \cref{main4localbip}, by ruling out the possibility that $G$ contains a sparse pair of vertices which is the spoke of an odd wheel.

\begin{lemma}\label{lemma45wheel}
	Let $G$ be a locally bipartite graph with $\delta(G) > 1/2 \cdot \abs{G}$ and which does not contain $H_{0}$. Then $G$ does not contain a sparse pair $u, v$ with $uv$ being the missing spoke of a 5-wheel.
\end{lemma}

\begin{Proof}
	Suppose the conclusion does not hold: label the configuration as follows where $u, v$ is a sparse pair.
	
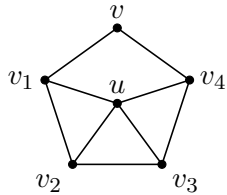
\begin{figure}[H]
		\centering
		\begin{tikzpicture}
			\foreach \pt in {0,1,...,4} 
			{
				\tkzDefPoint(\pt*360/5 + 90:1){v_\pt}
			} 
			\tkzDefPoint(0,0){u}
			\tkzDrawPolySeg(v_0,v_1,v_2,v_3,v_4,v_0) 
			\foreach \pt in {1,...,4} 
			{
				\tkzDrawSegment(u,v_\pt)
			} 
			\tkzDrawPoints(v_0,v_...,v_4)
			\tkzDrawPoint(u)
			\tkzLabelPoint[above](v_0){$v$}
			\tkzLabelPoint[left](v_1){$v_{1}$}
			\tkzLabelPoint[below left](v_2){$v_{2}$}
			\tkzLabelPoint[below right](v_3){$v_{3}$}
			\tkzLabelPoint[right](v_4){$v_{4}$}
			\tkzLabelPoint[above](u){$u$}
		\end{tikzpicture}
	\end{figure}
	As $u, v$ is sparse, $v_{1}$ is not adjacent to $v_{4}$ and so $uv_{1}vv_{4}$ is an induced 4-cycle in $G$. By \cref{lemma4sparse}, the pair $v_{1}, v_{4}$ must be dense. But the pair $v_1, v_3$ is also dense and so $D_{v_{1}} = \set{x \colon v_1, x \text{ is dense}}$ contains the edge $v_{3} v_{4}$. This, however, contradicts \cref{lemma4dense}.
\end{Proof}

\begin{lemma}\label{lemma42nbs}
	Let $G$ be a locally bipartite graph with $\delta(G) > 6/11 \cdot \abs{G}$ which does not contain $H_{0}$. Then $G$ does not contain a sparse pair $u, v$ with $uv$ being the missing spoke of an odd wheel.
\end{lemma}

\begin{Proof}
	Suppose that $u, v$ is a sparse pair which is the missing spoke of a $(2k + 1)$-wheel, where we assume that $k$ is minimal. By \cref{lemma45wheel}, $k$ is at least 3. Label the configuration as follows and write $v_{0}$ for $v$ (we consider indices modulo $2k + 1$). Let $C = \set{v, v_{1}, \dotsc, v_{2k}}$.
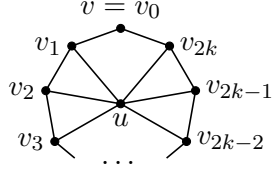
	\begin{figure}[H]
		\centering
		\begin{tikzpicture}
			\foreach \pt in {0,1,2,3,4,5,6,7,8} 
			{
				\tkzDefPoint(\pt*360/9 + 90:1){v_\pt}
			}
			\tkzDefPoint(0,0){u}
			\tkzDefMidPoint(v_3,v_4) \tkzGetPoint{A}
			\tkzDefMidPoint(v_5,v_6) \tkzGetPoint{B}
			\tkzDefMidPoint(v_4,v_5) \tkzGetPoint{C}
			
			\tkzDrawPolySeg(v_6,v_7,v_8,v_0,v_1,v_2,v_3)
			\foreach \pt in {1,2,3,6,7,8} 
			{
				\tkzDrawSegment(u,v_\pt)
			}
			\tkzDrawSegment(v_3,A)
			\tkzDrawSegment(B,v_6)
			
			\tkzDrawPoints(v_0,v_1,v_2,v_3,v_6,v_7,v_8)
			\tkzDrawPoint(u)
			\tkzLabelPoint[above](v_0){$v = v_{0}$}
			\tkzLabelPoint[left](v_1){$v_{1}$}
			\tkzLabelPoint[left](v_2){$v_{2}$}
			\tkzLabelPoint[left](v_3){$v_{3}$}
			\tkzLabelPoint[right](v_6){$v_{2k - 2}$}
			\tkzLabelPoint[right](v_7){$v_{2k - 1}$}
			\tkzLabelPoint[right](v_8){$v_{2k}$}
			\tkzLabelPoint[below](u){$u$}
			\tkzLabelPoint[above](C){$\dotsc$}
		\end{tikzpicture}
	\end{figure}
	Consider a vertex $x$ that is adjacent to $u$. Suppose that $x$ is adjacent to two vertices in $C$ which are not two apart: $x$ is adjacent to $v_{i}$ and $v_{i + r}$ where $r \in \set{1, 3, 4, \dotsc, k}$. Firstly if $r = 1$, then either $G$ contains the $K_{4}$ $uxv_{i}v_{i + 1}$ or $G_{u, v}$ contains an edge (if one of $v_{i}$ or $v_{i + 1}$ is $v$) which contradicts the sparsity of $u, v$. Secondly if $r > 1$ is odd, then $C' = x v_{i} v_{i + 1} \dotsb v_{i + r}$ is an odd cycle which is shorter than $C$. Either $C'$ is in $G_{u}$ (if $v \not \in C'$) contradicting the local bipartiteness of $G$ or we have found a shorter odd cycle than $C$ which satisfies the properties of $C$ (if $v \in C'$). Finally if $r > 2$ is even, then $C' = xv_{i + r}v_{i + r + 1} \dotsb v_{i - 1} v_{i}$ is an odd cycle which is shorter than $C$. Again we either obtain an odd cycle in $G_{u}$ or contradict the minimality of $C$. Hence every neighbour of $u$ has at most two neighbours in $C$.
	
	All vertices have at most $2k$ neighbours in $C$ as otherwise $G$ contains a $(2k + 1)$-wheel. Hence,
	\begin{align*}
		(2k + 1) \delta(G) & \leqslant e(G, C) \leqslant 2 d(u) + 2k(\abs{G} - d(u)) \\
		& = 2k \abs{G} - (2k - 2) d(u) \leqslant 2k \abs{G} - (2k - 2) \delta(G),
	\end{align*}
	so
	\begin{equation*}
		\frac{6}{11} < \frac{\delta(G)}{\abs{G}} \leqslant \frac{2k}{4k - 1},
	\end{equation*}
	which implies that $k < 3$, a contradiction.
\end{Proof}

\subsection{The proof of \texorpdfstring{\cref{main4localbip}}{Theorem 4}}\label{sec:proof4localbip}

Here we will prove \cref{main4localbip} which is restated below for convenience. The argument proceeds as follows. We start with an $H_{0}$-free, locally bipartite graph $G$ with $\delta(G) > 6/11 \cdot \abs{G}$ and wish to show that $G$ is 3-colourable. We may assume that $G$ is edge-maximal (so if it is 3-colourable, then it will in fact be complete tripartite). Edge-maximality and the previous two sections will allow us to classify in what configurations sparse pairs arise.

We then take $I$ to be a largest independent set in $G$. It is enough to show that all edges between $I$ and $G \setminus I$ are present as then $G \setminus I$ is bipartite (since $G$ is locally bipartite) and so $G$ is 3-colourable. Recall \cref{lemma4DuI} which provides information on how $I$ interacts with the rest of $G$. It says that if $u \in I$ and $v \not \in I$ are not adjacent, then the pair $u, v$ is sparse. In fact, we can extract more. The following definition will be helpful.

\begin{definition}[quasidense]
	A pair of vertices $u, v$ is \defn{quasidense} if there is a sequence of vertices $u = d_{1}, d_{2}, \dotsc, d_{k} = v$ such that all pairs $d_{i}, d_{i + 1}$ are dense \textnormal{(}$i = 1, 2, \dotsc, k - 1$\textnormal{)}.
\end{definition}

\Cref{lemma4DuI} immediately implies that if the pair $u, v$ is quasidense and $u \in I$, then $v \in I$ also. So if there are vertices $u \in I$, $v \not \in I$ with $u$ not adjacent to $v$, then the pair $u, v$ is sparse and, furthermore, not quasidense. This will contradict our classification of the configurations in which sparse pairs appear.

\mainforlocalbip*

\begin{Proof}
	Let $G$ be an $H_{0}$-free, locally bipartite graph that satisfies $\delta(G) > 6/11 \cdot \abs{G}$. We are required to show that $G$ is 3-colourable. We may assume that $G$ is edge-maximal: for any sparse pair $u, v$ of $G$, the addition of edge $uv$ to $G$ introduces an odd wheel or a copy of $H_{0}$. By \cref{lemma4H}, the addition of $uv$ to $G$ introduces an odd wheel, a copy of $H_{2}^{+}$, or a copy of $\overline{C}_{7}$ (note that $G$ itself does not contain these).
	
	Firstly, if the addition of $uv$ introduces an odd wheel, then, by \cref{lemma42nbs}, $uv$ must be a rim of that wheel -- this case is depicted in \cref{fig:Woddrimuv} below. Secondly, if the addition of $uv$ introduces a copy of $\overline{C}_{7}$, then that copy of $\overline{C}_{7}$ less the edge $uv$ must not contain $H_{0}$ -- this case is depicted in \cref{fig:C7uv34}. Finally, if the addition of $uv$ introduces a copy of $H_{2}^{+}$, then that copy of $H_{2}^{+}$ less the edge $uv$ must not contain $H_{0}$ -- this case is depicted in \cref{fig:H21uv01,fig:H21uv12,fig:H21uv23,fig:H21uv34,fig:H21uv13} below.
	
	Thus, in $G$, any sparse pair $u, v$ must appear in one of the following configurations (with the labels of $u$ and $v$ possibly swapped).
	
	\begin{figure}[H]
		\centering
		\begin{subfigure}{.24\textwidth}
			\centering
			\begin{tikzpicture}
				\foreach \pt in {0,1,...,10} 
				{
					\tkzDefPoint(\pt*360/11 + 90:1){v_\pt}
				}
				\tkzDefPoint(0,0){u}
				\tkzDefMidPoint(v_1,v_2) \tkzGetPoint{A}
				\tkzDefMidPoint(v_9,v_10) \tkzGetPoint{B}
				\tkzDefMidPoint(v_1,v_10) \tkzGetPoint{C}
				
				\tkzDrawPolySeg(v_2,v_3,v_4,v_5)
				\tkzDrawPolySeg(v_6,v_7,v_8,v_9)
				\foreach \pt in {2,3,...,9} 
				{
					\tkzDrawSegment(u,v_\pt)
				}
				\tkzDrawSegment(v_2,A)
				\tkzDrawSegment(B,v_9)
				
				\tkzDrawPoints(v_2,v_...,v_9)
				\tkzDrawPoint(u)
				\tkzLabelPoint[below](v_5){$u$}
				\tkzLabelPoint[below](v_6){$v$}
				\tkzLabelPoint[below](C){$\dotsc$}
			\end{tikzpicture}
			\subcaption{}\label{fig:Woddrimuv}
		\end{subfigure}
		\begin{subfigure}{.24\textwidth}
			\centering
			\begin{tikzpicture}
				\foreach \pt in {0,1,...,6} 
				{
					\tkzDefPoint(\pt*360/7 + 90:1){v_\pt}
				} 
				\tkzDrawPolySeg(v_4,v_5,v_6,v_0,v_1,v_2,v_3)
				\tkzDrawPolySeg(v_0,v_2,v_4,v_6,v_1,v_3,v_5,v_0)
				\tkzDrawPoints(v_0,v_...,v_6)
				\tkzLabelPoint[below](v_3){$u$}
				\tkzLabelPoint[below](v_4){$v$}
			\end{tikzpicture}
			\subcaption{}\label{fig:C7uv34}
		\end{subfigure}
		\begin{subfigure}{.24\textwidth}
			\centering
			\begin{tikzpicture}
				\foreach \pt in {0,1,...,6} 
				{
					\tkzDefPoint(\pt*360/7 + 90:1){v_\pt}
				}
				\tkzDefPoint(0,0){u}
				\tkzDrawPolySeg(v_1,v_2,v_3,v_4,v_5,v_6,v_0)
				\tkzDrawPolySeg(v_1,v_3,v_5,v_0, v_2,v_4,v_6)
				\tkzDrawSegments(u,v_0 u,v_2 u,v_5)
				\tkzDrawPoints(v_0,v_...,v_6,u)
				\tkzLabelPoint[above](v_0){$u$}
				\tkzLabelPoint[above left](v_1){$v$}
			\end{tikzpicture}
			\subcaption{}\label{fig:H21uv01}
		\end{subfigure}
		\begin{subfigure}{.24\textwidth}
			\centering
			\begin{tikzpicture}
				\foreach \pt in {0,1,...,6} 
				{
					\tkzDefPoint(\pt*360/7 + 90:1){v_\pt}
				}
				\tkzDefPoint(0,0){u}
				\tkzDrawPolySeg(v_2,v_3,v_4,v_5,v_6,v_0,v_1)
				\tkzDrawPolySeg(v_1,v_3,v_5,v_0, v_2,v_4,v_6)
				\tkzDrawSegments(u,v_0 u,v_2 u,v_5)
				\tkzDrawPoints(v_0,v_...,v_6,u)
				\tkzLabelPoint[above left](v_1){$u$}
				\tkzLabelPoint[left](v_2){$v$}
			\end{tikzpicture}
			\subcaption{}\label{fig:H21uv12}
		\end{subfigure}
		
		\medskip
		
		\begin{subfigure}{.32\textwidth}
			\centering
			\begin{tikzpicture}
				\foreach \pt in {0,1,...,6} 
				{
					\tkzDefPoint(\pt*360/7 + 90:1){v_\pt}
				}
				\tkzDefPoint(0,0){u}
				\tkzDrawPolySeg(v_3,v_4,v_5,v_6,v_0,v_1,v_2)
				\tkzDrawPolySeg(v_1,v_3,v_5,v_0, v_2,v_4,v_6)
				\tkzDrawSegments(u,v_0 u,v_2 u,v_5)
				\tkzDrawPoints(v_0,v_...,v_6,u)
				\tkzLabelPoint[left](v_2){$u$}
				\tkzLabelPoint[below](v_3){$v$}
			\end{tikzpicture}
			\subcaption{}\label{fig:H21uv23}
		\end{subfigure}
		\begin{subfigure}{.32\textwidth}
			\centering
			\begin{tikzpicture}
				\foreach \pt in {0,1,...,6} 
				{
					\tkzDefPoint(\pt*360/7 + 90:1){v_\pt}
				}
				\tkzDefPoint(0,0){u}
				\tkzDrawPolySeg(v_4,v_5,v_6,v_0,v_1,v_2,v_3)
				\tkzDrawPolySeg(v_1,v_3,v_5,v_0, v_2,v_4,v_6)
				\tkzDrawSegments(u,v_0 u,v_2 u,v_5)
				\tkzDrawPoints(v_0,v_...,v_6,u)
				\tkzLabelPoint[below](v_3){$u$}
				\tkzLabelPoint[below](v_4){$v$}
			\end{tikzpicture}
			\subcaption{}\label{fig:H21uv34}
		\end{subfigure}
		\begin{subfigure}{.32\textwidth}
			\centering
			\begin{tikzpicture}
				\foreach \pt in {0,1,...,6} 
				{
					\tkzDefPoint(\pt*360/7 + 90:1){v_\pt}
				}
				\tkzDefPoint(0,0){u} 
				\tkzDrawPolySeg(v_0,v_1,v_2,v_3,v_4,v_5,v_6,v_0)
				\tkzDrawPolySeg(v_3,v_5,v_0, v_2,v_4,v_6)
				\tkzDrawSegments(u,v_0 u,v_2 u,v_5)
				\tkzDrawPoints(v_0,v_...,v_6,u)
				\tkzLabelPoint[above left](v_1){$u$}
				\tkzLabelPoint[below](v_3){$v$}
			\end{tikzpicture}
			\subcaption{}\label{fig:H21uv13}
		\end{subfigure}
		\caption{Configurations in which a sparse pair $u, v$ may appear (labels $u$ and $v$ possibly swapped).}\label{fig:cases}
	\end{figure}
	
	Consider a largest independent set in $G$: an independent set $I$ of size $\alpha(G)$. We will now show that all edges between $I$ and $G \setminus I$ are present. Fix a vertex $u \in I$ and let $v$ be any other vertex which is not adjacent to $u$. It suffices to show that $v \in I$. If the pair $u, v$ is (quasi)dense, then $v \in I$ so we may assume that $u, v$ is sparse (and not quasidense). Thus $u, v$ appears in one of the configurations given in \cref{fig:cases} (with labels $u$ and $v$ possibly swapped). However, the pair $u, v$ is quasidense in \cref{fig:Woddrimuv,fig:C7uv34,fig:H21uv01,fig:H21uv34,fig:H21uv13}. Hence we may assume that $u, v$ appear in one of \cref{fig:H21uv12,fig:H21uv23}. We consider these two configurations together (ignoring the central vertex). For ease we label some more of the vertices as follows.
	
	
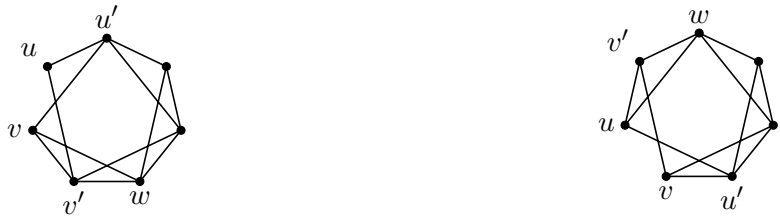
\begin{figure}[H]
		\centering
		\begin{subfigure}{.45\textwidth}
			\centering
			\begin{tikzpicture}
				\foreach \pt in {0,1,...,6} 
				{
					\tkzDefPoint(\pt*360/7 + 90:1){v_\pt}
				} 
				\tkzDrawPolySeg(v_2,v_3,v_4,v_5,v_6,v_0,v_1)
				\tkzDrawPolySeg(v_1,v_3,v_5,v_0, v_2,v_4,v_6)
				\tkzDrawPoints(v_0,v_...,v_6)
				\tkzLabelPoint[above left](v_1){$u$}
				\tkzLabelPoint[left](v_2){$v$}
				\tkzLabelPoint[above](v_0){$u'$}
				\tkzLabelPoint[below](v_3){$v'$}
				\tkzLabelPoint[below](v_4){$w$}
			\end{tikzpicture}
		\end{subfigure}
		\begin{subfigure}{.45\textwidth}
			\centering
			\begin{tikzpicture}
				\foreach \pt in {0,1,...,6} 
				{
					\tkzDefPoint(\pt*360/7 + 90:1){v_\pt}
				} 
				\tkzDrawPolySeg(v_3,v_4,v_5,v_6,v_0,v_1,v_2)
				\tkzDrawPolySeg(v_1,v_3,v_5,v_0, v_2,v_4,v_6)
				\tkzDrawPoints(v_0,v_...,v_6)
				\tkzLabelPoint[left](v_2){$u$}
				\tkzLabelPoint[below](v_3){$v$}
				\tkzLabelPoint[above left](v_1){$v'$}
				\tkzLabelPoint[below](v_4){$u'$}
				\tkzLabelPoint[above](v_0){$w$}
			\end{tikzpicture}
		\end{subfigure}
	\end{figure}
	\addtocounter{figure}{-1}
	In both cases, the pair $u', w$ is dense and so, by \cref{lemma4dense}, the pair $u', v'$ is not dense. However, $u'v'$ is not an edge, as the pair $u, v$ is sparse, and so $u', v'$ is a sparse pair. But then $uu'vv'$ is an induced 4-cycle in which both non-edges are sparse which contradicts \cref{lemma4sparse}.
	
	Thus all edges between $I$ and $G \setminus I$ are present. Let $u \in I$, so $G_{u} = G[V(G) \setminus I]$. But $G$ is locally bipartite, so $G[V(G) \setminus I]$ is bipartite. Using a third colour for the independent set $I$ gives a 3-colouring of $G$.
\end{Proof}

\section{Homomorphism results}\label{sec:hom}

In this section, we will prove \cref{hom2C,hom2H,hom2Hepsilon}, which we restate here for convenience.

\homtoC*

\homtoH*

\homtoHepsilon*

The proof of \cref{hom2C} (which appears in \cref{sec:homtoC}) takes a copy of $\overline{C}_{7}$ in $G$ and builds structure around it, focussing initially on those vertices with four neighbours in the copy of $\overline{C}_{7}$ (of which there are many -- more than $9/11 \cdot \abs{G}$ in fact) and then tacking the rest onto these. The proof of \cref{hom2H} (which appears in \cref{sec:homtoH}) is eminently similar but longer. In place of a copy of $\overline{C}_{7}$ we take a copy of $H_{2}^{+}$ (if one of these is not present, then, by \cref{main4localbip,lemma4H}, $G$ is either 3-colourable or contains $\overline{C}_{7}$ and so we are done by \cref{hom2C} -- note that $\overline{C}_{7}$ is 4-colourable). Around this copy of $H_{2}^{+}$, structure is built in an analogous way to the proof of \cref{hom2C} with the aim of showing that there is a homomorphism $G \to H_{2}^{+}$ (which we believe there is). We will not fully complete this endeavour but will show that there is a homomorphism from $G$ to the following graph which is $H_{2}^{+}$ with four extra edges.

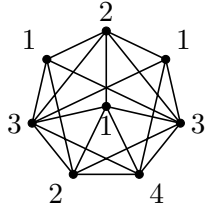
\begin{figure}[H]
	\centering
	\begin{tikzpicture}
		\foreach \pt in {0,1,...,6} 
		{
			\tkzDefPoint(\pt*360/7 + 90:1){v_\pt}
		} 
		\tkzDefPoint(0,0){u}
		\tkzDrawPolySeg(v_0,v_1,v_2,v_3,v_4,v_5,v_6,v_0)
		\tkzDrawPolySeg(v_1,v_3,v_5,v_0, v_2,v_4,v_6)
		\tkzDrawSegments(u,v_0 u,v_2 u,v_5 u,v_3 u,v_4 v_2,v_6 v_1,v_5)
		\tkzDrawPoints(v_0,v_...,v_6)
		\tkzDrawPoint(u)
		\tkzLabelPoint[below](u){1}
		\tkzLabelPoint[above](v_0){2}
		\tkzLabelPoint[above left](v_1){1}
		\tkzLabelPoint[left](v_2){3}
		\tkzLabelPoint[below left](v_3){2}
		\tkzLabelPoint[below right](v_4){4}
		\tkzLabelPoint[right](v_5){3}
		\tkzLabelPoint[above right](v_6){1}
	\end{tikzpicture}
\end{figure}

Thankfully, this graph is 4-colourable (colouring shown in the diagram) and so we have \cref{hom2H}. Finally, the proof of \cref{hom2Hepsilon} (which appears in \cref{sec:homtoHepsilon}) uses all the machinery developed in the proof of \cref{hom2H} and makes use of taking $\varepsilon$ sufficiently small so that the overall structure of $G$ is very similar to that of the weighted $H_{2}^{+}$ shown in \cref{fig:weightings} (which has minimum degree $5/9$).

\subsection{Proof of \texorpdfstring{\cref{hom2C}}{Theorem 6}}\label{sec:homtoC}

In this subsection we prove \cref{hom2C}. Fix a locally bipartite graph $G$ with $\delta(G) > 6/11 \cdot \abs{G}$ that contains a copy of $\overline{C}_{7}$ which we label as follows. We will always consider indices modulo seven.
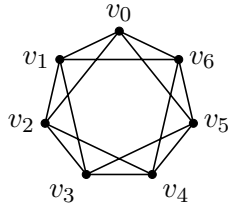
\begin{figure}[H]
	\centering
	\begin{tikzpicture}
		\foreach \pt in {0,1,...,6} 
		{
			\tkzDefPoint(\pt*360/7 + 90:1){v_\pt}
		} 
		\tkzDrawPolySeg(v_0,v_1,v_2,v_3,v_4,v_5,v_6,v_0) 
		\tkzDrawPolySeg(v_0,v_2,v_4,v_6,v_1,v_3,v_5,v_0)
		\tkzDrawPoints(v_0,v_...,v_6)
		\tkzLabelPoint[above](v_0){$v_{0}$}
		\tkzLabelPoint[left](v_1){$v_{1}$}
		\tkzLabelPoint[left](v_2){$v_{2}$}
		\tkzLabelPoint[below left](v_3){$v_{3}$}
		\tkzLabelPoint[below right](v_4){$v_{4}$}
		\tkzLabelPoint[right](v_5){$v_{5}$}
		\tkzLabelPoint[right](v_6){$v_{6}$}
	\end{tikzpicture}
\end{figure}
Let
\begin{align*}
	D & = \set{x \in V(G) \colon x \text{ is adjacent to four of } v_{0}, v_{1}, \dotsc, v_{6}}, \\
	R & = \set{x \in V(G) \colon x \text{ is adjacent to at most three of } v_{0}, v_{1}, \dotsc, v_{6}}.
\end{align*}
By \cref{rmk:5nbs}, no vertex is adjacent to five of the $v_{i}$ so $D \cup R$ partitions $V(G)$. More precisely, note that no vertex is adjacent to three consecutive $v_{i}$ (otherwise there is a $K_{4}$) nor to all of $v_{i - 2}, v_{i}, v_{i + 2}$ (otherwise there is a 5-wheel centred at $v_{i}$). In particular, if we let
\begin{equation*}
	D_{i} = \Gamma(v_{i - 2}, v_{i - 1}, v_{i + 1}, v_{i + 2}),
\end{equation*}
then $D_{0} \cup D_{1} \cup \dotsb \cup D_{6}$ partition $D$. We have a simple upper bound on the size of $R$.
\begin{align}
	7 \delta(G) & \leqslant e(\set{v_{0}, v_{1}, \dotsc, v_{6}}, G) \leqslant 4 \abs{D} + 3 \abs{R} = 4 \abs{G} - \abs{R} \nonumber \\
	\Rightarrow \abs{R} & \leqslant 4 \abs{G} - 7 \delta(G). \label{eq:R}
\end{align}
Also note that $D_{i} \cup D_{i + 3}$ is independent for all $i$: if there is an edge $dd'$ inside $D_{i} \cup D_{i + 3}$, then $dv_{i + 1}v_{i + 2}d'$ is a $K_{4}$. In particular, there is a homomorphism $G[D] \to \overline{C}_{7}$. Our aim is to get a handle on $R$.
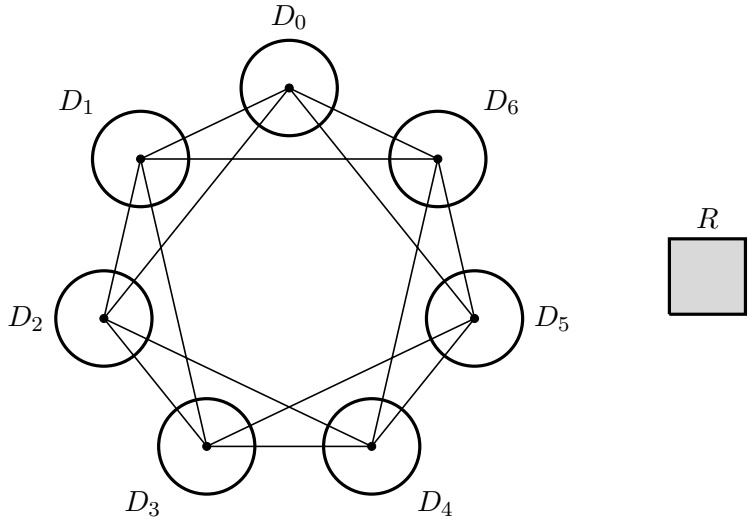
\begin{figure}[H]
	\centering
	\begin{tikzpicture}
		\foreach \pt in {0,1,...,6} 
		{
			\tkzDefPoint(\pt*360/7 + 90:2.5){v_\pt}
			\tkzDefShiftPoint[v_\pt](0:18pt){u_\pt}
			\tkzDrawCircle[black,very thick](v_\pt,u_\pt)
		}
		\tkzDrawPolySeg(v_0,v_1,v_2,v_3,v_4,v_5,v_6,v_0) 
		\tkzDrawPolySeg(v_0,v_2,v_4,v_6,v_1,v_3,v_5,v_0)
		\tkzDrawPoints(v_0,v_...,v_6)
		\tkzDefShiftPoint[v_0](90:18pt){D0}
		\tkzDefShiftPoint[v_1](135:18pt){D1}
		\tkzDefShiftPoint[v_2](180:18pt){D2}
		\tkzDefShiftPoint[v_3](225:18pt){D3}
		\tkzDefShiftPoint[v_4](315:18pt){D4}
		\tkzDefShiftPoint[v_5](0:18pt){D5}
		\tkzDefShiftPoint[v_6](45:18pt){D6}
		\tkzLabelPoint[above](D0){$D_{0}$}
		\tkzLabelPoint[above left](D1){$D_{1}$}
		\tkzLabelPoint[left](D2){$D_{2}$}
		\tkzLabelPoint[below left](D3){$D_{3}$}
		\tkzLabelPoint[below right](D4){$D_{4}$}
		\tkzLabelPoint[right](D5){$D_{5}$}
		\tkzLabelPoint[above right](D6){$D_{6}$}
		
		\tkzDefPoint(5,-0.5){A}
		\tkzDefPoint(5,0.5){B}
		\tkzDefPoint(6,0.5){C}
		\tkzDefPoint(6,-0.5){D}
		\tkzDrawPolySeg[fill=gray!70,very thick](A,B,C,D,A)
		\tkzDefMidPoint(B,C) \tkzGetPoint{R}
		\tkzLabelPoint[above](R){$R$}
	\end{tikzpicture}
\end{figure}
We will use the following lemma frequently.

\begin{lemma}\label{hom:lemma}
	Let $X \subset V(G)$ be a set of four vertices. Either there is $x \in R$ adjacent to all of $X$ or there is $x \in D$ with at least three neighbours in $X$.
\end{lemma}

\begin{Proof}
	Since $\delta(G) > 6/11 \cdot \abs{G}$, we have $4 \delta(G) > 6 \abs{G} - 7 \delta(G)$. We will use this inequality multiple times in this section. This together with inequality \eqref{eq:R} gives
	\begin{equation*}
		e(X, G) \geqslant 4 \delta(G) > 6 \abs{G} - 7 \delta(G) \geqslant 2 \abs{G} + \abs{R} = 2 \abs{D} + 3 \abs{R}.
	\end{equation*}
	But $D \cup R$ partition $V(G)$ so either some vertex in $D$ has more than two neighbours in $X$ or some vertex in $R$ has more than three neighbours in $X$.
\end{Proof}

Our first two claims show that the collections of $v_{i}$ to which vertices can be adjacent are similar to the collections of the $D_{i}$ in which vertices can have neighbours.

\begin{claim}\label{hom:3consec}
	For all $i$, no vertex has a neighbour in each of $D_{i - 1}$, $D_{i}$, $D_{i + 1}$.
\end{claim}

\begin{Proof}
	If not, without loss of generality we may choose $d_{6}, d_{0}, d_{1}$ in $D_{6}, D_{0}, D_{1}$ respectively with common neighbour $u$ such that $e(\set{d_{6}, d_{0}}) + e(\set{d_{0}, d_{1}})$ is maximal. We now apply \cref{hom:lemma} to $\set{u, d_{6}, d_{0}, d_{1}}$.
	
	Suppose some $x$ is adjacent to all of $u, d_{6}, d_{0}, d_{1}$. Now apply \cref{hom:lemma} to $X = \set{u, x, d_{6}, d_{1}}$: as $uxd_{6}$ and $uxd_{1}$ are triangles, no vertex is adjacent to all of $X$ and furthermore, any vertex with three neighbours in $X$ must be adjacent to both $d_{6}$ and $d_{1}$. In particular, some $d' \in D$ is adjacent to $d_{6}$, $d_{1}$ and to one of $u$, $x$. But then $d' \in D_{0}$, so, in our choice of $d_{6}, d_{0}, d_{1}, u$ at the start, we could swap $d'$ for $d_{0}$ and $u$ for whichever of $u$ and $x$ is adjacent to $d'$. This contradicts the maximality unless $d_{0}$ is adjacent to both $d_{6}$ and $d_{1}$. But then $d_{0}d_{1}ux$ is a $K_{4}$.
	
	Hence, in fact, there is some $d \in D$ adjacent to three of $u, d_{6}, d_{0}, d_{1}$. No vertex in $D$ is adjacent to all of $d_{6}, d_{0}, d_{1}$, so $d$ is adjacent to $u$. By symmetry we may assume $d$ is adjacent to $d_{6}$. If $d$ is adjacent to $d_{0}$ as well, then $d$ is adjacent to both $v_{6}$ and $v_{0}$. But then $ud_{6}v_{0}v_{6}d_{0}$ is an odd circuit in $G_{d}$.
	
	Thus $d$ is adjacent to $u, d_{6}$ and $d_{1}$. But then, $d \in D_{0}$ and so $ud_{6}v_{1}v_{6}d_{1}$ is an odd circuit in $G_{d}$.
\end{Proof}

\begin{claim}\label{hom:3spaced}
	For all $i$, no vertex has a neighbour in each of $D_{i - 2}$, $D_{i}$, $D_{i + 2}$.
\end{claim}

\begin{Proof}
	If not, without loss of generality we may choose $d_{5}, d_{0}, d_{2}$ in $D_{5}, D_{0}, D_{2}$ respectively with common neighbour $u$ such that $e(\set{d_{5}, d_{0}}) + e(\set{d_{0}, d_{2}})$ is maximal. No vertex in $D$ has a neighbour in each of $D_{5}$, $D_{0}$ and $D_{2}$, so $u \in R$. Apply \cref{hom:lemma} to $\set{u, d_{5}, d_{0}, d_{2}}$.
	
	Suppose some $x$ is adjacent to all of $u, d_{5}, d_{0}, d_{2}$. Now apply \cref{hom:lemma} to $\set{u, x, d_{0}, d_{2}}$: as $uxd_{0}$ and $uxd_{2}$ are triangles, there must be some $d' \in D$ is adjacent to $d_{0}$, $d_{2}$ and one of $u$, $x$. But then $d' \in D_{1}$ and so one of $u$, $x$ has a neighbour in each of $D_{0}$, $D_{1}$ and $D_{2}$ which contradicts \cref{hom:3consec}.
	
	Hence, there is some $d \in D$ adjacent to three of $u, d_{5}, d_{0}, d_{2}$. No vertex in $D$ is adjacent to all of $d_{5}, d_{0}, d_{2}$ so $d$ is adjacent to $u$. By symmetry, we may assume $d$ is adjacent to $d_{2}$. If $d$ is adjacent to $d_{0}$ as well, then $d \in D_{1}$, so $u$ has a neighbour in each of $D_{0}, D_{1}, D_{2}$ contradicting \cref{hom:3consec}.
	
	Thus $d$ is adjacent to $u, d_{5}$ and $d_{2}$, so $d \in D_{0} \cup D_{3} \cup D_{4}$. If $d \in D_{0}$, then $ud_{5}v_{6}v_{1}d_{2}$ is an odd circuit in $G_{d}$. Hence, we may assume by symmetry that $d \in D_{3}$. Write $d_{3}$ for $d$.
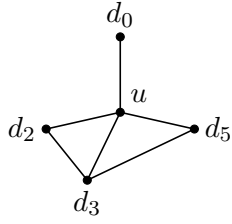
	\begin{figure}[H]
		\centering
		\begin{tikzpicture}
			\foreach \pt in {0,2,3,5} 
			{
				\tkzDefPoint(\pt*360/7 + 90:1){d_\pt}
			}
			\tkzDefPoint(0,0){u}
			\tkzDrawSegments(u,d_0 u,d_2 u,d_3 u,d_5 d_2,d_3 d_3,d_5)
			\tkzDrawPoints(u,d_0,d_2,d_3,d_5)
			\tkzLabelPoint[above right](u){$u$}
			\tkzLabelPoint[above](d_0){$d_{0}$}
			\tkzLabelPoint[left](d_2){$d_{2}$}
			\tkzLabelPoint[below](d_3){$d_{3}$}
			\tkzLabelPoint[right](d_5){$d_{5}$}
		\end{tikzpicture}
	\end{figure}
	We now show there is some $d'_{0} \in D_{0}$ adjacent to both $u$ and $d_{5}$. Apply \cref{hom:lemma} to $\set{u, d_{5}, v_{6}, d_{0}}$: by \cref{hom:3consec}, no vertex is adjacent to all of $d_{5}, v_{6}, d_{0}$ so there is $d'' \in D$ adjacent to $u$ and to two of $d_{5}, v_{6}, d_{0}$.
	\begin{itemize}[noitemsep]
		\item If $d''$ is adjacent to $d_{5}$ and $d_{0}$, then $d'' \in D_{6}$, so $u$ has a neighbour in each of $D_{5}, D_{6}, D_{0}$, contrary to \cref{hom:3consec}.
		\item If $d''$ is adjacent to $d_{5}$ and $v_{6}$, then $d'' \in D_{4} \cup D_{0}$. But if $d'' \in D_{4}$, then $u$ has a neighbour in each of $D_{2}, D_{3}, D_{4}$, contrary to \cref{hom:3consec}, so $d'' \in D_{0}$. We may take $d'_{0} = d''$.
		\item If $d''$ is adjacent to $v_{6}$ and $d_{0}$, then $d'' \in D_{5} \cup D_{1}$. But if $d'' \in D_{1}$, then $u$ has a neighbour in each of $D_{0}, D_{1}, D_{2}$, contrary to \cref{hom:3consec}, so $d'' \in D_{5}$. By maximality at the start we must have $d_{5}$ adjacent to $d_{0}$. We may take $d'_{0} = d_{0}$.
	\end{itemize}
	Thus there is some $d'_{0} \in D_{0}$ adjacent to both $u$ and $d_{5}$. But then $d'_{0}ud_{3}v_{4}v_{6}$ is an odd circuit in $G_{d_{5}}$.
\end{Proof}

From the previous two claims it follows that for every vertex $v$ there is an $i$ such that
\begin{equation*}
	\Gamma(v) \cap D \subset \Gamma(v_{i}) \cap D =  D_{i - 2} \cup D_{i - 1} \cup D_{i + 1} \cup D_{i + 2}.
\end{equation*}
For $i = 0, 1, \dotsc, 6$ choose
\begin{equation*}
	R_{i} \subset \set{v \in R \colon \Gamma(v) \cap D \subset \Gamma(v_{i}) \cap D},
\end{equation*}
so that $R_{0} \cup R_{1} \cup \dotsb \cup R_{6}$ is a partition of $R$. There may be some flexibility in the choice of the $R_{i}$ (e.g.\ if $\Gamma(v) \cap D \subset D_{0} \cup D_{3} \cup D_{4}$, then we could take $v$ in $R_{2}$ or $R_{5}$) -- we will make use of this later. For now we just take any arbitrary choice. Note, by definition, that
\begin{equation*}
	e(R_{i}, D_{i} \cup D_{i - 3} \cup D_{i + 3}) = 0.
\end{equation*}
For each $i$, let $T_{i} = D_{i} \cup R_{i}$ -- note that these partition $V(G)$. We can give a lower bound for the size of $T_{i}$. For each $i$,
\begin{align*}
	d(v_{i - 1}) + d(v_{i + 1}) & = \abs{D_{i}} + \abs{D} + \abs{R \cap \Gamma(v_{i - 1})} + \abs{R \cap \Gamma(v_{i + 1})} \\
	& = \abs{D_{i}} + \abs{D} + \abs{R \cap (\Gamma(v_{i - 1}) \cup \Gamma(v_{i + 1}))} + \abs{R \cap \Gamma(v_{i - 1}, v_{i + 1})} \\
	& \leqslant \abs{D_{i}} + \abs{D} + \abs{R} + \abs{R_{i}} = \abs{G} + \abs{T_{i}},
\end{align*}
so
\begin{equation}\label{eq:DR}
	\abs{T_{i}} \geqslant 2 \delta(G) - \abs{G}.
\end{equation}

We will eventually show that $T_{i} \cup T_{i + 3}$ is independent for all $i$ and so the map sending all vertices in $T_{i}$ to $v_{i}$ is a homomorphism from $G$ to $\overline{C}_{7}$.

\begin{claim}\label{hom:DR}
	Every $d \in D_{i}$ and $u \in T_{i  + 1}$ have a common neighbour in $D$. Similarly, every $d \in D_{i}$ and $u \in T_{i - 1}$ have a common neighbour in $D$.
\end{claim}

\begin{Proof}
	By symmetry it suffices to prove this for $d_{0} \in D_{0}$ and $u \in T_{1}$. As $e(D_{0} \cup T_{1}, D_{4}) = 0$, we have $\Gamma(u) \cup \Gamma(d_{0}) \subset V(G) \setminus D_{4}$ so
	\begin{align*}
		\abs{\Gamma(d_{0}) \cap \Gamma(u)} & = d(d_{0}) + d(u) - \abs{\Gamma(d_{0}) \cup \Gamma(u)} \geqslant 2 \delta(G) + \abs{D_{4}} - \abs{G} \\
		& \geqslant 4 \delta(G) - 2 \abs{G} - \abs{R_{4}} > 4 \abs{G} - 7 \delta(G) - \abs{R_{4}} \geqslant \abs{R} - \abs{R_{4}} \\
		& \geqslant \abs{\Gamma(d_{0}) \cap R},
	\end{align*}
	where we have used inequality \eqref{eq:DR}, $\delta(G) > 6/11 \cdot \abs{G}$, inequality \eqref{eq:R} and $e(d_{0}, R_{4}) = 0$ respectively for the final four inequalities. In particular, $d_{0}$ and $u$ have a common neighbour $d_{u} \in \Gamma(d_{0}) \cap D$. 
\end{Proof}

\begin{claim}\label{hom:indep1}
	For all $i$ and $d \in D_{i}$, the sets $\Gamma(d) \cap T_{i - 1}$ and $\Gamma(d) \cap T_{i + 1}$ are independent. 
\end{claim}

\begin{Proof}
	By symmetry it suffices to prove this for $\Gamma(d_{0}) \cap T_{1}$ where $d_{0} \in D_{0}$. Suppose that $\Gamma(d_{0}) \cap T_{1}$ contains an edge $uv$. By \cref{hom:DR}, $d_{0}$ and $u$ have a common neighbour $d_{u} \in D$. As $d_{u}$ is adjacent to $u \in T_{1}$, it must also be adjacent to $v_{1}$. Similarly there is $d_{v} \in D$ adjacent to $d_{0}, v, v_{1}$. But then $v_{1}d_{u}uvd_{v}$ is an odd circuit in $G_{d_{0}}$.
\end{Proof}

\begin{claim}\label{hom:DRindep}
	For all $i$, $T_{i}$ is independent.
\end{claim}

\begin{Proof}
	By symmetry it suffices to prove this for $i = 0$. Suppose that $uv$ is an edge in $T_{0}$. We already have that $D_{0}$ is independent and $e(D_{0}, R_{0}) = 0$ so $u, v \in R_{0}$. We may choose the edge $uv$ in $R_{0}$ so that $e(\set{u, v, v_{1}, v_{6}})$ is maximal. Apply \cref{hom:lemma} to $\set{u, v, v_{1}, v_{6}}$.
	
	Suppose some $x \in R$ is adjacent to all of $u, v, v_{1}, v_{6}$. Then $x \in R_{0}$. By maximality of $e(\set{u, v, v_{1}, v_{6}})$, we must have had $u$ and $v$ adjacent to both $v_{1}$, $v_{6}$ and so $uvv_{1}v_{6}$ is a $K_{4}$.
	
	Hence, there is some $d \in D$ with three neighbours amongst $u, v, v_{1}, v_{6}$. Either $d$ is adjacent to both $u$, $v$ or to both $v_{1}$, $v_{6}$. If $d$ is adjacent to both $v_{1}$, $v_{6}$, then $d \in D_{0}$. But then $d$ is adjacent to neither $u$ nor $v$, as $e(D_{0}, R_{0}) = 0$. Hence $d$ is adjacent to both $u$, $v$. By \cref{hom:indep1}, $d \not \in D_{1} \cup D_{6}$ so $d \in D_{2} \cup D_{5}$. By symmetry, we may assume that $d = d_{2}  \in D_{2}$.
	
	Apply \cref{hom:lemma} to $\set{u, v, d_{2}, v_{5}}$: $d_{2}uv$ is a triangle so there is $d' \in D$ adjacent to $v_{5}$ and to two of $u, v, d_{2}$. As $d'$ is adjacent to $v_{5}$ and at least one of $u, v \in R_{0}$, $d'$ is in $D_{6}$. But then $d'$ is not adjacent to $d_{2}$ and so is to both $u$ and $v$. However, edge $uv$ lies in $\Gamma(d') \cap R_{0}$, contrary to \cref{hom:indep1}.
\end{Proof}

This shows that there is a homomorphism $G \to K_{7}$ and so $G$ is 7-colourable. Before proceeding it will help to give structure to $G_{d}$ for each $d \in D$. 

\begin{claim}\label{claim:Gd}
	For any $i \in \set{0, 1, 2, \dotsc, 6}$ and any $d \in D_{i}$, $G_{d}$ is connected bipartite. Furthermore, there is a bipartition of $G_{d}$ into two vertex classes $A_{d}$, $B_{d}$ which satisfy $(T_{i - 1} \cup D_{i + 2}) \cap \Gamma(d) \subset A_{d}$, $(T_{i + 1} \cup D_{i - 2}) \cap \Gamma(d) \subset B_{d}$ and at least one of $R_{i + 2} \cap \Gamma(d) \subset A_{d}$, $R_{i - 2} \cap \Gamma(d) \subset B_{d}$ occurs.
\end{claim}

\begin{Proof}
	We may assume $i = 0$. Fix $d \in D_{0}$ and define for $j = 5, 6, 1, 2$,
	\begin{align*}
		D_{j}^{d} & = D_{j} \cap \Gamma(d), \\
		R_{j}^{d} & = R_{j} \cap \Gamma(d), \\
		T_{j}^{d} & = T_{j} \cap \Gamma(d) = D_{j}^{d} \cup R_{j}^{d},
	\end{align*}
	and note that the $T_{j}^{d}$ partition $V(G_{d})$. Since $v_{6} \in G_{d}$, we can define
	\begin{align*}
		A_{d} & = \set{x \in G_{d} \colon \operatorname{dist}_{G_{d}}(x, v_{6}) \text{ is even}}, \\
		B_{d} & = \set{x \in G_{d} \colon \operatorname{dist}_{G_{d}}(x, v_{6}) \text{ is odd}}.
	\end{align*}
	$G$ is locally bipartite, so $G_{d}$ is bipartite and so $A_{d}$ and $B_{d}$ are independent sets. Now
	\begin{itemize}[noitemsep]
		\item $v_{6} \in A_{d}$, $v_{1} \in B_{d}$.
		\item $v_{6}$ is adjacent to all of $D_{5}^{d} \cup D_{1}^{d}$, so $D_{5}^{d} \cup D_{1}^{d} \subset B_{d}$.
		\item $v_{1}$ is adjacent to all of $D_{6}^{d} \cup D_{2}^{d}$, so $D_{6}^{d} \cup D_{2}^{d} \subset A_{d}$.
	\end{itemize}
	We next show that $R_{6}^{d} \subset A_{d}$. Let $x \in R_{6}^{d}$: by \cref{hom:DR}, $d$ and $x$ have a common neighbour $d' \in D$. As $d \in D_{0}$ and $x \in R_{6}$, $d'$ must be in $D_{1} \cup D_{5}$. Hence, $d' \in B_{d}$ and so $x \in A_{d}$. Similarly $R_{1}^{d} \subset B_{d}$.
	
	We now show that at least one of $R_{2}^{d} \subset A_{d}$, $R_{5}^{d} \subset B_{d}$ occurs. If not, then there is $u \in R_{2}^{d} \setminus A_{d}$ and $v \in R_{5}^{d} \setminus B_{d}$. Focus on $u$: $u \not \in A_{d}$ so $\Gamma_{G_{d}}(u) \subset V(G_{d}) - B_{d} \subset T_{6}^{d} \cup T_{2}^{d} \cup R_{5}^{d}$. But $u \in R_{2}$, the set $T_{2}$ is independent and $e(R_{2}, D_{6}) = 0$, so, in fact,
	\begin{equation*}
		\Gamma_{G_{d}}(u) \subset R_{5}^{d} \cup R_{6}^{d}.
	\end{equation*}
	Similarly,
	\begin{equation*}
		\Gamma_{G_{d}}(v) \subset R_{1}^{d} \cup R_{2}^{d}.
	\end{equation*}
	In particular,
	\begin{equation*}
		\abs{\Gamma_{G_{d}}(u)} + \abs{\Gamma_{G_{d}}(v)} \leqslant \lvert R_{1}^{d} \rvert + \lvert R_{2}^{d} \rvert + \lvert R_{5}^{d} \rvert + \lvert R_{6}^{d} \rvert \leqslant \abs{R}.
	\end{equation*}
	But then, by inequality \eqref{eq:R},
	\begin{equation*}
		4 \abs{G} - 7 \delta(G) \geqslant \abs{R} \geqslant d(d, u) + d(d, v) \geqslant 4 \delta(G) - 2\abs{G},
	\end{equation*}
	which contradicts $\delta(G) > 6/11 \cdot \abs{G}$.
	
	Finally we need to show that $G_{d}$ is connected. We will do this by showing $A_{d} \cup B_{d} = V(G_{d})$. We already have $T_{1}^{d} \cup T_{6}^{d} \cup D_{2}^{d} \cup D_{5}^{d} \subset A_{d} \cup B_{d}$ and at least one of $R_{2}^{d}$, $R_{5}^{d}$ is a subset of $A_{d} \cup B_{d}$ -- we need only show that the other one is too. By symmetry, we may assume $R_{2}^{d} \subset A_{d} \cup B_{d}$. Fix $x \in R_{5}^{d}$. Now $d_{G_{d}}(x) = d(x, d) \geqslant 2 \delta(G) - \abs{G} > 0$, so $x$ has some neighbour in $G_{d}$. But $R_{5}$ is an independent set, so $x$ has a neighbour in $A_{d} \cup B_{d}$ and so $x \in A_{d} \cup B_{d}$.
\end{Proof}

We are finally in a position to show that there is a homomorphism $G \to \overline{C}_{7}$. It is here that we will make use of the flexibility in the choice of the $R_{i}$.

\begin{claim}\label{claim:RR}
	It is possible to choose the $R_{j}$ so that the sets $T_{i} \cup T_{i + 3}$ are all independent.
\end{claim}

\begin{Proof}
	Note that $T_{i}$, $T_{i + 3}$, $D_{i} \cup D_{i + 3}$ are all independent and $e(D_{i}, R_{i + 3}) = e(D_{i + 3}, R_{i}) = 0$ so it suffices to show that it is possible to ensure $e(R_{i}, R_{i + 3}) = 0$ for all $i$. We choose the $R_{i}$ so that
	\begin{equation*}
		S = \sum_{i = 0}^{6} e(R_{i}, R_{i + 3})
	\end{equation*}
	is minimal. Suppose that $S$ is not zero: by symmetry, we may assume there is some $u \in R_{2}$, $v \in R_{6}$ with $u$ adjacent to $v$. Apply \cref{hom:lemma} to $\set{u, v, v_{0}, v_{1}}$. Note that any common neighbour of $v_{0}$, $v_{1}$ is in $T_{2} \cup T_{6}$ so is adjacent to at most one of $u$, $v$. Moreover, any common neighbour of $v_{0}$, $v_{1}$ which lies in $D$ is in $D_{2} \cup D_{6}$ so is adjacent to neither $u$ nor $v$. Hence there is $d \in D$ which is adjacent to both $u$, $v$ and to one of $v_{0}$, $v_{1}$. By symmetry, we may assume $d$ is adjacent to $v_{1}$ and so $d \in D_{0}$. That is, there is at least one $d \in D_{0}$ adjacent to both $u$ and $v$.
	
	For any $d \in D_{0} \cap \Gamma(u, v)$, consider the bipartition of $G_{d}$ given by the previous claim:
	\begin{itemize}[noitemsep]
		\item $(T_{6} \cup D_{2}) \cap \Gamma(d) \subset A_{d}$.
		\item $(T_{1} \cup D_{5}) \cap \Gamma(d) \subset B_{d}$.
		\item At least one of $R_{2} \cap \Gamma(d) \subset A_{d}$, $R_{5} \cap \Gamma(d) \subset B_{d}$ occurs.
	\end{itemize}
	Now $v \in R_{6}$ is adjacent to $d$, so $v \in A_{d}$. Also $u$ is adjacent to $v$ and $d$, so $u \in B_{d}$. But $u \in R_{2} \cap \Gamma(d)$ so $R_{5} \cap \Gamma(d) \subset B_{d}$ occurs. Now $\Gamma_{G_{d}}(u) \subset A_{d} \subset (T_{6} \cup T_{2}) \cap \Gamma(d)$. But $u \in R_{2}$, the set $T_{2}$ is independent and $e(R_{2}, D_{6}) = 0$, so
	\begin{equation}\label{eq:homcont}
		\Gamma_{G_{d}}(u) \subset R_{6} \cap \Gamma(d).
	\end{equation}
	Note that this holds for any choice of $d \in D_{0} \cap \Gamma(u, v)$.
	
	We first deal with the case where $u$ has some neighbour in $D_{1}$. Pick any $d_{0} \in \Gamma(u, v) \cap D_{0}$, $d_{1} \in \Gamma(u) \cap D_{1}$ and apply \cref{hom:lemma} to $\set{u, v, d_{0}, d_{1}}$.
	\begin{figure}[H]
		\centering
		\begin{tikzpicture}
			\tkzDefPoint(-0.8,0){u}
			\tkzDefPoint(0.8,0){v}
			\tkzDefPoint(-0.5,0.5){d1}
			\tkzDefPoint(0.5,0.5){d0}
			\tkzDrawPolySeg(u,v,d0,u)
			\tkzDrawSegment(u,d1)
			\tkzDrawPoints(u,v,d0,d1)
			\tkzLabelPoint[below](u){$u$}
			\tkzLabelPoint[below](v){$v$}
			\tkzLabelPoint[above](d0){$d_{0}$}
			\tkzLabelPoint[above](d1){$d_{1}$}
		\end{tikzpicture}
	\end{figure}
	Vertices $d_{0}uv$ form a triangle so some $d \in D$ is adjacent to $d_{1}$ and to two of $d_{0}, u, v$. If $d$ is adjacent to $d_{0}$, then $d \in D_{2} \cup D_{6}$, so $d$ is adjacent to neither $u$ nor $v$. Hence $d \in \Gamma(u, v, d_{1}) \cap D$, so $d \in D_{0}$. Thus $d \in \Gamma(u,v) \cap D_{0}$, $d_{1} \in \Gamma(u) \cap D_{1}$ and $d$ is adjacent to $d_{1}$. But then $\Gamma_{G_{d}}(u)$ contains $d_{1} \not \in R_{6}$ contradicting \eqref{eq:homcont}.
	
	We are finally left with the case where $u$ has no neighbours in $D_{1}$. This means we could have put $u$ in $R_{5}$ rather than $R_{2}$ when we chose the $R_{i}$. In particular, by the minimality of $S$,
	\begin{equation*}
		e(u, R_{5}) + e(u, R_{6}) \leqslant e(u, R_{1}) + e(u, R_{2}),
	\end{equation*}
	hence,
	\begin{equation*}
		2\bigl(e(u, R_{5}) + e(u, R_{6})\bigr) \leqslant e(u, R_{1} \cup R_{2} \cup R_{5} \cup R_{6}) \leqslant \abs{R}.
	\end{equation*}
	Pick any $d \in \Gamma(u, v) \cap D_{0}$: as $\Gamma_{G_{d}}(u) \subset R_{6} \cap \Gamma(d)$, we have
	\begin{equation*}
		d_{G_{d}}(u) \leqslant e(u, R_{6}) \leqslant 1/2 \cdot \abs{R}.
	\end{equation*}
	Thus
	\begin{equation*}
		\abs{R} \geqslant 2 d(d, u) \geqslant 4 \delta(G) - 2 \abs{G} > 4 \abs{G} - 7 \delta(G) \geqslant \abs{R},
	\end{equation*}
	where we used $\delta(G) > 6/11 \cdot \abs{G}$ and inequality \eqref{eq:R} for the final two inequalities. This is a contradiction and so $S = 0$, as required.
\end{Proof}

\begin{Proof}[of \cref{hom2C}]
	For $i \in \set{1, \dotsc, 7}$, let $R_i$ be given by \cref{claim:RR} and $T_i = D_i \cup R_i$. Then the $T_i$ partition $V(G)$ and there are no edges within $T_i \cup T_{i + 3}$. Thus the map which sends all vertices in $T_i$ to a vertex $v_i$ is a homomorphism from $G$ to a copy of $\overline{C}_7$.
\end{Proof}

\subsection{Proof of \texorpdfstring{\cref{hom2H}}{Theorem 7}}\label{sec:homtoH}

In this subsection we prove \cref{hom2H}. This has many similarities with the proof of \cref{hom2C}, although not having the full symmetry of $\overline{C}_{7}$ available adds some technicalities.

Fix a locally bipartite graph $G$ with $\delta(G) > 6/11 \cdot \abs{G}$. By \cref{main4localbip,lemma4H}, $G$ is either 3-colourable, contains $\overline{C}_{7}$ or contains $H_{2}^{+}$. In the first two cases we are done (using \cref{hom2C}), so we assume that $G$ does not contain a copy of $\overline{C}_{7}$ but does contain a copy of $H_{2}^{+}$ (and so also a copy of $H_{2}$).

We say \defn{$a_{0}a_{1} \dotsc a_{6}$ is a copy of $H_{2}$ in $G$} to mean that the following configuration appears in $G$. We will continue to use \cref{rmk:5nbs}: no vertex is adjacent to five of the vertices which form a copy of $H_{2}$. We will always consider indices modulo seven.
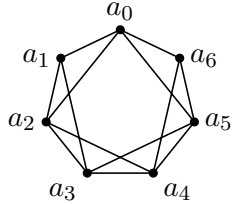
\begin{figure}[H]
	\centering
	\begin{tikzpicture}
		\foreach \pt in {0,1,...,6} 
		{
			\tkzDefPoint(\pt*360/7 + 90:1){v_\pt}
		} 
		\tkzDrawPolySeg(v_0,v_1,v_2,v_3,v_4,v_5,v_6,v_0)
		\tkzDrawPolySeg(v_1,v_3,v_5,v_0, v_2,v_4,v_6)
		\tkzDrawPoints(v_0,v_...,v_6)
		\tkzLabelPoint[above](v_0){$a_{0}$}
		\tkzLabelPoint[left](v_1){$a_{1}$}
		\tkzLabelPoint[left](v_2){$a_{2}$}
		\tkzLabelPoint[below left](v_3){$a_{3}$}
		\tkzLabelPoint[below right](v_4){$a_{4}$}
		\tkzLabelPoint[right](v_5){$a_{5}$}
		\tkzLabelPoint[right](v_6){$a_{6}$}
	\end{tikzpicture}
\end{figure}
Our first few claims will nail down to which $a_{i}$ other vertices may be adjacent. This will eventually allow us to define the sets $D_{i}$ in a similar way to the proof of \cref{hom2C}.

\begin{claim}\label{claim':H2toH21}
	Let $a_{0}a_{1} \dotsc a_{6}$ be a copy of $H_{2}$ in $G$. Then there is a vertex $u \not \in \set{a_{0}, a_{1}, \dotsc, a_{6}}$ adjacent to $a_{5}$, $a_{0}$ and $a_{2}$ \textnormal{(}i.e.\ any copy of $H_{2}$ `extends' to a copy of $H_{2}^{+}$\textnormal{)}. Furthermore, $\abs{\Gamma(a_{5}, a_{0}, a_{2})} \geqslant 11 \delta(G) - 6 \abs{G}$.
\end{claim}

\begin{Proof}
	Consider a copy of $H_{2}$ in $G$ with vertices $X = \{a_{0}, a_{1}, \dotsc, a_{6}\}$. We assign weight $\omega \colon X \to \ZZN$ as shown in the diagram below (recall this notation from \cref{sec:notation}). For each vertex $v \in G$, let $f(v)$ be the total weight of the neighbours of $v$ in $X$.
	\begin{figure}[H]
		\centering
		\begin{tikzpicture}
			\foreach \pt in {0,1,...,6} 
			{
				\tkzDefPoint(\pt*360/7 + 90:1){v_\pt}
			} 
			\tkzDrawPolySeg(v_0,v_1,v_2,v_3,v_4,v_5,v_6,v_0)
			\tkzDrawPolySeg(v_1,v_3,v_5,v_0, v_2,v_4,v_6)
			\tkzDrawPoints(v_0,v_...,v_6)
			\tkzLabelPoint[above](v_0){$a_{0} \colon 3$}
			\tkzLabelPoint[left](v_1){$a_{1} \colon 1$}
			\tkzLabelPoint[left](v_2){$a_{2} \colon 2$}
			\tkzLabelPoint[below left](v_3){$a_{3} \colon 1$}
			\tkzLabelPoint[below right](v_4){$a_{4} \colon 1$}
			\tkzLabelPoint[right](v_5){$a_{5} \colon 2$}
			\tkzLabelPoint[right](v_6){$a_{6} \colon 1$}
		\end{tikzpicture}
	\end{figure}
	The final part of the proof of \cref{claim:H02H21C7} showed that any vertex with $f(v) \geqslant 7$ is either adjacent to all of $a_{5}$, $a_{0}$, $a_{2}$ or is in a copy of $\overline{C}_{7}$. As $G$ does not contain $\overline{C}_{7}$, $\Gamma(a_{5}, a_{0}, a_{2})$ is exactly the set of vertices with $f(v) \geqslant 7$. However, any vertex in $\Gamma(a_{5}, a_{0}, a_{2})$ is adjacent to no other vertex of $X$ as $G$ is locally bipartite. Hence, $\Gamma(a_{5}, a_{0}, a_{2})$ is exactly the set of vertices with $f(v) = 7$ and all other vertices have $f(v) \leqslant 6$. In particular,
	\begin{equation*}
		11 \delta(G) \leqslant \sum_{v \in G} f(v) \leqslant 7 \abs{\Gamma(a_{5}, a_{0}, a_{2})} + 6[\abs{G} - \abs{\Gamma(a_{5}, a_{0}, a_{2})}] = \abs{\Gamma(a_{5}, a_{0}, a_{2})} + 6 \abs{G}.
	\end{equation*}
	Thus $\abs{\Gamma(a_{5}, a_{0}, a_{2})} \geqslant 11 \delta(G) - 6 \abs{G}$ and so $\Gamma(a_{5}, a_{0}, a_{2})$ is non-empty. As $G$ is locally bipartite, no vertex in $X$ is adjacent to all of $a_{5}$, $a_{0}$, $a_{2}$, so the copy of $H_{2}$ extends to a copy of $H^{+}_{2}$.
\end{Proof}

\begin{claim}\label{claim':6013}
	Let $a_{0}a_{1} \dotsc a_{6}$ be a copy of $H_{2}$ in $G$. Then no vertex is adjacent to all of $a_{6}, a_{0}, a_{1}, a_{3}$.
\end{claim}

\begin{Proof}
	Suppose some vertex $a$ is adjacent to all of $a_{6}, a_{0}, a_{1}, a_{3}$. All vertices are adjacent to at most four of the $a_{i}$ so $a$ cannot be one of the $a_{i}$. Let $A = \set{a, a_{0}, a_{1}, \dotsc, a_{6}}$ and gives weights, $\omega$, to the vertices of $A$ as shown.
	\begin{figure}[H]
		\centering
		\begin{tikzpicture}
			\foreach \pt in {0,1,...,6} 
			{
				\tkzDefPoint(\pt*360/7 + 90:1.5){v_\pt}
			}
			\tkzDefPoint(0,0){v}
			\tkzDrawPolySeg(v_0,v_1,v_2,v_3,v_4,v_5,v_6,v_0)
			\tkzDrawPolySeg(v_1,v_3,v_5,v_0, v_2,v_4,v_6)
			\tkzDrawSegments(v,v_0 v,v_1 v,v_3 v,v_6)
			\tkzDrawPoints(v_0,v_...,v_6,v)
			\tkzLabelPoint[above](v_0){$a_{0} \colon 2$}
			\tkzLabelPoint[left](v_1){$a_{1} \colon 1$}
			\tkzLabelPoint[left](v_2){$a_{2} \colon 1$}
			\tkzLabelPoint[below left](v_3){$a_{3} \colon 2$}
			\tkzLabelPoint[below right](v_4){$a_{4} \colon 1$}
			\tkzLabelPoint[right](v_5){$a_{5} \colon 2$}
			\tkzLabelPoint[right](v_6){$a_{6} \colon 1$}
			\tkzLabelPoint[below right = -3pt](v){$a \colon 1$}
		\end{tikzpicture}
		\caption{The configuration of \cref{claim':6013}.}\label{fig:AnotherCounterexample}
	\end{figure}
	For a vertex $v$, let $f(v)$ be the total weight of the neighbours of $v$ in $A$. Now,
	\begin{equation*}
		\sum_{v \in V(G)} f(v) = \omega(A, G) \geqslant 11 \delta(G) > 6 \abs{G},
	\end{equation*}
	so some vertex $v$ has $f(v) \geqslant 7$. Vertex $v$ is not adjacent to all of $a_{0}, a_{3}, a_{5}$ else $va_{0}a_{6}a_{4}a_{3}$ is an odd circuit in $G_{a_{5}}$ and $v$ is not adjacent to all of $a, a_{1}, a_{2}, a_{4}, a_{6}$ as these form a 5-cycle. Thus $v$ is adjacent to exactly two of $a_{0}, a_{3}, a_{5}$ and at least three of $a, a_{1}, a_{2}, a_{4}, a_{6}$. As $v$ is adjacent to at most four of the $a_{i}$, $v$ must be adjacent to $a$.
	\begin{itemize}[noitemsep]
		\item If $v$ is adjacent to $a_{0}$, then $v$ is not adjacent to $a_{1}$ (else $vaa_{0}a_{1}$ is a $K_{4}$), and $v$ is not adjacent to $a_{6}$ (else $vaa_{6}a_{0}$ is a $K_{4}$). Thus $v$ is adjacent to $a_{2}$ and $a_{4}$. But then $va_{0}a_{1}a_{3}a_{4}$ is an odd circuit in $G_{a_{2}}$.
		\item If $v$ is adjacent to both $a_{3}$, $a_{5}$, then $v$ is not adjacent to $a_{4}$ (else $va_{3}a_{4}a_{5}$ is a $K_{4}$), and $v$ is not adjacent to $a_{1}$ (else $vaa_{1}a_{3}$ is a $K_{4}$). Thus $v$ is adjacent to $a_{2}$ and $a_{6}$. But then $G_{a}$ contains the odd circuit $va_{6}a_{0}a_{1}a_{3}$. \qedhere
	\end{itemize}
\end{Proof}

\begin{claim}\label{claim':601}
	Let $a_{0}a_{1} \dotsc a_{6}$ be a copy of $H_{2}$ in $G$. Then no vertex is adjacent to all of $a_{6}, a_{0}, a_{1}$.
\end{claim}

\begin{Proof}
	Suppose some vertex $a$ is adjacent to all of $a_{6}, a_{0}, a_{1}$. All vertices are adjacent to at most four of the $a_{i}$ so $a$ cannot be one of the $a_{i}$. Let $A = \set{a, a_{0}, a_{1}, \dotsc, a_{6}}$ and gives weights, $\omega$, to the vertices of $A$ as shown.
	
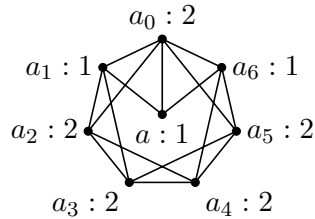
\begin{figure}[H]
		\centering
		\begin{tikzpicture}
			\foreach \pt in {0,1,...,6} 
			{
				\tkzDefPoint(\pt*360/7 + 90:1){v_\pt}
			}
			\tkzDefPoint(0,0){v}
			\tkzDrawPolySeg(v_0,v_1,v_2,v_3,v_4,v_5,v_6,v_0)
			\tkzDrawPolySeg(v_1,v_3,v_5,v_0, v_2,v_4,v_6)
			\tkzDrawSegments(v,v_0 v,v_1 v,v_6)
			\tkzDrawPoints(v_0,v_...,v_6,v)
			\tkzLabelPoint[above](v_0){$a_{0} \colon 2$}
			\tkzLabelPoint[left](v_1){$a_{1} \colon 1$}
			\tkzLabelPoint[left](v_2){$a_{2} \colon 2$}
			\tkzLabelPoint[below left](v_3){$a_{3} \colon 2$}
			\tkzLabelPoint[below right](v_4){$a_{4} \colon 2$}
			\tkzLabelPoint[right](v_5){$a_{5} \colon 2$}
			\tkzLabelPoint[right](v_6){$a_{6} \colon 1$}
			\tkzLabelPoint[below](v){$a \colon 1$}
		\end{tikzpicture}
	\end{figure}
	For a vertex $v$, let $f(v)$ be the total weight of the neighbours of $v$ in $A$. Now,
	\begin{equation*}
		\sum_{v \in V(G)} f(v) = \omega(A, G) \geqslant 13 \delta(G) > 7 \abs{G},
	\end{equation*}
	so some vertex $v$ has $f(v) \geqslant 8$. Vertex $v$ must be adjacent to at least three of $a_{0}, a_{2}, a_{3}, a_{4}, a_{5}$. First suppose that $v$ is adjacent to at least four of $a_{0}, a_{2}, a_{3}, a_{4}, a_{5}$. Vertex $v$ cannot be adjacent to all of $a_{2}, a_{3}, a_{4}$ (else $va_{2}a_{3}a_{4}$ is a $K_{4}$) so $v$ is adjacent to $a_{0}$ and $a_{5}$. Similarly $v$ is adjacent to $a_{2}$. By symmetry, we may assume that $v$ is adjacent to $a_{3}$. But then $va_{0}a_{1}a_{4}a_{3}$ is an odd circuit in $G_{a_{5}}$. Thus $v$ is adjacent to exactly three of $a_{0}, a_{2}, a_{3}, a_{4}, a_{5}$ and so at least two of $a, a_{1}, a_{6}$. As $v$ is adjacent to at most four of the $a_{i}$, $v$ must be adjacent to $a$.
	\begin{itemize}[noitemsep]
		\item If $v$ is adjacent to $a_{0}$, then $v$ cannot be adjacent to $a_{1}$ (else $vaa_{0}a_{1}$ is a $K_{4}$) and $v$ cannot be adjacent to $a_{6}$ (else $vaa_{6}a_{0}$ is a $K_{4}$). Thus $v$ is adjacent to only one of $a, a_{1}, a_{6}$ -- a contradiction.
		\item Otherwise $v$ is not adjacent to $a_{0}$. Certainly $v$ is not adjacent to all of $a_{2}, a_{3}, a_{4}$ (else $va_{2}a_{3}a_{4}$ is a $K_{4}$) so $v$ must be adjacent to $a_{5}$. Similarly $v$ must be adjacent to $a_{2}$. By symmetry, we may assume that $v$ is adjacent to $a_{4}$. Then $v$ is not adjacent to $a_{6}$ (else $va_{4}a_{5}a_{6}$ is a $K_{4}$) so $v$ is adjacent to $a$ and $a_{1}$. Thus $v$ is adjacent to $a_{1}$, $a_{2}, a_{4}, a_{5}$ and $a$. In particular, $v$ is not $a$ nor any of the $a_{i}$ except possibly $a_{3}$. But then $a_{0}a_{1}a_{2}va_{4}a_{5}a_{6}$ is a copy of $H_{2}$ in $G$ and $a$ is adjacent to all of $a_{6}, a_{0}, a_{1}, v$ which contradicts \cref{claim':6013}. \qedhere 
	\end{itemize}
\end{Proof}

\begin{claim}\label{claim':1346}
	Let $a_{0}a_{1} \dotsc a_{6}$ be a copy of $H_{2}$ in $G$. Then no vertex is adjacent to all of $a_{1},a_{3}, a_{4}, a_{6}$.
\end{claim}

\begin{Proof}
	Suppose vertex $a$ is adjacent to all of $a_{1}$, $a_{3}, a_{4}, a_{6}$. All vertices are adjacent to at most four of the $a_{i}$, so $a$ cannot be an $a_{i}$. Let $Z = \set{a, a_{6}, a_{0}, a_{1}}$.
	
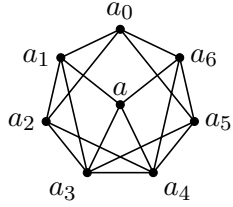
\begin{figure}[H]
		\centering
		\begin{tikzpicture}
			\foreach \pt in {0,1,...,6} 
			{
				\tkzDefPoint(\pt*360/7 + 90:1){v_\pt}
			} 
			\tkzDefPoint(0,0){u}
			\tkzDrawPolySeg(v_0,v_1,v_2,v_3,v_4,v_5,v_6,v_0)
			\tkzDrawPolySeg(v_1,v_3,v_5,v_0,v_2,v_4,v_6)
			\tkzDrawPoints(v_0,v_...,v_6,u)
			\foreach \pt in {1,3,4,6}
			{
				\tkzDrawSegment(u,v_\pt)
			}
			\tkzLabelPoint[above](v_0){$a_{0}$}
			\tkzLabelPoint[left](v_1){$a_{1}$}
			\tkzLabelPoint[left](v_2){$a_{2}$}
			\tkzLabelPoint[below left](v_3){$a_{3}$}
			\tkzLabelPoint[below right](v_4){$a_{4}$}
			\tkzLabelPoint[right](v_5){$a_{5}$}
			\tkzLabelPoint[right](v_6){$a_{6}$}
			\tkzLabelPoint[above](u){$a$}
		\end{tikzpicture}
	\end{figure}
	We claim that each vertex has at most two neighbours in $Z$. By \cref{claim':601}, no vertex is adjacent to all of $a_{6}, a_{0}, a_{1}$. If a vertex $v$ is adjacent to all of $a, a_{0}, a_{1}$, then $va_{0}a_{2}a_{3}a$ is an odd circuit in $G_{a_{1}}$ while if $v$ is adjacent to all of $a, a_{6}, a_{0}$, then $va_{0}a_{5}a_{4}a$ is an odd circuit in $G_{a_{6}}$. Finally if $v$ is adjacent to all of $a, a_{6}, a_{1}$, then $va_{1}a_{3}a_{4}a_{6}$ is an odd circuit in $G_{a}$. Thus
	\begin{equation*}
		4 \delta(G) \leqslant e(Z, G) \leqslant 2 \abs{G},
	\end{equation*}
	which contradicts $\delta(G) > 6/11 \cdot \abs{G}$.
\end{Proof}

We can now show that any vertex with four neighbours in a copy of $H_{2}$ `looks like' one of the vertices of the $H_{2}$.

\begin{claim}\label{claim':02345}
	Let $a_{0}a_{1} \dotsc a_{6}$ be a copy of $H_{2}$ in $G$. Suppose a vertex $v$ is adjacent to at least four of the $a_{i}$. Then there is $i \in \set{0, 2, 3, 4, 5}$ such that
	\begin{equation*}
		\Gamma(v) \cap \set{a_{0}, a_{1}, \dotsc, a_{6}} = \Gamma(a_{i}) \cap \set{a_{0}, a_{1}, \dotsc, a_{6}}.
	\end{equation*}
\end{claim}

\begin{Proof}
	Fix a vertex $v$ which is adjacent to at least four of the $a_{i}$. Firstly there is no $i$ with $v$ adjacent to all of $a_{i - 1}, a_{i}, a_{i + 1}$ else there is a $K_{4}$, or $v$ is adjacent to all of $a_{6}, a_{0}, a_{1}$ contradicting \cref{claim':601}. Now suppose there is an $i$ with $v$ adjacent to all of $a_{i - 2}, a_{i}, a_{i + 2}$. If $i = 2, 3, 4, 5$, then $va_{i - 2}a_{i - 1}a_{i + 1}a_{i + 2}$ is an odd circuit in $G_{a_{i}}$. If $i$ is 1 or 6, then, by symmetry, we may assume $i = 1$: vertex $v$ is adjacent to all of $a_{6}, a_{1}, a_{3}$. Now $v$ is adjacent to none of $a_{0}$ (by \cref{claim':601}), $a_{2}$ (else there is a $K_{4}$) or $a_{4}$ (by \cref{claim':1346}). Thus $v$ is adjacent to $a_{5}$. But then $v$ is adjacent to $a_{1}, a_{3}, a_{5}$ which is the already discounted case of $i = 3$. Finally if $i = 0$, then $v$ is adjacent to all of $a_{5}, a_{0}, a_{2}$ so $v$ is adjacent to neither $a_{1}$ nor $a_{6}$ (else there is a $K_{4}$). By symmetry, we may assume that $v$ is adjacent to $a_{3}$. But then $v$ is adjacent to all of $a_{3}, a_{5}, a_{0}$ which is the already discounted case of $i = 5$.
	
	Hence there is an $i$ with $v$ adjacent to all of $a_{i - 2}, a_{i - 1}, a_{i + 1}, a_{i + 2}$ and no other $a_{j}$. Now $i$ is not 1 as otherwise $a_{0}va_{2}a_{3}a_{4}a_{5}a_{6}$ is a copy of $\overline{C}_{7}$. Similarly $i$ is not 6. For all other $i$, $\Gamma(v) \cap \set{a_{0}, a_{1}, \dotsc, a_{6}} = \Gamma(a_{i}) \cap \set{a_{0}, a_{1}, \dotsc, a_{6}}$.
\end{Proof}

Fix some copy, $v_{0}v_{1} \dotsc v_{6}$, of $H_{2}$ in $G$. We are ready to build some structure around this copy of $H_{2}$. Let
\begin{equation*}
	\begin{aligned}
		D_{i} & = \Gamma(v_{i - 2}, v_{i - 1}, v_{i + 1}, v_{i + 2})  & \text{for } i & = 0, 2, 3, 4, 5, \\
		D_{1} & = \Gamma(v_{0}, v_{2}, v_{3}), & D_{6} & = \Gamma(v_{4}, v_{5}, v_{0}), \\
		D & = \cup_{i = 0}^{6} D_{i}, & R & = V(G) \setminus D.
	\end{aligned}
\end{equation*}
By \cref{claim':02345}, no vertex is adjacent to five of the $v_{i}$, and so the $D_{i}$ are pairwise disjoint. Also from \cref{claim':02345}, the vertices adjacent to exactly four of the $v_{i}$ are those in
\begin{equation*}
	D^{\ast} \coloneqq D_{0} \cup D_{2} \cup D_{3} \cup D_{4} \cup D_{5},
\end{equation*}
and all other vertices are adjacent to at most three of the $v_{i}$. Thus we can give a simple upper bound on the size of $R \cup D_{1} \cup D_{6} = V(G) \setminus D^{\ast}$.
\begin{align}
	7 \delta(G) \leqslant e(\set{v_{0}, v_{1}, \dotsc, v_{6}}, G) & \leqslant 4 \abs{D^{\ast}} + 3 \abs{R \cup D_{1} \cup D_{6}} = 4 \abs{G} - \abs{R \cup D_{1} \cup D_{6}} \nonumber \\
	\Rightarrow \abs{R \cup D_{1} \cup D_{6}} & \leqslant 4 \abs{G} - 7 \delta(G). \label{eq:R'}
\end{align}
Each $D_{i}$ is independent (if $dd'$ is an edge in $D_{i}$, then at least one of $dd'v_{i - 2}v_{i - 1}$, $dd'v_{i + 1}v_{i + 2}$ is a $K_{4}$). Suppose there is an edge $dd'$ between $D_{i}$ and $D_{i + 3}$. If $i \neq 5, 6$, then $dv_{i + 1}v_{i + 2}d'$ is a $K_{4}$. If $i = 5$, then $d'v_{3}v_{4}v_{6}v_{0}$ is a 5-cycle in $G_{d}$ and if $i = 6$, then $dv_{0}v_{1}v_{3}v_{4}$ is a 5-cycle in $G_{d'}$. Hence $D_{i} \cup D_{i + 3}$ is an independent set for all $i$. Finally, if $d_{1} \in D_{1}$ and $d_{6} \in D_{6}$ are adjacent, then $v_{0}d_{1}v_{2}v_{3}v_{4}v_{5}d_{6}$ is a copy of $\overline{C}_{7}$. Thus there is a homomorphism $G[D] \to H_{2}$. Our aim is to get a handle on $R$.

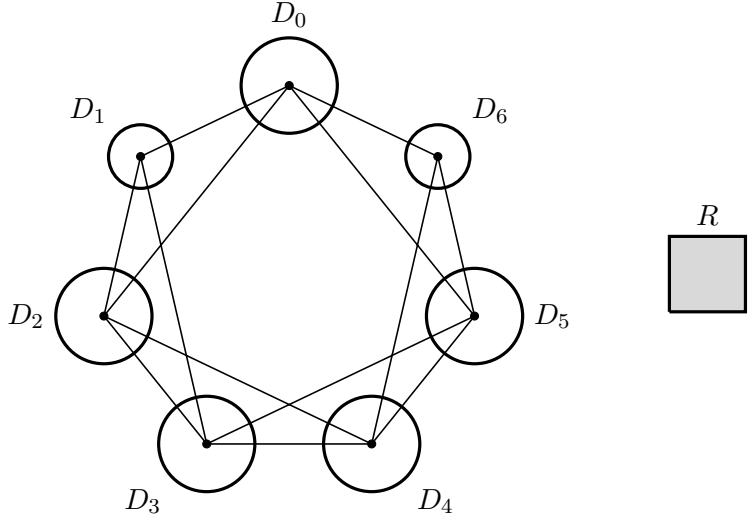
\begin{figure}[H]
	\centering
	\begin{tikzpicture}
		\foreach \pt in {0,2,3,4,5} 
		{
			\tkzDefPoint(\pt*360/7 + 90:2.5){v_\pt}
			\tkzDefShiftPoint[v_\pt](0:18pt){u_\pt}
			\tkzDrawCircle[black,very thick](v_\pt,u_\pt)
		}
		\foreach \pt in {1,6} 
		{
			\tkzDefPoint(\pt*360/7 + 90:2.5){v_\pt}
			\tkzDefShiftPoint[v_\pt](0:12pt){u_\pt}
			\tkzDrawCircle[black,very thick](v_\pt,u_\pt)
		}
		\tkzDrawPolySeg(v_0,v_1,v_2,v_3,v_4,v_5,v_6,v_0) 
		\tkzDrawPolySeg(v_1,v_3,v_5,v_0,v_2,v_4,v_6)
		\tkzDrawPoints(v_0,v_...,v_6)
		\tkzDefShiftPoint[v_0](90:18pt){D0}
		\tkzDefShiftPoint[v_1](135:12pt){D1}
		\tkzDefShiftPoint[v_2](180:18pt){D2}
		\tkzDefShiftPoint[v_3](225:18pt){D3}
		\tkzDefShiftPoint[v_4](315:18pt){D4}
		\tkzDefShiftPoint[v_5](0:18pt){D5}
		\tkzDefShiftPoint[v_6](45:12pt){D6}
		\tkzLabelPoint[above](D0){$D_{0}$}
		\tkzLabelPoint[above left](D1){$D_{1}$}
		\tkzLabelPoint[left](D2){$D_{2}$}
		\tkzLabelPoint[below left](D3){$D_{3}$}
		\tkzLabelPoint[below right](D4){$D_{4}$}
		\tkzLabelPoint[right](D5){$D_{5}$}
		\tkzLabelPoint[above right](D6){$D_{6}$}
		
		\tkzDefPoint(5,-0.5){A}
		\tkzDefPoint(5,0.5){B}
		\tkzDefPoint(6,0.5){C}
		\tkzDefPoint(6,-0.5){D}
		\tkzDrawPolySeg[fill=gray!70,very thick](A,B,C,D,A)
		\tkzDefMidPoint(B,C) \tkzGetPoint{R}
		\tkzLabelPoint[above](R){$R$}
	\end{tikzpicture}
\end{figure}

The following lemma corresponds to \cref{hom:lemma} and is just as useful. Its proof is identical with inequality \eqref{eq:R'} in place of \eqref{eq:R}, $D^{\ast}$ in place of $D$ and $R \cup D_{1} \cup D_{6}$ in place of $R$.

\begin{lemma}\label{hom':lemma}
	Let $X \subset V(G)$ be a set of four vertices. Either there is $x \in R \cup D_{1} \cup D_{6}$ adjacent to all of $X$ or there is $x \in D^{\ast}$ with at least three neighbours in $X$.
\end{lemma}

Our first three claims show that the collections of $v_{i}$ to which vertices can be adjacent are similar to the collections of the $D_{i}$ in which vertices can have neighbours. This will eventually allow us to define $R_{i}$ in a similar way to the proof of \cref{hom2C}.

\begin{claim}\label{claim':consec}
	For all $i$, no vertex has a neighbour in each of $D_{i - 1}$, $D_{i}$, $D_{i + 1}$.
\end{claim}

\begin{Proof}
	If not, choose $d_{i - 1}, d_{i}, d_{i + 1}$ in $D_{i - 1}, D_{i}, D_{i + 1}$ respectively with common neighbour $u$ such that $e(\set{d_{i - 1}, d_{i}}) + e(\set{d_{i}, d_{i + 1}})$ is maximal. Proceeding exactly as in the proof of \cref{hom:3consec} (with \cref{hom':lemma} in place of \cref{hom:lemma}) shows there is some $d \in D^{\ast}$ adjacent to $u, d_{i - 1}$ and $d_{i + 1}$. Then $d \in D_{i}$.
	\begin{itemize}[noitemsep]
		\item If $i = 0$, then $dd_{1}v_{2}v_{3}v_{4}v_{5}d_{6}$ form a copy of $H_{2}$, while $u$ is adjacent to $d_{6}$, $d$ and $d_{1}$, contrary to \cref{claim':601}.
		\item If $i \neq 0$, then $v_{i - 1}$ is adjacent to $d_{i + 1}, v_{i + 1}$ and $v_{i + 1}$ is adjacent to $d_{i - 1}, v_{i - 1}$ so $ud_{i - 1}v_{i + 1}v_{i - 1}d_{i + 1}$ is an odd circuit in $G_{d}$. \qedhere
	\end{itemize}
\end{Proof}

\begin{claim}\label{claim':spaced}
	For all $i \neq 0$, no vertex has a neighbour in each of $D_{i - 2}$, $D_{i}$, $D_{i + 2}$.
\end{claim}

\begin{Proof}
	If not, let $i \neq 0$ and choose $d_{i - 2}, d_{i}, d_{i + 2}$ in $D_{i - 2}, D_{i}, D_{i + 2}$ respectively with common neighbour $u$ such that $e(\set{d_{i - 2}, d_{i}}) + e(\set{d_{i}, d_{i + 2}})$ is maximal. 
	
	Proceeding exactly as in the proof of \cref{hom:3spaced} (with \cref{hom':lemma} for \cref{hom:lemma} and \cref{claim':consec} for \cref{hom:3consec}) shows that $u \in R$ and there is some $d \in D^{\ast}$ adjacent to $u, d_{i - 2}$ and $d_{i + 2}$. In particular, $d \in D_{i} \cup D_{i + 3} \cup D_{i - 3}$. Suppose $d \in D_{i}$: $i \neq 0$ so $v_{i - 1}, v_{i + 1}$ are adjacent so $ud_{i - 2}v_{i - 1}v_{i + 1}d_{i + 2}$ is an odd circuit in $G_{d}$. Thus $d \in D_{i - 3} \cup D_{i + 3}$. By symmetry we may assume that $d \in D_{i + 3}$. Write $d_{i + 3}$ for $d$.
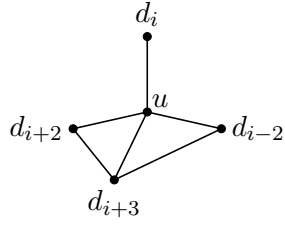
	\begin{figure}[H]
		\centering
		\begin{tikzpicture}
			\foreach \pt in {0,2,3,5} 
			{
				\tkzDefPoint(\pt*360/7 + 90:1){v_\pt}
			} 
			\tkzDefPoint(0,0){u}
			\tkzDrawSegments(u,v_0 u,v_2 u,v_3 u,v_5 v_2,v_3 v_3,v_5)
			\tkzDrawPoints(u,v_0,v_2,v_3,v_5)
			\tkzLabelPoint[above right = -1pt](u){$u$}
			\tkzLabelPoint[right](v_5){$d_{i - 2}$}
			\tkzLabelPoint[above](v_0){$d_{i}$}
			\tkzLabelPoint[left](v_2){$d_{i + 2}$}
			\tkzLabelPoint[below](v_3){$d_{i + 3}$}
		\end{tikzpicture}
	\end{figure}
	\noindent Continuing as in the proof of \cref{hom:3spaced} shows that there is $d'_{i} \in D_{i}$ adjacent to both $u$ and $d_{i - 2}$.
	
	If $i \neq 2$, then $v_{i - 1}$ and $v_{i - 3}$ are adjacent, so $d'_{i}ud_{i + 3}v_{i - 3}v_{i - 1}$ is an odd circuit in $G_{d_{i - 2}}$. Finally if $i = 2$, then $v_{3}d_{4}d_{5}ud_{0}d'_{2}v_{1}$ is a copy of $H_{2}$ in $G$ (note that all the vertices are distinct: vertex $u \in R$ and the others are in distinct $D_{i}$). However, $v_{2}$ is adjacent to all of $v_{1}$, $v_{3}$, $d_{4}$ contrary to \cref{claim':601}.
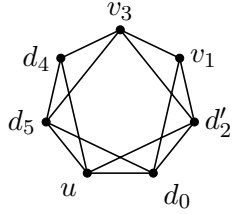
	\begin{figure}[H]
		\centering
		\begin{tikzpicture}
			\foreach \pt in {0,1,...,6} 
			{
				\tkzDefPoint(\pt*360/7 + 90:1){v_\pt}
			} 
			\tkzDrawPolySeg(v_0,v_1,v_2,v_3,v_4,v_5,v_6,v_0)
			\tkzDrawPolySeg(v_1,v_3,v_5,v_0, v_2,v_4,v_6)
			\tkzDrawPoints(v_0,v_...,v_6)
			\tkzLabelPoint[above](v_0){$v_{3}$}
			\tkzLabelPoint[left](v_1){$d_{4}$}
			\tkzLabelPoint[left](v_2){$d_{5}$}
			\tkzLabelPoint[below left](v_3){$u$}
			\tkzLabelPoint[below right](v_4){$d_{0}$}
			\tkzLabelPoint[right](v_5){$d'_{2}$}
			\tkzLabelPoint[right](v_6){$v_{1}$}
		\end{tikzpicture}
	\end{figure}
	\vspace{-12pt}
\end{Proof}

\begin{claim}\label{claim':602}
	No vertex has a neighbour in each of $D_{6}$, $D_{0}$, $D_{2}$ and no vertex has a neighbour in each of $D_{5}$, $D_{0}$, $D_{1}$.
\end{claim}

\begin{Proof}
	Suppose not -- by symmetry we may assume that some vertex has a neighbour in each of $D_{6}, D_{0}, D_{2}$. Choose $d_{6} \in D_{6}, d_{0} \in D_{0}$ and $d_{2} \in D_{2}$ with common neighbour $u$ such that $e(\set{d_{6}, d_{0}}) + e(\set{d_{0}, d_{2}})$ is maximal. We first show that $d_{0}$ is adjacent to both $d_{2}$ and $d_{6}$. Apply \cref{hom':lemma} to $\set{u, d_{6}, d_{0}, d_{2}}$.
	
	Suppose some $x$ is adjacent to all of $u, d_{6}, d_{0}, d_{2}$. Now apply \cref{hom':lemma} to $\set{u, x, d_{0}, d_{2}}$: $uxd_{0}$ and $uxd_{2}$ are triangles so some vertex in $D^{\ast}$ is adjacent to $d_{0}$ and $d_{2}$ and one of $u, x$. But no vertex in $D^{\ast}$ has a neighbour in each of $D_{0}, D_{2}$.
	
	Hence, there is some $d \in D^{\ast}$ adjacent to three of $u, d_{6}, d_{0}, d_{2}$. No vertex in $D^{\ast}$ has a neighbour in each of $D_{0}, D_{2}$ so $d$ is adjacent to $u$, $d_{6}$ and one of $d_{0}$, $d_{2}$. If $d$ is adjacent to $d_{0}$, then $d \in D_{5}$. But then $u$ has a neighbour in each of $D_{5}, D_{6}, D_{0}$ contrary to \cref{claim':consec}. Thus $d$ is adjacent to $d_{2}$ and so $d \in D_{0} \cup D_{4}$. If $d \in D_{4}$, then $u$ has a neighbour in each of $D_{2}, D_{4}, D_{6}$, contradicting \cref{claim':spaced}. Hence $d \in D_{0}$. By maximality at the start we must have $d_{0}$ adjacent to both $d_{2}$ and $d_{6}$.
	
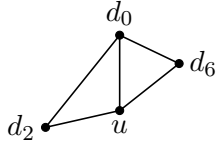
\begin{figure}[H]
		\centering
		\begin{tikzpicture}
			\foreach \pt in {0,2,6} 
			{
				\tkzDefPoint(\pt*360/7 + 90:1){v_\pt}
			} 
			\tkzDefPoint(0,0){u}
			\tkzDrawSegments(u,v_0 u,v_2 u,v_6 v_6,v_0 v_0,v_2)
			\tkzDrawPoints(u,v_0,v_2,v_6)
			\tkzLabelPoint[below](u){$u$}
			\tkzLabelPoint[right](v_6){$d_{6}$}
			\tkzLabelPoint[above](v_0){$d_{0}$}
			\tkzLabelPoint[left](v_2){$d_{2}$}
		\end{tikzpicture}
	\end{figure}
	Now $d_{0}v_{1}d_{2}v_{3}v_{4}v_{5}d_{6}$ is a copy of $H_{2}$ so can be extended, by \cref{claim':H2toH21}, to a copy of $H_{2}^{+}$. That is, there is some other vertex, $v$, adjacent to $v_{5}, d_{0}, d_{2}$. But then $G_{d_{0}}$ contains the odd circuit $ud_{2}vv_{5}d_{6}$.
\end{Proof}

\begin{corollary}
	Every $v$ in $G$ satisfies one of the following properties.
	\begin{itemize}[noitemsep]
		\item $v$ has a neighbour in each of $D_{5}, D_{0}, D_{2}$ and $\Gamma(v) \cap D \subset D_{5} \cup D_{0} \cup D_{2}$.
		\item There is an $i$ such that $\Gamma(v) \cap D \subset \Gamma(v_{i}) \cap D$.
	\end{itemize}
\end{corollary}

\begin{Proof}
	Fix a vertex $v$. First suppose that $v$ has a neighbour in each of $D_{5}, D_{0}, D_{2}$. By \cref{claim':consec}, $v$ has no neighbours in $D_{1} \cup D_{6}$. By \cref{claim':spaced}, $v$ has no neighbours in $D_{3} \cup D_{4}$. Thus $\Gamma(v) \cap D \subset D_{5} \cup D_{0} \cup D_{2}$.
	
	Otherwise $v$ does not have a neighbour in each of $D_{5}, D_{0}, D_{2}$. By \cref{claim':consec,claim':spaced}, there is an $i$ with $\Gamma(v) \cap D \subset D_{i - 2} \cup D_{i - 1} \cup D_{i + 1} \cup D_{i + 2}$. If $i \neq 1, 6$, then $D_{i - 2} \cup D_{i - 1} \cup D_{i + 1} \cup D_{i + 2} = \Gamma(v_{i}) \cap D$ and so we are done. Otherwise we may assume, by symmetry, that $i = 1$: $\Gamma(v) \cap D \subset D_{6} \cup D_{0} \cup D_{2} \cup D_{3}$. By \cref{claim':602}, we have one of the following.
	\begin{itemize}[noitemsep]
		\item $\Gamma(v) \cap D \subset D_{0} \cup D_{2} \cup D_{3} = \Gamma(v_{1}) \cap D$,
		\item $\Gamma(v) \cap D \subset D_{6} \cup D_{2} \cup D_{3} \subset \Gamma(v_{4}) \cap D$,
		\item $\Gamma(v) \cap D \subset D_{6} \cup D_{0} \cup D_{3} \subset \Gamma(v_{5}) \cap D$. \qedhere
	\end{itemize}
\end{Proof}

This corollary gives structure to $R$. Firstly let
\begin{equation*}
	R_{502} = \set{v \in R \colon v \text{ has a neighbour in each of } D_{5}, D_{0}, D_{2}}.
\end{equation*}
Then, for $i = 0, 1, \dotsc, 6$ choose
\begin{equation*}
	R_{i} \subset \set{v \in R \colon \Gamma(v) \cap D \subset \Gamma(v_{i}) \cup D},
\end{equation*}
so that $R_{502} \cup R_{0} \cup R_{1} \cup \dotsb \cup R_{6}$ is a partition of $R$. There may be some flexibility in the choice of the $R_{i}$ -- we will make use of this later. For now we just take any arbitrary choice. Note, for each $i$, that
\begin{equation*}
	e(R_{i}, D_{i} \cup D_{i - 3} \cup D_{i + 3}) = 0,
\end{equation*}
and also that
\begin{equation*}
	e(R_{1}, D_{6}) = e(R_{6}, D_{1}) = e(R_{502}, D_{1} \cup D_{3} \cup D_{4} \cup D_{6}) = 0.
\end{equation*}
For $i = 0, 1, \dotsc, 6$, let $T_{i} = D_{i} \cup R_{i}$. We may give a lower bound for the size of $T_{i}$. Firstly,
\begin{align*}
	d(v_{1}) + d(v_{6}) & = \abs{D_{0}} + \abs{D} - \abs{D_{1}} - \abs{D_{6}} + \abs{R \cap (\Gamma(v_{1}) \cup \Gamma(v_{6}))} + \abs{R \cap \Gamma(v_{1}, v_{6})} \\
	& \leqslant \abs{D_{0}} + \abs{D} - \abs{D_{1}} - \abs{D_{6}} + \abs{R} - \abs{R_{1}} - \abs{R_{6}} - \abs{R_{502}} + \abs{R_{0}},
\end{align*}
so
\begin{equation}\label{eq':DR0}
	\abs{T_{0}} \geqslant 2\delta(G) - \abs{G} + \abs{T_{1}} + \abs{T_{6}} + \abs{R_{502}}.
\end{equation}
Next,
\begin{align*}
	d(v_{0}) + d(v_{2}) & = \abs{D_{1}} + \abs{D} + \abs{R \cap (\Gamma(v_{0}) \cup \Gamma(v_{2}))} + \abs{R \cap \Gamma(v_{0}, v_{2})} \\
	& \leqslant \abs{D_{1}} + \abs{D} + \abs{R} + \abs{R_{1}} + \abs{R_{502}},
\end{align*}
and so (using symmetry for the second inequality)
\begin{equation}\label{eq':DR1}
	\begin{aligned}
		\abs{T_{1}} + \abs{R_{502}} & \geqslant 2 \delta(G) - \abs{G}, \\
		\abs{T_{6}} + \abs{R_{502}} & \geqslant 2 \delta(G) - \abs{G}.
	\end{aligned}
\end{equation}
A similar argument applied to $d(v_{1}) + d(v_{3})$ gives
\begin{equation}\label{eq':DR2}
	\begin{aligned}
		\abs{T_{2}} & \geqslant 2 \delta(G) - \abs{G} + \abs{T_{6}} + \abs{R_{502}}, \\
		\abs{T_{5}} & \geqslant 2 \delta(G) - \abs{G} + \abs{T_{1}} + \abs{R_{502}},
	\end{aligned}
\end{equation}
and one applied to $d(v_{2}) + d(v_{4})$ gives
\begin{equation}\label{eq':DR3}
	\begin{aligned}
		\abs{T_{3}} & \geqslant 2 \delta(G) - \abs{G}, \\
		\abs{T_{4}} & \geqslant 2 \delta(G) - \abs{G}.
	\end{aligned}
\end{equation}
In particular, for all $i \neq 1, 6$,
\begin{equation}\label{eq':DR}
	\abs{T_{i}} \geqslant 2 \delta(G) - \abs{G},
\end{equation}
and, using inequality \eqref{eq':DR1} combined with \eqref{eq':DR0} and \eqref{eq':DR2}, we have for $i = 0,2,5$,
\begin{equation}\label{eq':DR025}
	\abs{T_{i}} \geqslant 4 \delta(G) - 2 \abs{G}.
\end{equation}
The next three claims are technical in nature but will speed up what follows.

\begin{claim}\label{claim':DRD}
	Every $d \in D_{i}$ and every $u \in R_{i + 1}$ have a common neighbour in $D^{\ast}$ unless $i = 2$, in which case they have a common neighbour in $D^{\ast} \cup R_{502}$. Similarly, every $d \in D_{i}$ and every $u \in R_{i - 1}$ have a common neighbour in $D^{\ast}$ unless $i = 5$, in which case they have a common neighbour in $D^{\ast} \cup R_{502}$.
\end{claim}

\begin{Proof}
	It is enough to prove the assertion when $d \in D_{i}$ and $u \in R_{i + 1}$, the other assertion following symmetrically. Now $\Gamma(d) \cup \Gamma(u) \subset V(G) \setminus D_{i - 3}$, so
	\begin{equation*}
		\abs{\Gamma(d) \cap \Gamma(u)} = d(d) + d(u) - \abs{\Gamma(d) \cup \Gamma(u)} \geqslant 2 \delta(G) - \abs{G} + \abs{D_{i - 3}}.
	\end{equation*}
	Now, for $i = 0, 1, 3, 5, 6$, inequality \eqref{eq':DR} gives
	\begin{align*}
		\abs{\Gamma(d) \cap \Gamma(u)} & \geqslant 2 \delta(G) - \abs{G} + \abs{D_{i - 3}} \geqslant 4 \delta(G) - 2 \abs{G} - \abs{R_{i - 3}} \\
		& > 4 \abs{G} - 7 \delta(G) - \abs{R_{i - 3}} \geqslant \abs{R \cup D_{1} \cup D_{6}} - \abs{R_{i - 3}} \\
		& \geqslant \abs{\Gamma(d) \cap (R \cup D_{1} \cup D_{6})},
	\end{align*}
	where we used $\delta(G) > 6/11 \cdot \abs{G}$, inequality \eqref{eq:R'} and $e(D_{i}, R_{i - 3}) = 0$ in the third, fourth and fifth inequalities respectively. Hence there is a common neighbour of $d$ and $u$ which is not in $R \cup D_{1} \cup D_{6}$, so is in $D^{\ast}$.
	
	For $i = 4$, inequality \eqref{eq':DR1} gives
	\begin{align*}
		\abs{\Gamma(d) \cap \Gamma(u)} & \geqslant 2 \delta(G) - \abs{G} + \abs{D_{1}} \geqslant 4 \delta(G) - 2 \abs{G} - \abs{R_{1}} - \abs{R_{502}} \\
		& > 4 \abs{G} - 7 \delta(G) - \abs{R_{1}} - \abs{R_{502}} \geqslant \abs{R \cup D_{1} \cup D_{6}} - \abs{R_{1}} - \abs{R_{502}} \\
		& \geqslant \abs{\Gamma(d) \cap (R \cup D_{1} \cup D_{6})},
	\end{align*}
	where we used $\delta(G) > 6/11 \cdot \abs{G}$, inequality \eqref{eq:R'} and $e(D_{4}, R_{1} \cup R_{502}) = 0$ in the third, fourth and fifth inequalities respectively.
	
	For $i = 2$, inequality \eqref{eq':DR1} gives
	\begin{align*}
		\abs{\Gamma(d) \cap \Gamma(u)} & \geqslant 2 \delta(G) - \abs{G} + \abs{D_{6}} \geqslant 4 \delta(G) - 2 \abs{G} - \abs{R_{6}} - \abs{R_{502}} \\
		& > \abs{R \cup D_{1} \cup D_{6}} - \abs{R_{6}} - \abs{R_{502}} \geqslant \abs{\Gamma(d) \cap ((R \setminus R_{502}) \cup D_{1} \cup D_{6})}. \qedhere
	\end{align*}
\end{Proof}

\begin{claim}\label{claim':DR502D}
	Every $d \in D_{0}$ and every $u \in R_{502}$ have a common neighbour in $D_{2} \cup D_{5}$. Every $d \in D_{2} \cup D_{5}$ and every $u \in R_{502}$ have a common neighbour in $D_{0}$.
\end{claim}

\begin{Proof}
	Fix $d \in D_{0}$ and $u \in R_{502}$. Now $\Gamma(d) \cup \Gamma(u) \subset V(G) \setminus D_{3}$, so, as in the previous claim,
	\begin{align*}
		\abs{\Gamma(d) \cap \Gamma(u)} & \geqslant 4 \delta(G) - 2\abs{G} - \abs{R_{3}} > \abs{R \cup D_{1} \cup D_{6}} - \abs{R_{3}} \\
		& \geqslant \abs{\Gamma(d) \cap (R \cup D_{1} \cup D_{6})},
	\end{align*}
	so $d, u$ have a common neighbour in $D^{\ast}$. But $d \in D_{0}$ so this common neighbour must be in $D_{2} \cup D_{5}$.
	
	Suppose, for contradiction, $d \in D_{2}$, $u \in R_{502}$ have no common neighbour in $D_{0}$: then $d, u$ have no common neighbour in $D$. Now $u \in R_{502}$ so $u$ has a neighbour $d_{0} \in D_{0}$, $d_{5} \in D_{5}$. Apply \cref{hom':lemma} to $\set{u, d, d_{0}, d_{5}}$: no vertex in $D^{\ast}$ is adjacent to $d, d_{0}$ or to $d_{0}, d_{5}$ or to $d, u$ so there is some $v$ adjacent to all of $u, d, d_{0}, d_{5}$.
	
	Now apply \cref{hom':lemma} to $\set{u, v, d_{0}, d_{5}}$: as $uvd_{0}$ and $uvd_{5}$ are triangles, some vertex in $D^{\ast}$ is adjacent to both $d_{0}$, $d_{5}$ (and one of $u$, $v$). But no vertex in $D^{\ast}$ has a neighbour in each of $D_{0}$, $D_{5}$.
\end{Proof}

\begin{claim}\label{claim':RRD}
	For each $i \in \set{0, 1, \dotsc, 6, 502}$, every two vertices in $R_{i}$ have a common neighbour in $D^{\ast}$.
\end{claim}

\begin{Proof}
	Before addressing the claim, we show that 
	\begin{equation}\label{eq':RTTT}
		\abs{R_{502}} + \lvert T_{j} \rvert + \lvert T_{j + 3}\rvert + \lvert T_{j - 3}\rvert \geqslant 3 (2\delta(G) - \abs{G})
	\end{equation}
	for all $j \in \{0, 1, \dotsc, 6\}$. Note that at most one of $j, j - 3, j + 3$ can be 1 or 6. If none of $j, j - 3, j + 3$ is 1 or 6, then \eqref{eq':RTTT} follows from inequality \eqref{eq':DR}. If exactly one of them is 1 or 6, then \eqref{eq':RTTT} follows from inequalities \eqref{eq':DR1} and \eqref{eq':DR}.
	
	We first deal with $i \in \set{0, 1, \dotsc, 6}$. If $u, v \in R_{i}$ have no common neighbour in $D^{\ast}$, then $\Gamma(u) \cap D^{\ast}$ and $\Gamma(v) \cap D^{\ast}$ are disjoint subsets of $D_{i - 2} \cup D_{i - 1} \cup D_{i + 1} \cup D_{i + 2}$. Now, by inequality \eqref{eq:R'},
	\begin{equation*}
		\abs{\Gamma(u) \cap D^{\ast}} \geqslant d(u) - \abs{R \cup D_{1} \cup D_{6}} \geqslant 8 \delta(G) - 4 \abs{G},
	\end{equation*}
	so $\abs{D_{i - 2} \cup D_{i - 1} \cup D_{i + 1} \cup D_{i + 2}} \geqslant 16 \delta(G) - 8 \abs{G}$. This, together with inequality \eqref{eq':RTTT}, gives
	\begin{align*}
		\abs{G} & \geqslant \abs{D_{i - 2} \cup D_{i - 1} \cup D_{i + 1} \cup D_{i + 2}} + \abs{R_{502}} + \abs{T_{i}} + \abs{T_{i + 3}} + \abs{T_{i - 3}} \\
		& \geqslant (16 \delta(G) - 8 \abs{G}) + 3(2 \delta(G) - \abs{G}) = 22 \delta(G) - 11 \abs{G},
	\end{align*}
	which contradicts $\delta(G) > 6/11 \cdot \abs{G}$.
	
	Now we deal with $i = 502$. If $u, v \in R_{502}$ have no common neighbour in $D^{\ast}$, then, as above,
	\begin{equation*}
		\abs{D_{5} \cup D_{0} \cup D_{2}} \geqslant 16 \delta(G) - 8 \abs{G}.
	\end{equation*}
	But then combining this with inequality~\eqref{eq':RTTT} for $j = 1$ gives
	\begin{align*}
		\abs{G} & \geqslant \abs{D_{5} \cup D_{0} \cup D_{2}} + \abs{R_{502}} + \abs{T_{1}} + \abs{T_{4}} + \abs{T_{3}} \\
		& \geqslant 16 \delta(G) - 8 \abs{G} + 3(2 \delta(G) - \abs{G}) = 22 \delta(G) - 11 \abs{G},
	\end{align*}
	which contradicts $\delta(G) > 6/11 \cdot \abs{G}$.
\end{Proof}

We now aim to show that $T_{i} = D_{i} \cup R_{i}$ is independent for all $i$.

\begin{claim}\label{claim':DTindep}
	For all $i$ and $d \in D_{i}$\textnormal{:} $\Gamma(d) \cap T_{i - 1}$ and $\Gamma(d) \cap T_{i + 1}$ are independent.
\end{claim}

\begin{Proof}
	Suppose there is $d_i \in D_{i}$ such that $\Gamma(d_i) \cap T_{i + 1}$ contains the edge $uv$. As $D_{i + 1}$ is independent and $e(D_{i + 1}, R_{i + 1}) = 0$, both $u, v \in R_{i + 1}$.
	
	We first deal with the case when $i$ is not 2. Then, by \cref{claim':DRD}, there is $d_{u} \in D^{\ast}$ adjacent to both $u, d_{i}$. As $d_{u}$ is adjacent to $u \in R_{i + 1}$, $d_{u}$ is adjacent to $v_{i + 1}$. Similarly there is $d_{v} \in D^{\ast}$ adjacent to $v, d_{i}, v_{i + 1}$. But then $v_{i + 1}d_{u}uvd_{v}$ is an odd circuit in $G_{d_{i}}$.
	
	Now suppose $i = 2$ and write $d = d_{2}$. By \cref{claim':DRD}, $d_{2}$ and $u$ have a common neighbour $x_{u} \in D^{\ast} \cup R_{502}$. If $x_{u} \in D^{\ast}$, then as $x_{u}$ is adjacent to $u \in R_{3}$, $x_{u}$ is adjacent to $v_{3}$, while if $x_{u} \in R_{502}$, then, by \cref{claim':DR502D}, $x_{u}$ and $d_{2}$ have a common neighbour $d' \in D_{0}$. Taking $d_{u} = v_{3}$ in the former case and $d_{u} = d'$ in the latter, we see that $x_{u}$ and $d_{2}$ have a common neighbour $d_{u} \in D$ which is adjacent to $v_{1}$. Similarly $d_{2}$ and $v$ have a common neighbour $x_{v}$ such that $x_{v}$ and $d_{2}$ have a common neighbour $d_{v}$ which is adjacent to $v_{1}$. But then $v_{1}d_{u}x_{u}uvx_{v}d_{v}$ is an odd circuit in $G_{d_{2}}$.
\end{Proof}

\begin{claim}\label{claim':Tindep}
	For all $i$, $T_{i}$ is independent.
\end{claim}

\begin{Proof}
	Suppose $uv$ is an edge in $T_{i}$ -- as $D_{i}$ is independent and $e(D_{i}, R_{i}) = 0$ we have $u,v \in R_{i}$. By \cref{claim':RRD}, $u, v$ have a common neighbour in $d \in D^{\ast}$. By \cref{claim':DTindep}, $d \not \in D_{i - 1} \cup D_{i + 1}$ so $d \in D_{i - 2} \cup D_{i + 2}$. We complete the argument assuming that $d \in D_{i - 2}$. The other case is analogous. Write $d = d_{i - 2}$.
	
	Apply \cref{hom':lemma} to $\set{u, v, d_{i - 2}, v_{i + 2}}$: as $uvd_{i - 2}$ is a triangle, there is $d' \in D^{\ast}$ adjacent to $v_{i + 2}$ and to two of $u, v, d_{i - 2}$. In particular, $d'$ is adjacent to $v_{i + 2}$ and to at least one of $u, v \in R_{i}$, so $d' \in D_{i + 1}$. But then $d'$ cannot be adjacent to $d_{i - 2}$ and so is adjacent to both $u$ and $v$. However, edge $uv$ lies in $\Gamma(d') \cap R_{i}$, contradicting \cref{claim':DTindep}.
\end{Proof}

\begin{claim}\label{claim':R502indep}
	$R_{502}$ is independent.
\end{claim}

\begin{Proof}
	Suppose $uv$ is an edge in $R_{502}$. By \cref{claim':RRD}, $u, v$ have a common neighbour $d \in D_{5} \cup D_{0} \cup D_{2}$.
	
	First suppose that $d \in D_{0}$. By \cref{claim':DR502D}, $d, u$ have a common neighbour $d_{u} \in D_{2} \cup D_{5}$ and $d, v$ have a common neighbour $d_{v} \in D_{2} \cup D_{5}$. We may assume that $d_{u} \in D_{2}$.
	\begin{itemize}[noitemsep]
		\item If $d_{v} \in D_{2}$, then $v_{1}d_{u}uvd_{v}$ is an odd circuit in $G_{d}$.
		\item If $d_{v} \in D_{5}$, then $dv_{1}d_{u}v_{3}v_{4}d_{v}v_{6}$ is a copy of $H_{2}$ so, by \cref{claim':H2toH21}, can be extended to a copy of $H_{2}^{+}$: there is $x$ adjacent to all of $d, d_{u}, d_{v}$. But then $xd_{u}uvd_{v}$ is an odd circuit in $G_{d}$.
	\end{itemize}
	Next suppose that $d \in D_{2} \cup D_{5}$. By symmetry, we may assume that $d \in D_{2}$. By \cref{claim':DR502D}, $d, u$ have a common neighbour $d_{u} \in D_{0}$ and $d, v$ have a common neighbour $d_{v} \in D_{0}$. But then $v_{1}d_{u}uvd_{v}$ is an odd circuit in $G_{d}$.
\end{Proof}

We have made good progress: we now know that there is a homomorphism $G \to K_{8}$ and so $G$ is 8-colourable.

\begin{claim}\label{claim':R502R1R6}
	$e(R_{502}, T_{1} \cup T_{6}) = 0$.
\end{claim}

\begin{Proof}
	If not, then we may assume there is an edge $uv$ with $u \in R_{502}$ and $v \in T_{1}$. As $e(R_{502}, D_{1}) = 0$, we have $v  \in R_1$. We first show that $u, v$ have a common neighbour in $D_{0} \cup D_{2}$. Apply \cref{hom':lemma} to $\set{u, v, v_{0}, v_{2}}$: any common neighbour of $v_{0}, v_{2}$ is in $T_{1} \cup R_{502}$ so is adjacent to at most one of $u, v$. Hence there is $d \in D^{\ast}$ adjacent to three of $u, v, v_{0}, v_{2}$. No vertex of $D^{\ast}$ is adjacent to both $v_{0}, v_{2}$ so $d$ is a common neighbour of $u, v$. As $d$ is adjacent to $u$ and $v$, $d \in D_{0} \cup D_{2}$.
	
	First suppose $d \in D_{0}$. By \cref{claim':DRD}, $d, v$ have a common neighbour $d_{v} \in D^{\ast}$: $d_{v}$ is adjacent to both $d, v$ so $d_{v} \in D_{2}$. By \cref{claim':DR502D}, $d, u$ have a common neighbour $d_{u} \in D_{2} \cup D_{5}$.
	\begin{itemize}[noitemsep]
		\item If $d_{u} \in D_{2}$, then $v_{1} d_{u} u v d_{v}$ is an odd circuit in $G_{d}$.
		\item If $d_{u} \in D_{5}$, then $dv_{1}d_{v}v_{3}v_{4}d_{u}v_{6}$ is a copy of $H_{2}$ so, by \cref{claim':H2toH21}, there is a vertex $x$ adjacent to all of $d, d_{u}, d_{v}$. But then $xd_{u}uvd_{v}$ is an odd circuit in $G_{d}$.
	\end{itemize}
	
	Now suppose $d \in D_{2}$. By \cref{claim':DR502D}, $d, u$ have a common neighbour $d_{u} \in D_{0}$. By \cref{claim':DRD}, $d, v$ have a common neighbour $d_{v} \in D^{\ast}$: $d_{v}$ is adjacent to both $d, v$ so $d_{v} \in D_{0} \cup D_{3}$. But then $v_{1}d_{u}uvd_{v}$ is an odd circuit in $G_{d}$.
\end{Proof}

\begin{claim}\label{claim':R1R6}
	$e(T_{1}, T_{6}) = 0$.
\end{claim}

\begin{Proof}
	If not, then there is an edge $uv$ with $u \in T_{1}$, $v \in T_{6}$. As $e(T_1, D_6) = e(T_6, D_1)$, we have $u \in R_{1}$ and $v \in R_{6}$. We first show that $u, v$ have a common neighbour $d \in D_{0}$. Apply \cref{hom':lemma} to $\set{v_{0}, v_{2}, u, v}$: any common neighbour of $v_{0}, v_{2}$ is in $T_{1} \cup R_{502}$ so is not adjacent to $u$. Hence there is $d \in D^{\ast}$ adjacent to three of $v_{0}, v_{2}, u, v$. No vertex in $D^{\ast}$ is adjacent to $v_{0}, v_{2}$ so $d$ is adjacent to $u, v$. But $u \in R_{1}, v \in R_{6}$ so $d \in D_{0}$.
	
	By \cref{claim':DRD}, $d, u$ have a common neighbour $d_{u} \in D^{\ast}$ and $d, v$ have a common neighbour $d_{v} \in D^{\ast}$. As $d \in D_{0}$ and $u \in R_{1}$, we have $d_{u} \in D_{2}$. Similarly, $d_{v} \in D_{5}$. Now $dv_{1}d_{u}v_{3}v_{4}d_{v}v_{6}$ is a copy of $H_{2}$, so, by \cref{claim':H2toH21}, there is a vertex $x$ adjacent to all of $d, d_{u}, d_{v}$. But then $xd_{u}uvd_{v}$ is an odd circuit in $G_{d}$.
\end{Proof}

Before proceeding it will help to give structure to $G_{d}$ for each $d \in D_{1} \cup D_{2} \cup \dotsb \cup D_{6}$. This corresponds to \cref{claim:Gd} in the proof of \cref{hom2C}.

\begin{claim}\label{claim':Gd}
	For each $i \in \set{1, 2, \dotsc, 6}$ and every $d \in D_{i}$, $G_{d}$ is connected bipartite. Furthermore, there is a bipartition of $G_{d}$ into two vertex classes $A_{d}$, $B_{d}$ which satisfy $(T_{i - 1} \cup D_{i + 2}) \cap \Gamma(d) \subset A_{d}$, $(T_{i + 1} \cup D_{i - 2}) \cap \Gamma(d) \subset B_{d}$ and at least one of $R_{i + 2} \cap \Gamma(d) \subset A_{d}$, $R_{i - 2} \cap \Gamma(d) \subset B_{d}$ occurs. If $i = 2$, then $R_{502} \cap \Gamma(d) \subset A_{d}$ and if $i = 5$, then $R_{502} \cap \Gamma(d) \subset B_{d}$.
\end{claim}

\begin{Proof}
	Fix $d \in D_{i}$ and define for $j = i - 2, i - 1, i + 1, i + 2$,
	\begin{align*}
		D_{j}^{d} & = D_{j} \cap \Gamma(d), \\
		R_{j}^{d} & = R_{j} \cap \Gamma(d), \\
		R_{502}^{d} & = R_{502} \cap \Gamma(d),
	\end{align*}
	and note that these partition $V(G_{d})$ (and some of them can be empty). Also let $T_{j}^{d} = T_{j} \cap \Gamma(d)$. Vertex $v_{i - 1} \in G_{d}$. We let
	\begin{align*}
		A_{d} & = \set{x \in G_{d} \colon \operatorname{dist}_{G_{d}}(x, v_{i - 1}) \text{ is even}}, \\
		B_{d} & = \set{x \in G_{d} \colon \operatorname{dist}_{G_{d}}(x, v_{i - 1}) \text{ is odd}}.
	\end{align*}
	$G$ is locally bipartite so $G_{d}$ is bipartite and so $A_{d}$ and $B_{d}$ are independent sets. Now, as $i$ is not 0,
	\begin{itemize}[noitemsep]
		\item $v_{i - 1} \in A_{d}$, $v_{i + 1} \in B_{d}$.
		\item $v_{i - 1}$ is adjacent to all of $D_{i - 2}^{d} \cup D_{i + 1}^{d}$, so $D_{i - 2}^{d} \cup D_{i + 1}^{d} \subset B_{d}$.
		\item $v_{i + 1}$ is adjacent to all of $D_{i - 1}^{d} \cup D_{i + 2}^{d}$, so $D_{i - 1}^{d} \cup D_{i + 2}^{d} \subset A_{d}$.
	\end{itemize}
	If $i = 2$, then, by \cref{claim':DR502D}, any $x \in R_{502}^{d}$ has a neighbour in $D_{0}^{d} \subset B_{d}$, so $R_{502}^{d} \subset A_{d}$. If $i = 5$, then, by \cref{claim':DR502D}, any $x \in R_{502}^{d}$ has a neighbour in $D_{0}^{d} \subset A_{d}$ so $R_{502}^{d} \subset B_{d}$. For other $i$, $R_{502}^{d}$ is empty.
	
	We next show that $R_{i - 1}^{d} \subset A_{d}$. Fix $x \in R_{i - 1}^{d}$ -- it suffices to show $x \in A_{d}$. Suppose $x, d$ have a common neighbour in $d' \in D$. As $x \in R_{i - 1}$ and $d \in D_{i}$, $d'$ must be in $D_{i - 2} \cup D_{i + 1}$. Hence $d' \in B_{d}$ and so $x \in A_{d}$. On the other hand if $x, d$ do not have a common neighbour in $D$, then \cref{claim':DRD} guarantees that $i = 5$ and $x$ has a neighbour in $R_{502}^{d} \subset B_{d}$, so $x \in A_{d}$. Similarly $R_{i + 1}^{d} \subset B_{d}$.
	
	We now show that at least one of $R_{i + 2}^{d} \subset A_{d}$, $R_{i - 2}^{d} \subset B_{d}$ occurs. If not, then there is $u \in R_{i + 2}^{d} \setminus A_{d}$ and $v \in R_{i - 2}^{d} \setminus B_{d}$. Focus on $u$: $u \not \in A_{d}$ so $\Gamma_{G_{d}}(u) \subset V(G_{d}) - B_{d} \subset T_{i - 1}^{d} \cup T_{i + 2}^{d} \cup R_{i - 2}^{d} \cup R_{502}^{d}$. But $u \in R_{i + 2}$, the set $T_{i + 2}$ is independent and $e(R_{i + 2}, D_{i - 1}) = 0$, so, in fact,
	\begin{equation*}
		\Gamma_{G_{d}}(u) \subset R_{i - 1}^{d} \cup R_{i - 2}^{d} \cup R_{502}^{d}.
	\end{equation*}
	Similarly
	\begin{equation*}
		\Gamma_{G_{d}}(v) \subset R_{i + 1}^{d} \cup R_{i + 2}^{d} \cup R_{502}^{d}.
	\end{equation*}
	If $u$ and $v$ both have a neighbour in $R_{502}^{d}$, then $R_{502}^{d}$ would be non-empty, so $i = 2, 5$ and either $R_{502}^{d} \subset A_{d}$ or $R_{502}^{d} \subset B_{d}$. The former contradicts $v \not \in B_{d}$ and the latter contradicts $u \not \in A_{d}$. Thus, at most one of $u, v$ has a neighbour in $R_{502}^{d}$. In particular,
	\begin{equation*}
		\abs{\Gamma_{G_{d}}(u)} + \abs{\Gamma_{G_{d}}(v)} \leqslant \lvert R_{i - 2}^{d} \rvert + \lvert R_{i - 1}^{d} \rvert + \lvert R_{i + 1}^{d} \rvert + \lvert R_{i + 2}^{d} \rvert + \lvert R_{502}^{d} \rvert \leqslant \abs{R}.
	\end{equation*}
	But then inequality \eqref{eq:R'} gives
	\begin{equation*}
		4 \abs{G} - 7 \delta(G) \geqslant \abs{R} \geqslant d(d, u) + d(d, v) \geqslant 4 \delta(G) - 2\abs{G},
	\end{equation*}
	which contradicts $\delta(G) > 6/11 \cdot \abs{G}$.
	
	Finally the proof that $G_{d}$ is connected is identical to that part of the proof of \cref{claim:Gd}.
\end{Proof}

To prove that there is a homomorphism $G \to H_{2}^{+}$ we would need to show that $e(R_{i}, R_{i + 3}) = 0$ for all $i$ and $e(R_{3} \cup R_{4}, R_{502}) = 0$. We make a start.

\begin{claim}\label{claim':RR024}
	$e(T_{i}, T_{i + 3}) = 0$ for $i = 0, 2, 4$.
\end{claim}

\begin{Proof}
	Suppose not: there is an edge $uv$ with $u \in T_{i + 3}$, $v \in T_{i}$. As $e(T_i, D_{i} \cup D_{i - 3} \cup D_{i + 3}) = 0$, we must have $u \in R_{i + 3}$ and $v \in R_{i}$. We first show that $u$, $v$ have a common neighbour $d \in D^{\ast} \cap (D_{i + 1} \cup D_{i + 2})$. Apply \cref{hom':lemma} to $\set{u, v, v_{i + 1}, v_{i + 2}}$: any common neighbour of $v_{i + 1}, v_{i + 2}$ is in $T_{i} \cup T_{i + 3}$ so is adjacent to at most one of $u, v$. Moreover, any common neighbour of $v_{i + 1}, v_{i + 2}$ in $D^{\ast}$ is in $D_{i} \cup D_{i + 3}$ and so is adjacent to neither $u$ nor $v$. Hence there is $d \in D^{\ast}$ adjacent to both $u, v$ and one of $v_{i + 1}, v_{i + 2}$ -- in particular, $d \in D_{i + 1} \cup D_{i + 2}$. 
	
	If $i = 2$, we may take $d \in D_{3}$, by symmetry. If $i = 0, 4$ we may assume, by symmetry that $i = 4$. Since $d \in D^{\ast}$, we have $d \in D_{5}$. In conclusion, we have adjacent vertices $u \in R_{i + 3}$, $v \in R_{i}$ with common neighbour $d \in D_{i + 1}$ where $i$ is 2 or 4. Consider the bipartition of $G_{d}$ given by \cref{claim':Gd}:
	\begin{itemize}[noitemsep]
		\item $(T_{i} \cup D_{i + 3}) \cap \Gamma(d) \subset A_{d}$.
		\item $(T_{i + 2} \cup D_{i - 1} \cup R_{502}) \cap \Gamma(d) \subset B_{d}$.
		\item At least one of $R_{i + 3} \cap \Gamma(d) \subset A_{d}$ or $R_{i - 1} \cap \Gamma(d) \subset B_{d}$ occurs.
	\end{itemize}
	As $v \in R_{i}$, we have $v \in A_{d}$ and so $u \in B_{d}$. But $u \in R_{i + 3} \cap \Gamma(d)$, so $R_{i - 1} \cap \Gamma(d) \subset B_{d}$ occurs. Now $\Gamma_{G_{d}}(u) \subset A_{d} \subset (T_{i} \cup T_{i + 3}) \cap \Gamma(d)$. But $u \in R_{i + 3}$, the set $T_{i + 3}$ is independent and $e(R_{i + 3}, D_{i}) = 0$, so
	\begin{equation*}
		\Gamma_{G_{d}}(u) \subset R_{i} \cap \Gamma(d).
	\end{equation*}
	Thus, $\abs{R_{i}} \geqslant d(d, u) \geqslant 2 \delta(G) - \abs{G}$. But then, using inequalities \eqref{eq:R'} and \eqref{eq':DR1},
	\begin{equation*}
		4 \abs{G} - 7 \delta(G) \geqslant \abs{R \cup D_{1} \cup D_{6}} \geqslant \abs{R_{i}} + \abs{T_{1} \cup R_{502}} \geqslant 4 \delta(G) - 2 \abs{G},
	\end{equation*}
	which contradicts $\delta(G) > 6/11 \cdot \abs{G}$.
\end{Proof}	

We have been flexible about the $R_{i}$ and so all of our results thus far hold for any $R_{i}$ satisfying their definition. Now is the time to make a further choice. We choose the $R_{i}$ so that
\begin{equation}\label{eq':S}
	S = \sum_{i = 0}^{6} e(R_{i}, R_{i + 3}) + e(R_{3} \cup R_{4}, R_{502})
\end{equation}
is minimal.

\begin{claim}\label{claim':RR13}
	$e(T_{i}, T_{i + 3}) = 0$ for $i = 1, 3$.
\end{claim}

\begin{Proof}
	By symmetry it suffices to prove this for $i = 3$. Suppose we have $u \in T_{6}$, $v \in T_{3}$ with $u$ adjacent to $v$. As $e(D_{3}, T_{6}) = e(D_{6}, T_{3}) = 0$, we must have $u \in R_{6}$ and $v \in R_{3}$. Just as in the proof of \cref{claim':RR024}, $u$ and $v$ have a common neighbour $d \in D_{4} \cup D_{5}$. When $d \in D_{4}$, the argument of \cref{claim':RR024} works again (with $i = 3$). We deal with the more difficult $d \in D_{5}$ case. For any $d \in D_{5} \cap \Gamma(u, v)$, consider the bipartition of $G_{d}$ given by \cref{claim':Gd}:
	\begin{itemize}[noitemsep]
		\item $(T_{4} \cup D_{0}) \cap \Gamma(d) \subset A_{d}$.
		\item $(T_{6} \cup D_{3} \cup R_{502}) \cap \Gamma(d) \subset B_{d}$.
		\item At least one of $R_{0} \cap \Gamma(d) \subset A_{d}$ or $R_{3} \cap \Gamma(d) \subset B_{d}$ occurs.
	\end{itemize}
	As $u \in R_{6}$, $u \in B_{d}$ and so $v \in A_{d}$. But $v \in R_{3} \cap \Gamma(d)$, so $R_{0} \cap \Gamma(d) \subset A_{d}$ occurs. Now $\Gamma_{G_{d}}(v) \subset B_{d} \subset T_{6} \cup T_{3} \cup R_{502}$. But $v \in R_{3}$, the set $T_{3}$ is independent and $e(R_{3}, D_{6}) = 0$, so
	\begin{equation*}
		\Gamma_{G_{d}}(v) \subset (R_{6} \cup R_{502}) \cap \Gamma(d).
	\end{equation*}
	Note that this holds for any choice of $d \in D_{5} \cap \Gamma(u, v)$.
	
	We first deal with the case where $v$ has at least one neighbour in $D_{4}$. Pick any $d_{5} \in \Gamma(u, v) \cap D_{5}$, $d_{4} \in \Gamma(v) \cap D_{4}$.
	
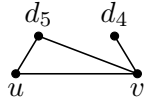
\begin{figure}[H]
		\centering
		\begin{tikzpicture}
			\tkzDefPoint(-0.8,0){u}
			\tkzDefPoint(0.8,0){v}
			\tkzDefPoint(-0.5,0.5){d1}
			\tkzDefPoint(0.5,0.5){d0}
			\tkzDrawPolySeg(u,v,d1,u)
			\tkzDrawSegment(v,d0)
			\tkzDrawPoints(u,v,d0,d1)
			\tkzLabelPoint[below](u){$u$}
			\tkzLabelPoint[below](v){$v$}
			\tkzLabelPoint[above](d0){$d_{4}$}
			\tkzLabelPoint[above](d1){$d_{5}$}
		\end{tikzpicture}
	\end{figure}	
	Apply \cref{hom':lemma} to $\set{u,v,d_{4},d_{5}}$: the vertices $d_{5}, u, v$ form a triangle so some $d \in D^{\ast}$ is adjacent to $d_{4}$ and to two of $d_{5}, u, v$. If $d$ is adjacent to $d_{5}$, then $d \in D_{3} \cup D_{6}$, so $d$ is adjacent to neither $u$ nor $v$. Hence $d \in \Gamma(u, v, d_{4}) \cap D^{\ast}$, so $d \in D_{5}$. But then $d \in \Gamma(u, v) \cap D_{5}$ and $\Gamma_{G_{d}}(v)$ contains $d_{4} \not \in (R_{6} \cup R_{502}) \cap \Gamma(d)$, a contradiction.
	
	We are finally left with the case where $v$ has no neighbours in $D_{4}$: when we were choosing the $R_{i}$ we could have put $u$ in $R_{0}$. In particular, \cref{claim':RR024} gives $e(u, R_{4}) = 0$. Also $R_{3}$ is independent, so $e(u, R_{3}) = 0$. Thus if we put $u \in R_{0}$, then $u$ would contribute 0 to $S$ while currently it contributes at least 1 (the edge $uv$ contributes to $e(R_{6}, R_{3})$). This contradicts the minimality of $S$.
\end{Proof}

\begin{Proof}[of \cref{hom2H}]
	Let $G$ be a locally bipartite graph with $\delta(G) > 6/11 \cdot \abs{G}$. By \cref{main4localbip,lemma4H}, $G$ is either 3-colourable, contains $\overline{C}_{7}$ or contains $H_{2}^{+}$. In the first two cases we are done (using \cref{hom2C} -- note that $\overline{C}_{7}$ is 4-colourable). Hence, we may assume that $G$ does not contain a copy of $\overline{C}_{7}$ but does contain a copy of $H_{2}^{+}$. We can thus follow the argument of this subsection, defining the $D_{i}$ and $R_{i}$ and establishing all the claims concerning them.
	
	Note that the $T_{i}$ and $R_{502}$ together partition $V(G)$. By \cref{claim':Tindep}, each $T_{i}$ is independent and, by \cref{claim':R502indep}, $R_{502}$ is independent. By \cref{claim':RR13,claim':RR024}, $e(T_i, T_{i + 3}) = 0$ for $i = 0, 1, \dotsc, 4$ and, by \cref{claim':R1R6}, $e(T_1, T_6) = 0$. Finally, $e(R_{502}, T_1 \cup T_6) = 0$ by \cref{claim':R502R1R6}.
	
	Hence, identifying $T_{i}$ with $v_{i}$ and $R_{502}$ with a single vertex gives a homomorphism from $G$ to the following graph which is $H_{2}^{+}$ with four extra edges. This graph is 4-colourable (colouring shown in the diagram) and so $\chi(G) \leqslant 4$.
	
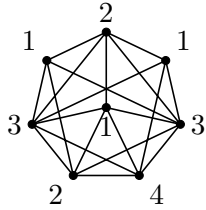
\vspace{-6pt}
	\begin{figure}[H]
		\centering
		\begin{tikzpicture}
			\foreach \pt in {0,1,...,6} 
			{
				\tkzDefPoint(\pt*360/7 + 90:1){v_\pt}
			} 
			\tkzDefPoint(0,0){u}
			\tkzDrawPolySeg(v_0,v_1,v_2,v_3,v_4,v_5,v_6,v_0)
			\tkzDrawPolySeg(v_1,v_3,v_5,v_0, v_2,v_4,v_6)
			\tkzDrawSegments(u,v_0 u,v_2 u,v_5 u,v_3 u,v_4 v_2,v_6 v_1,v_5)
			\tkzDrawPoints(v_0,v_...,v_6)
			\tkzDrawPoint(u)
			\tkzLabelPoint[below](u){1}
			\tkzLabelPoint[above](v_0){2}
			\tkzLabelPoint[above left](v_1){1}
			\tkzLabelPoint[left](v_2){3}
			\tkzLabelPoint[below left](v_3){2}
			\tkzLabelPoint[below right](v_4){4}
			\tkzLabelPoint[right](v_5){3}
			\tkzLabelPoint[above right](v_6){1}
		\end{tikzpicture}
	\end{figure}
	\vspace{-12pt}
\end{Proof}

\subsection{Proof of \texorpdfstring{\cref{hom2Hepsilon}}{Theorem 8}}\label{sec:homtoHepsilon}

In this subsection we prove \cref{hom2Hepsilon}. We remind the reader that there is now a stronger minimum degree condition on $G$: $\delta(G) \geqslant (5/9 - \varepsilon) \cdot \abs{G}$ for some small positive $\varepsilon$. We are also given that $G$ does not contain $\overline{C}_7$. We are in the same position as at the start of the proof of \cref{hom2H} but with a stronger minimum degree condition that we will leverage to give $G$ greater structure. Hence, we may use all of the machinery from our proof of \cref{hom2H} and, in particular, we only need to show that $e(R_{3} \cup R_{4}, R_{502}) = e(R_{1}, R_{5}) = e(R_{2}, R_{6}) = 0$. As before, we choose the $R_{i}$ so that $S$, as given in equation~\eqref{eq':S}, is minimal. Using inequalities \eqref{eq':DR3} and \eqref{eq':DR025}, we have
\begin{align*}
	\abs{T_{3}}, \abs{T_{4}} & \geqslant 2 \delta(G) - \abs{G} > (1/9 - 2 \varepsilon) \abs{G}, \\
	\abs{T_{0}}, \abs{T_{2}}, \abs{T_{5}} & \geqslant 4 \delta(G) - 2 \abs{G} > (2/9 - 4 \varepsilon) \abs{G}.
\end{align*}
Also, by \cref{claim':H2toH21}, we have
\begin{equation*}
	\abs{R_{502}} \geqslant \abs{\Gamma(v_{5}, v_{0}, v_{2})} \geqslant 11 \delta(G) - 6 \abs{G} > (1/9 - 11 \varepsilon) \abs{G}.
\end{equation*}
Now $2(1/9 - 2 \varepsilon) + 3(2/9 - 4 \varepsilon) + (1/9 - 11 \varepsilon) = 1 - 27 \varepsilon$, so in fact we have
\begin{align*}
	\abs{T_{3}}, \abs{T_{4}}, \abs{R_{502}} & = (1/9 - \cO(\epsilon)) \abs{G}, \\
	\abs{T_{0}}, \abs{T_{2}}, \abs{T_{5}} & = (2/9 - \cO(\epsilon)) \abs{G}.
\end{align*}
Throughout we will use $\cO(\epsilon)$ to denote a quantity for which there is an absolute positive constant $C$ (in particular, independent of $G$ and $\varepsilon$) such that the quantity lies between $-C \varepsilon$ and $C \varepsilon$. By inequality \eqref{eq:R'},
\begin{equation*}
	\abs{D_{1} \cup D_{6} \cup R} \leqslant 4 \abs{G} - 7 \delta(G) < (1/9 + 7 \varepsilon) \abs{G}.
\end{equation*}
Putting all this together (and noting that $R_{502} \subset R$) we have
\begin{align}
	\abs{R \setminus R_{502}}, \abs{D_{1}}, \abs{D_{6}} & = \cO(\epsilon) \abs{G}, \nonumber \\
	\abs{D_{3}}, \abs{D_{4}}, \abs{R_{502}} & = (1/9 + \cO(\epsilon)) \abs{G}, \label{eq':5/9} \\
	\abs{D_{0}}, \abs{D_{2}}, \abs{D_{5}} & = (2/9 + \cO(\epsilon)) \abs{G}. \nonumber
\end{align}
Note that these numbers match the weighting of $H_{2}^{+}$ given in \cref{fig:weightings}. That was a weighting of $H_{2}^{+}$ with minimum degree attaining $5/9$.

\begin{claim}\label{claim'':R1R5}
	Provided $\varepsilon > 0$ is sufficiently small, $e(R_{1}, R_{5}) = e(R_{2}, R_{6}) = 0$.
\end{claim}

\begin{Proof}
	Suppose this is false. By symmetry we may take $r_{2} \in R_{2}$ and $r_{6} \in R_{6}$ where $r_{2}r_{6}$ is an edge. We first claim that $r_{2}$ has a neighbour $t_{1} \in T_{1}$. Indeed, if $r_{2}$ does not, then, when we chose the $R_{i}$, we could have put $r_{2}$ in $R_{5}$. Thus, by the minimality of $S$,
	\begin{equation*}
		e(r_{2}, R_{5}) + e(r_{2}, R_{6}) \leqslant e(r_{2}, R_{1}) + e(r_{2}, R_{2}).
	\end{equation*}
	However, the left-hand is positive (the edge $r_{2}r_{6}$ contributes to it), while the right-hand side is zero ($R_{2}$ is independent and $r_{2}$ has no neighbours in $R_{1}$ by assumption). Thus $r_{2}$ does indeed have a neighbour in $t_{1} \in T_{1}$.
	
	Now, $\Gamma(r_{2}) \subset T_{0} \cup T_{1} \cup T_{3} \cup T_{4} \cup R_{502} \cup R_{6}$ and this union has size $(5/9 + \cO(\epsilon)) \abs{G}$, so $r_{2}$ has at most $\cO(\epsilon) \abs{G}$ non-neighbours in $D_{0} \cup D_{3} \cup D_{4}$. Also, $\Gamma(t_{1}) \subset T_{0} \cup T_{2} \cup T_{3} \cup R_{5}$ and this union has size $(5/9 + \cO(\epsilon)) \abs{G}$, so $t_{1}$ has at most $\cO(\epsilon) \abs{G}$ non-neighbours in $D_{0} \cup D_{3}$. Similarly, $r_{6}$ has at most $\cO(\epsilon) \abs{G}$ non-neighbours in $D_{0} \cup D_{4}$. But $D_{0}, D_{3}$ both have size at least $(1/9 + \cO(\epsilon)) \abs{G}$, so, provided $\varepsilon$ is small enough, there is $d_{0} \in D_{0}$ adjacent to all of $r_{2}, t_{1}, r_{6}$ and there is $d_{3} \in D_{3}$ adjacent to both $t_{1}, r_{2}$.
	
	Finally, $\Gamma(d_{3}) \subset T_{1} \cup T_{2} \cup T_{4} \cup T_{5}$ and this union has size $(5/9 + \cO(\epsilon)) \abs{G}$, so $d_{3}$ has at most $\cO(\epsilon) \abs{G}$ non-neighbours in $D_{4}$. But $D_{4}$ has size $(1/9 + \cO(\epsilon)) \abs{G}$, so, provided $\varepsilon$ is small enough, there is $d_{4} \in D_{4}$ adjacent to all of $r_{2}, d_{3}, r_{6}$. But then $d_{0}t_{1}d_{3}d_{4}r_{6}$ is a 5-cycle in $G_{r_{2}}$.
\end{Proof}

\begin{claim}\label{claim'':R3R4R502}
	Provided $\varepsilon > 0$ is sufficiently small, $e(R_{3} \cup R_{4}, R_{502}) = 0$.
\end{claim}

\begin{Proof}
	Suppose this is false. By symmetry we may take some $r_{3} \in R_{3}$ that has at least one neighbour in $R_{502}$. We first claim that $r_{3}$ has at least one neighbour in $D_{4}$. Indeed, if $r_{3}$ does not, then, when we chose the $R_{i}$, we could have put $r_{3}$ in $R_{0}$. Thus, by the minimality of $S$,
	\begin{equation*}
		e(r_{3}, R_{0}) + e(r_{3}, R_{6}) + e(r_{3}, R_{502}) \leqslant e(r_{3}, R_{3}) + e(r_{3}, R_{4}).
	\end{equation*}
	But we showed in \cref{claim':RR024} that $e(R_{0}, R_{4}) = 0$ and we did this before we made the choice to minimise $S$. In particular, as we could have put $r_{3}$ in $R_{0}$ we know that $e(r_{3}, R_{4}) = 0$. Also $R_{3}$ is independent so $e(r_{3}, R_{3}) = 0$. But then, the right-hand side of the inequality is zero, while the left-hand side is positive ($r_{3}$ has at least one neighbour in $R_{502}$).
	
	Thus $r_{3}$ has at least one neighbour in $R_{502}$ and at least one neighbour in $D_{4}$. We may write,
	\begin{equation*}
		\abs{\Gamma(r_{3}) \cap D_{4}} = c_{4} \abs{G}, \quad \abs{\Gamma(r_{3}) \cap R_{502}} = c_{502} \abs{G},
	\end{equation*}
	where $0 < c_{4}, c_{502} \leqslant 1/9 + \cO(\epsilon)$ (using our knowledge of $\abs{D_{4}}$, $\abs{R_{502}}$). Also,
	\begin{align*}
		\abs{\Gamma(r_{3}) \cap D_{5}} & \geqslant d(r_{3}) - \abs{\Gamma(r_{3}) \cap R_{502}} - \abs{\Gamma(r_{3}) \cap D_{4}} - \abs{T_{2}} - \abs{T_{1}} - \abs{R \setminus R_{502}} \\
		& \geqslant (1/3 - c_{4} - c_{502} + \cO(\epsilon)) \abs{G}.
	\end{align*}
	But $1/3 - c_{4} - c_{502} \geqslant 1/9 - \cO(\epsilon)$, so, provided $\varepsilon$ is sufficiently small, $r_{3}$ has at least one neighbour in $D_{5}$.
	
	We next claim that the configuration in \cref{fig:5/9} appears with $d_{5} \in D_{5}$, $d_{4} \in D_{4}$, $r_{502} \in R_{502}$.	
	\begin{figure}[H]
		\centering
		\begin{tikzpicture}
			\foreach \pt in {3,4,5} 
			{
				\tkzDefPoint(\pt*360/7 + 90:1){v_\pt}
			} 
			\tkzDefPoint(0,0){u}
			\tkzDrawPolySeg(v_3,u,v_5,v_4,v_3)
			\tkzDrawSegment(v_3,v_5)
			\tkzDrawPoints(v_3,v_4,v_5,u)
			\tkzLabelPoint[above left](u){$r_{502}$}
			\tkzLabelPoint[right](v_5){$d_{5}$}
			\tkzLabelPoint[below left](v_3){$r_{3}$}
			\tkzLabelPoint[below right](v_4){$d_{4}$}
		\end{tikzpicture}
		\caption{The configuration of \cref{claim'':R3R4R502}.} \label{fig:5/9}
	\end{figure}
	First suppose that $c_{4}, c_{502} \geqslant 1/27$. Pick a neighbour $d_{5} \in D_{5}$ of $r_{3}$. Now $\Gamma(d_{5}) \subset T_{3} \cup T_{4} \cup T_{6} \cup T_{0} \cup R_{502}$ and this union has size $(5/9 + \cO(\epsilon)) \abs{G}$, so $d_{5}$ has at most $\cO(\epsilon) \abs{G}$ non-neighbours in $D_{4} \cup R_{502}$. But $\Gamma(r_{3}) \cap D_{4}$ and $\Gamma(r_{3}) \cap R_{502}$ both have size at least $1/27 \cdot \abs{G}$, so, provided $\varepsilon$ is small enough, $d_{5}$ has a neighbour in each of $\Gamma(r_{3}) \cap D_{4}$, $\Gamma(r_{3}) \cap R_{502}$ giving the configuration in \cref{fig:5/9}.
	
	Otherwise $\min \set{c_{4}, c_{502}} < 1/27$ and so
	\begin{equation*}
		\abs{\Gamma(r_{3}) \cap D_{5}} \geqslant (1/3 - c_{4} - c_{502} + \cO(\epsilon)) \abs{G} \geqslant (1/3 - 1/9 - 1/27 + \cO(\epsilon)) \abs{G} = (5/27 + \cO(\epsilon)) \abs{G}.
	\end{equation*}
	Pick $r_{502} \in R_{502}$ and $d_{4} \in D_{4}$ both adjacent to $r_{3}$. Now, $\Gamma(d_{4}) \subset T_{2} \cup T_{3} \cup T_{5} \cup T_{6}$ and this union has size $(5/9 + \cO(\epsilon)) \abs{G}$, so $d_{4}$ has at most $\cO(\epsilon) \abs{G}$ non-neighbours in $D_{5}$. Also, $\Gamma(r_{502}) \subset T_{0} \cup T_{2} \cup T_{5} \cup R_{3} \cup R_{4}$ and this union has size $(2/3 + \cO(\epsilon)) \abs{G}$, so $r_{502}$ has at most $(1/9 + \cO(\epsilon)) \abs{G}$ non-neighbours in $D_{5}$. But $5/27 > 1/9$, so, provided $\varepsilon$ is sufficiently small, there is some $d_{5} \in \Gamma(r_{3}) \cap D_{5}$ adjacent to both $d_{4}$ and $r_{502}$.
	
	Hence, in all cases, the configuration in \cref{fig:5/9} appears with $d_{5} \in D_{5}, d_{4} \in D_{4}$ and $r_{502} \in R_{502}$. Consider the bipartition of $G_{d_{5}}$ given by \cref{claim':Gd}:
	\begin{itemize}[noitemsep]
		\item $(T_{4} \cup D_{0}) \cap \Gamma(d_{5}) \subset A_{d_{5}}$.
		\item $(T_{6} \cup D_{3} \cup R_{502}) \cap \Gamma(d_{5}) \subset B_{d_{5}}$.
	\end{itemize}
	In particular, $d_{4} \in A_{d_{5}}$ and $r_{502} \in B_{d_{5}}$. But then $r_{3} \in G_{d_{5}}$ has a neighbour in both $A_{d_{5}}$ and $B_{d_{5}}$, so can be in neither, which is a contradiction.
\end{Proof}

\begin{Proof}[of \cref{hom2Hepsilon}]
	Take $\varepsilon > 0$ sufficiently small so that \cref{claim'':R1R5,claim'':R3R4R502} hold. Note that the $T_i$ and $R_{502}$ together partition $V(G)$. By \cref{claim':Tindep}, each $T_i$ is independent and, by \cref{claim':R502indep}, $R_{502}$ is independent. By \cref{claim':RR024,claim':RR13,claim'':R1R5}, $e(T_i, T_{i + 3}) = 0$ for all $i$ and by \cref{claim':R1R6}, $e(T_1, T_6) = 0$. Finally, by \cref{claim':R502R1R6,claim'':R3R4R502}, $e(R_{502}, T_1 \cup T_3 \cup T_4 \cup T_6) = 0$.
	
	Hence, identifying $T_{i}$ with $v_{i}$ and $R_{502}$ with a single vertex gives a homomorphism from $G$ to $H_{2}^{+}$.
\end{Proof}

\section*{Acknowledgements}

The author is grateful to Andrew Thomason for his support and many interesting discussions. I am grateful to the anonymous referees for several helpful comments improving the exposition.

{
\fontsize{11pt}{12pt}
\selectfont
	
\hypersetup{linkcolor={red!70!black}}
\setlength{\parskip}{2pt plus 0.3ex minus 0.3ex}
	
\newcommand{\etalchar}[1]{$^{#1}$}

}
	
\end{document}